%&amstex     
\input amstex\documentstyle {amsppt}  
\pagewidth{12.5 cm}\pageheight{19 cm}\magnification\magstep1
\topmatter
\title Character sheaves on disconnected groups, X\endtitle
\author G. Lusztig\endauthor
\address Department of Mathematics, M.I.T., Cambridge, MA 02139\endaddress
\thanks Supported in part by the National Science Foundation\endthanks
\endtopmatter   
\document
\define\tind{\text{\rm tind}}
\define\tcv{\ti{\cv}}
\define\hPh{\hat\Ph}
\define\bfF{\bar{\fF}}
\define\tfA{\ti{\fA}}
\define\tua{\ti{\un{a}}}
\define\Irr{\text{\rm Irr}}
\define\tWW{\ti{\WW}}

\define\sgn{\text{\rm sgn}}
\define\ubK{\un{\bar K}}
\define\uhD{\un{\hD}}

\define\codim{\text{\rm codim}}

\define\Wb{\WW^\bul}

\define\tir{\ti r}

\define\ua{\un a}

\define\ue{\un e}

\define\uD{\un D}

\define\uK{\un K}

\define\bD{\bar D}

\define\hD{\hat D}

\define\da{\dagger}
\define\dsv{\dashv}

\define\po{\text{\rm pos}}

\define\si{\sim}
\define\wt{\widetilde}
\define\sqc{\sqcup}

\define\qua{\quad}

\define\tcl{\ti\cl}
\define\tca{\ti\ca}

\define\bA{\bar A}

\define\bG{\bar G}
\define\bK{\bar K}
\define\bC{\bar C}
\define\bY{\bar Y}
\define\bX{\bar X}
\define\bZ{\bar Z}

\define\bpi{\bar\p}

\define\lb{\linebreak}

\define\eSb{\endSb}

\define\op{\oplus}
   
\redefine\sp{\spadesuit}
\define\part{\partial}
\define\em{\emptyset}
\define\imp{\implies}
\define\ra{\rangle}
\define\n{\notin}
\define\iy{\infty}
\define\m{\mapsto}
\define\do{\dots}
\define\la{\langle}
\define\bsl{\backslash}

\define\lra{\leftrightarrow}

\define\sub{\subset}    
\define\bxt{\boxtimes}
\define\T{\times}
\define\ti{\tilde}
\define\nl{\newline}
\redefine\i{^{-1}}
\define\fra{\frac}
\define\un{\underline}
\define\ov{\overline}
\define\ot{\otimes}
\define\bbq{\bar{\QQ}_l}

\define\Ad{\text{\rm Ad}}
\define\Hom{\text{\rm Hom}}

\define\Aut{\text{\rm Aut}}
\define\Ind{\text{\rm Ind}}

\define\ind{\text{\rm ind}}

\define\res{\text{\rm res}}

\define\tr{\text{\rm tr}}

\define\supp{\text{\rm supp}}

\define\a{\alpha}
\redefine\b{\beta}
\redefine\c{\chi}
\define\g{\gamma}
\redefine\d{\delta}
\define\e{\epsilon}
\define\et{\eta}
\define\io{\iota}
\redefine\o{\omega}
\define\p{\pi}
\define\ph{\phi}
\define\ps{\psi}
\define\r{\rho}
\define\s{\sigma}
\redefine\t{\tau}
\define\th{\theta}

\redefine\l{\lambda}
\define\z{\zeta}
\define\x{\xi}

\define\vp{\varpi}
\define\vt{\vartheta}

\redefine\G{\Gamma}
\redefine\D{\Delta}
\define\Om{\Omega}

\redefine\L{\Lambda}
\define\Ph{\Phi}
\define\Ps{\Psi}

\redefine\aa{\bold a}

\define\boc{\bold c}
\define\dd{\bold d}
\define\ee{\bold e}

\define\kk{\bold k}

\define\qq{\bold q}

\redefine\ss{\bold s}
\redefine\tt{\bold t}

\redefine\AA{\bold A}

\define\CC{\bold C}

\define\EE{\bold E}
\define\FF{\bold F}

\define\II{\bold I}

\define\NN{\bold N}

\define\QQ{\bold Q}
\define\RR{\bold R}

\define\TT{\bold T}

\define\WW{\bold W}
\define\ZZ{\bold Z}

\define\ca{\Cal A}
\define\cb{\Cal B}
\define\cc{\Cal C}
\define\cd{\Cal D}
\define\ce{\Cal E}
\define\cf{\Cal F}
\define\cg{\Cal G}
\define\ch{\Cal H}

\define\cj{\Cal J}
\define\ck{\Cal K}
\define\cl{\Cal L}

\define\co{\Cal O}
\define\cp{\Cal P}

\define\car{\Cal R}

\define\ct{\Cal T}

\define\cv{\Cal V}

\define\cz{\Cal Z}
\define\cx{\Cal X}

\define\fe{\frak e}

\define\fs{\frak s}

\define\fA{\frak A}

\define\fE{\frak E}
\define\fF{\frak F}

\define\fK{\frak K}

\define\fP{\frak P}

\define\fS{\frak S}

\define\fU{\frak U}

\define\ta{\ti a}

\define\tc{\ti c}
\define\td{\ti d}

\define\ty{\ti y}

\define\tB{\ti B}

\define\tD{\ti D}
\define\tE{\ti E}

\define\tG{\ti G}
\define\tH{\ti H}

\define\tK{\ti K}

\define\tT{\ti T}

\define\tV{\ti V}
\define\tW{\ti W}

\define\tZ{\ti Z}
\define\tPh{\ti{\Ph}}
\define\sh{\sharp}
\define\Mod{\text{\rm Mod}}

\define\bp{\bar p}

\define\bg{\bar g}

\define\bS{\bar S}

\define\tce{\ti\ce}
\define\bul{\bullet}

\define\che{\check}
\define\cha{\che{\a}}
\define\bfU{\bar\fU}

\define\prq{\preceq}

\define\dbbq{\dot{\bbq}}
\define\ale{\aleph}
\define\BBD{BBD}
\define\DL{DL}
\define\KL{KL1}
\define\KLL{KL2}
\define\CLA{L7}
\define\ORT{L15}
\define\ORA{L14}
\define\CS{L3}
\define\UNE{L12}
\define\CDG{L9}
\define\OS{Os}
\define\SH{Sh}
\define\SP{S}

\head Introduction\endhead
Throughout this paper, $G$ denotes a fixed, not necessarily connected, reductive algebraic group over an 
algebraically closed field $\kk$ with a fixed connected component $D$ which generates $G$. This paper is a part 
of a series \cite{\CDG} which attempts to develop a theory of character sheaves on $D$.

Our main result here is the classification of "unipotent" character sheaves on $D$ (under a mild assumption on the
characteristic of $\kk$). This extends the results of \cite{\CS, IV,V} which applied to the case where $G=G^0$. 
While in the case of $G=G^0$ the classification of unipotent character sheaves is essentially the same as the 
classification of unipotent representations of a split connected reductive group over $\FF_q$, the classification
in the general case is essentially the same as the classification of unipotent representations of a not 
necessarily split connected reductive group over $\FF_q$ given in \cite{\ORA}. 

We now describe the content of the various sections in more detail. \S43 contains some preparatory material 
concerning (extended) Hecke algebra and two-sided cells which are used later in the study of unipotent character 
sheaves. In \S44 we study the unipotent character sheaves in connection with Weyl group representations and 
two-sided cells. (But it turns out that the method of associating a two-sided cell to a unipotent character sheaf
along the lines of \cite{\CS, III} is better for our purposes than the one in \S41.) A number of results in this 
section are conditional (they depend on a cleanness property and/or on a parity property); they will become 
unconditional in \S46. In \S45 we show that the problem of classifying the unipotent character sheaves on $D$ can
be reduced to the analogous problem in the case where $G^0$ is simple and $G$ has trivial centre. In \S46 we 
extend the results of \cite{\CS, IV,V} on the classification of unipotent character sheaves on $D$ from the case 
$G=G^0$ to the general case.

{\it Erratum to }\cite{\CDG, V}; in line 4 of 25.1 replace last $a$ by $s$.

{\it Erratum to }\cite{\CDG, VI}: on p.383 l.-25,-24 replace $Z$ by ${}'\bar Z^\ss$ and 
$\Delta^0_j$ by $\Delta_j$.

{\it Erratum to }\cite{\CDG, VII}: on p.248, l.4 of 35.5 replace $G^0F$ by $G^{0F}$.

{\it Erratum to }\cite{\CDG, VIII}: on p.346, l.14 replace the first $k$ by $k'$;
on p.350, l.3 and l.4 of 39.6 delete "The restriction of", "to";
on p.350, l.6 of 39.6 replace first $\s$ by $x$; the 5 lines preceding 39.8 
("If $n=3$ then ... is proved") should be replaced  by the following text:

"If $n=3$ then $W$ must be of type of type $D_4$, $\G$ is the alternating group in four letters $a,b,c,d$, 
$W^{(K)}$ is either $\{1\}$ or $\ZZ/2$ (with trivial $\G$-action) or the $\ZZ/2$-vector space spanned by 
$a,b,c,d$ with the obvious $\G$-action. It is enough to show that $E=E'\ot E''$ where $E'$ is a 
$\bfU[W^{(K)}\G]$-module defined over $\QQ$ and $E''$ is a $\bfU[W^{(K)}\G]$-module of dimension $1$ over $\bfU$.
If $W^{(K)}\G$ has order $\le2$, this follows from the fact that

(a) any simple $\bfU[\G]$-module is either defined over $\QQ$ or has dimension $1$. 
\nl
(Indeed, if it has dimension $>1$ then it is the restriction to $\G$ of the $3$-dimensional reflection 
representation of the symmetric group in four letters, which is defined over $\QQ$.) 

Now assume that $W^{(K)}$ has order $>2$. We can find a homomorphism $\e:W^{(K)}@>>>\bfU^*$ (with image in 
$\{1,-1\}$) whose stabilizer in $\G$ is denoted by $\G_\e$ and a simple $\bfU[\G_\e]$-module $E_0$ such that 
$E=\Ind_{W^{(K)}\G_\e}^{W^{(K)}\G}(E_\e\bxt E_0)$; here $E_\e$ is the one dimensional $\bfU[W^{(K)}]$-module 
defined by $\e$ (necessarily defined over $\QQ$). If $\G_\e=\G$ then $E=E_\e\bxt E_0$ where $E_0$ is as in (a) and
the desired result follows. If $\G_\e$ has order $2$ then $E_0$ is defined over $\QQ$ hence $E$ is defined over 
$\QQ$. If $\G_\e\ne\G$ and $\G_\e$ is not of order $2$ then $\G_\e$ is of order $3$, $E_0$ is the restriction to 
$\G_\e$ of a one dimensional $\bfU[\G]$-module $E''$ and we have $E=E'\ot E''$ where 
$E'=\Ind_{W^{(K)}\G_\e}^{W^{(K)}\G}(E_\e\bxt\bfU)$ is defined over $\QQ$. Hence the proposition holds in this 
case. The proposition is proved."

{\it Erratum to }\cite{\CDG, IX}: on p.354, l.-8 replace $V_\l$,$V_\l^D$ by $\Om_\l,\Om_\l^D$;
on p.354, l.-7 replace 34.4 by 34.2; on p.355, l.-8, l.-13 replace $V_\l$ by $\Om_\l$;
on p.359, first line of 40.8 replace $c_{y.\l}$ by $c_{y,\nu}$; on the preceding line replace
$in$ by $\in$; on p.361, l.9 insert "," before $\cl$; on p.363, l.6 before "Let" insert:
"Let $\dot{\cl}_w^{\sh}=IC(\bar Z^w_{\em,D},\dot{\cl}_w)$."; on p.365, second line of 41.4,
two $)$ are missing; on p.366, last displayed line of 41.4 replace $\dsv A$ by $\dsv\fe(A)$;
on p.368, l.2 remove "the condition that"; on p.369, l.7 a $)$ is missing;
on p.371, l.1 replace $H_n$ by $H$; on p.372, l.4 of 42.5 replace $\ot\ca$ by $\ot_{\ca}$;
on p.376, l.-22 replace $WW$ by $\WW$; on p.376, l.-10 replace $H_n^{D,\tca}D$ by $H_n^{D,\tca}$;
on p.377, l.-10 replace $vt$ by $\vt$; on p.378, l.6 replace $\Delta$ by $D$.

{\it Notation.} Let $\e:=\e_D$ be as in 26.2. If $X$ is an algebraic variety over $\kk$ and $K\in\cd(X)$ we write
$H^i(K)$ instead of ${}^pH^i(K)$. If $K\in\cd(X)$ we set $gr_1K=\sum_{i\in\ZZ}(-1)^iH^i(K)$, an element of the 
Grothendieck group of the category of perverse sheaves on $X$. The cardinal of a finite set $X$ is denote by 
$|X|$.

\head Contents\endhead
43. Preparatory results on Hecke algebras.

44. Unipotent character sheaves and two-sided cells.

45. Reductions.

46. Classification of unipotent character sheaves.

\head 43. Preparatory results on Hecke algebras\endhead   
\subhead 43.1\endsubhead
This section contains some preparatory material concerning (extended) Hecke algebra and two-sided cells which
will be used later in the study of unipotent character sheaves.

We fix an even integer $c\ge2$ which is divisible by $|G/G^0|$. Let $\G$ be a cyclic group of order $c$ with
generator $\vp$. Let $\tWW$ be the semidirect product of $\WW$ with $\G$ (with $\WW$ normal) where 
$\vp x\vp\i=\e(x)$ for $x\in\WW$. Note that the group $\WW^D$ in 34.2 is naturally a quotient of $\tWW$, via 
$x\vp^i\m x\uD^i$ with $x\in\WW$, $i\in\ZZ$. Let $\Irr(\WW)$ be the category whose objects are the simple 
(or equivalently, absolutely simple) $\QQ[\WW]$-modules. Let $\Irr^\e(\WW)$ be the category whose objects are the
simple $\QQ[\WW]$-modules $E_0$ such that $\tr(x,E_0)=\tr(\e(x),E_0)$ for all $x\in\WW$. Let $\Mod(\tWW)$ be the 
category whose objects are the $\QQ[\tWW]$-modules of finite dimension over $\QQ$. Let $\Irr(\tWW)$ be the 
subcategory of $\Mod(\tWW)$ consisting of those objects that remain simple on restriction to $\QQ[\WW]$. Let 
$\un\Irr(\tWW)$ be a set of representatives for the isomorphism classes in $\Irr(\tWW)$. Let $\io$ be the object 
of $\Irr(\tWW)$ whose underlying $\QQ$-vector space is $\QQ$ with $\WW$ acting trivially and $\vp$ acting as 
multiplication by $-1$. Note that if $E\in\Irr(\tWW)$ then $E|_{\QQ[\WW]}\in\Irr^\e(\WW)$. Conversely, we show:

(a) {\it for any $E_0\in\Irr^\e(\WW)$, the set $\{E\in\un\Irr(\tWW);E|_{\QQ[\WW]}\cong E_0\}$ has exactly two 
elements; one is isomorphic to the other tensored with $\io$.}
\nl
From \cite{\ORA, 3.2} we see that there exists a linear map of finite order $\g:E_0@>>>E_0$ such that 
$\g(x(e))=\e(x)(\g(x))$ for any $e\in E_0$, $x\in\WW$. (We use the following property of $\e$: if $s,s'\in\II$ are
such that $ss'$ has order $\ge4$ then $s,s'$ are in distinct orbits of $\e$ on $\II$.) Moreover, from the proof
in \cite{\ORA, 3.2} we see that $\g$ can be chosen so that $\g^{c'}=1$ where $c'$ is the order of the permutation
$\e:\WW@>>>\WW$. In particular we have $\g^c=1$. This proves that the set in (a) is nonempty. The remainder of (a)
is immediate.

Let $\fE$ be a subset of $\un\Irr(\tWW)$ such that $\{E|_{\QQ[\WW]};E\in\fE\}$ represents each isomorphism 
class in $\Irr^\e(\WW)$ exactly once. 

\subhead 43.2\endsubhead   
Recall the notation $\ca=\ZZ[v,v\i]$. Define $l:\tWW@>>>\NN$ by $l(x\vp^i)=l(x)$ for $x\in\WW,i\in\ZZ$; here 
$l:\WW@>>>\NN$ is the standard length function. Let $w_0$ be the longest element of $\WW$. Let $\tH$ be the 
$\ca$-algebra with $1$ with generators $\tT_w(w\in\tWW)$ and relations

$\tT_w\tT_{w'}=\tT_{ww'}$ for $w,w'\in\tWW$ with $l(ww')=l(w)+l(w')$,

$\tT_s^2=\tT_1+(v-v\i)\tT_s$ for $s\in\II$.
\nl
We have a surjective $\ca$-algebra homomorphism $\z:\tH@>>>H^D_1$, $\tT_{x\vp^i}\m\tT_{x\uD^i}$ for $x\in\WW$,
$i\in\ZZ$ where $H^D_1$ is the algebra $H^D_n$ in 34.4 (with $n=1$); thus, a number of properties of $\tH$ can be
deduced from the corresponding properties of $H^D_1$ in \S34. 

Let $\x\m\x^\da$ be the $\ca$-algebra isomorphism $\tH@>>>\tH$ such that $\tT_w^\da=(-1)^{l(w)}\tT_{w\i}\i$ for 
all $w\in\tWW$. Let $\bar{}:\ca@>>>\ca$ be the ring isomorphism such that $\ov{v^i}=v^{-i}$ for $i\in\ZZ$. Let 
$\bar{}:\tH@>>>\tH$, $\x\m\bar\x$ be the ring isomorphism such that $\ov{a\tT_w}=\ov{a}\tT_{w\i}\i$ for $w\in\tWW$,
$a\in\ca$; this isomorphism commutes with $\x\m\x^\da$. For $w\in\tWW$ we set
$$c_w=\sum_{y\in\WW;y\le x}v^{l(y)-l(x)}P_{y,x}(v^2)\tT_{y\vp^i}\in\tH,$$
$$\tc_w=\sum_{z\in\WW;x\le z}(-1)^{l(z)-l(x)}v^{l(x)-l(z)}P_{w_0z,w_0x}(v^2)\tT_{z\vp^i}\in\tH,$$
where $w=x\vp^i$ ($x\in\WW,i\in\ZZ$) and
$$P_{y,x}(\qq)=\sum_{j\in\ZZ}n_{y,x,j}\qq^{j/2},\qua n_{y,x,j}\in\ZZ$$
are the polynomials defined in \cite{\KL} for the Coxeter group $\WW$. Note that $n_{y,x,j}=0$ unless $j\in2\ZZ$
and $n_{x,x,j}=\d_{j,0}$. We have $c_w=c_x\tT_{\vp^i}$ and $\ov{c_w}=c_w$. It follows that
$$c_w^\da=\sum_{y\in\WW;y\le x}(-1)^{l(y)}v^{-l(y)+l(x)}P_{y,x}(v^{-2})\tT_{y\vp^i}\in\tH$$
and $\ov{c_w^\da}=c_w^\da$.

Let $\tH^v=\QQ(v)\ot_\ca\tH$, a $\QQ(v)$-algebra. Let $\tH^1=\QQ\ot_\ca\tH$ where $\QQ$ is regarded as an 
$\ca$-algebra under $v\m1$. We have $\tH^1=\QQ[\tWW]$ (with $\tT_w\in\tH^1$ identified with $w\in\QQ[\tWW]$ for 
$w\in\tWW$). Let $\x\m\x|_{v=1}$ be the ring homomorphism $\tH@>>>\tH^1$ given by $v\m1$, $\tT_w\m w$ for 
$w\in\tWW$.

Let $H,H^v,H^1$ be the algebras defined like $\tH,\tH^v,\tH^1$ by replacing $\tWW$ by $\WW$. We identify 
$H,H^v,H^1$ with subalgebras with $1$ of $\tH,\tH^v,\tH^1$ in an obvious way. We have $H^1=\QQ[\WW]$. Note that 
$H$ is the same as the algebra $H_n$ in 31.2 (with $n=1$).

For $x,y\in\WW$ we have $c_xc_y=\sum_{z\in\WW}r_{x,y}^zc_z$ with $r_{x,y}^z\in\ca$. There is a well defined 
function $\aa:\WW@>>>\NN$ such that for any $x,y,z\in\WW$ we have $r_{x,y}^z\in v^{\aa(z)}\ZZ[v\i]$ and for any 
$z\in\WW$ we have $r_{x,y}^z\n v^{\aa(z)-1}\ZZ[v\i]$ for some $x,y\in\WW$. For any $x,y,z\in\WW$ we define 
$\g_{x,y,z\i}\in\ZZ$ by $r_{x,y}^z=\g_{x,y,z\i}v^{\aa(z)}\mod v^{\aa(z)-1}\ZZ[v\i]$. 

We define a preorder $\prq$ on $\WW$ as follows: we say that $x'\prq x$ if there exists $x_1,x_2$ in $\WW$ such
that in the expansion (in $H$) $c_{x_1}c_xc_{x_2}=\sum_{y'\in\WW}r_{y'}c_{y'}$ with $r_{y'}\in\ca$ we have
$r_{x'}\ne0$. Let $\si$ be the equivalence relation on $\WW$ attached to $\prq$. The equivalence classes for $\si$
are called the two-sided cells of $\WW$. (See also \cite{\KL}.) We write $x\prec y$ instead of $x\prq y$, 
$x\not\si y$. It is known that $\aa:\WW@>>>\NN$ is constant on each two-sided cell. If $\boc,\boc'$ are two-sided
cells we write $\boc\prq\boc'$ instead of $x\prq x'$ for some/any $x\in\boc$, $x'\in\boc'$. This is a partial 
order on the set of two-sided cells; we also write $\boc\prec\boc'$ instead of $\boc\prq\boc',\boc\ne\boc'$.

The free abelian group $H^\iy$ with basis $\{t_x;x\in\WW\}$ is regarded as a ring with multiplication given by
$t_xt_y=\sum_{z\in\WW}\g_{x,y,z\i}t_z$ for $x,y\in\WW$. This ring has a unit element of the form 
$\sum_{\d\in\cd}t_\d$ where $\cd$ is a well defined subset of $\WW$. We have $H^\iy=\op_\boc H^\iy_\boc$ (as 
rings) where $\boc$ runs over the two-sided cells and $H^\iy_\boc$ is the subgroup of $H^\iy$ generated by 
$\{t_x;x\in\boc\}$. Let $\tH^\iy$ be the free abelian group with basis $\{t_x\vp^i;x\in\WW,i\in[0,c-1]\}$. We have
naturally $H^\iy\sub\tH^\iy$ ($t_x=t_xv^0$). The group ring $\ZZ[\G]$ is also naturally contained in $\tH^\iy$ by
$\vp^i\m\sum_{d\in\cd}t_dv^i$. We regard $\tH^\iy$ as a ring with $1$ so that $H^\iy$ and $\ZZ[\G]$ are subrings 
with $1$ and $\vp t_x\vp\i=t_{\vp(x)}$ for $x\in\WW$. We have a surjective ring homomorphism 
$\z^\iy:\tH^\iy@>>>H_1^{D,\iy}$, $t_x\vp^i\m t_{x\uD^i}$ for $x\in\WW,i\in\ZZ$ where $H_1^{D,\iy}$ is the ring 
$H_n^{D,\iy}$ (with $n=1$) in 34.12.

Define $\ca$-linear maps $\Ph:H@>>>\ca\ot H^\iy$, $\tPh:\tH@>>>\ca\ot\tH^\iy$ by
$\Ph(c_x^\da)=\sum_{z\in\WW,d\in\cd,\aa(d)=\aa(z)}r_{x,d}^zt_z$ for $x\in\WW$, 
$\tPh(c_{x\vp^i}^\da)=\Ph(c_x^\da)\vp^i$ for $x\in\WW,i\in\ZZ$. Now $\Ph,\tPh$ are homomorphisms of rings with 
$1$. 
We have a commutative diagram
$$\CD
\tH@>>>\ca\ot\tH^\iy  \\
@V\z VV       @V\z^\iy VV    \\
H_1^D@>>>\ca\ot H_1^{D,\iy}
\endCD$$
where the upper horizontal map is the composition of ${}^\da:\tH@>>>\tH$ with $\tPh$ and the lower horizontal map
is the map denoted by $\Ph$ in 34.1, 34.12 (which is not the same as the present $\Ph$).

For any field $k$ let $H^\iy_k=k\ot H^\iy$, $\tH^\iy_k=k\ot\tH^\iy$. Let $\Ph^v:H^v@>>>H^\iy_{\QQ(v)}$, 
$\tPh^v:\tH^v@>>>\tH^\iy_{\QQ(v)}$ be the $\QQ(v)$-algebra homomorphisms obtained from $\Ph,\tPh$ by extension of
scalars. Let $\Ph^1:H^1@>>>H^\iy_\QQ$, $\tPh^1:\tH^1@>>>\tH^\iy_\QQ$ be the $\QQ$-algebra homomorphisms obtained 
from $\Ph,\tPh$ by extension of scalars. Now $\Ph^v,\tPh,\Ph^1,\tPh^1$ are algebra isomorphisms. Since the 
$\QQ$-algebra $\QQ[\WW]=H^1$ is split semisimple, the same holds for the $\QQ$-algebra $H^\iy_\QQ$.

Now $\x\m\x^\da$ induces by extension of scalars a $\QQ(v)$-algebra isomorphism $\tH^v@>>>\tH^v$ and a 
$\QQ$-algebra isomorphism $\tH^1@>>>\tH^1$; these leave $H^v,H^1$ stable and are denoted again by $\x\m\x^\da$.

\subhead 43.3\endsubhead
Let $E_0\in\Irr(\WW)$. We can view $E_0$ as a simple $H^\iy_\QQ$-module $E_0^\iy$ via $\Ph^1$. Now 
$\QQ(v)\ot_\QQ E_0^\iy$ is naturally a simple $H^\iy_{\QQ(v)}$-module and this can be viewed as a simple 
$H^v$-module $E_0^v$ via $\Ph^v$. 

Let $E\in\Irr(\tWW)$. We can view $E$ as a simple $\tH^\iy_\QQ$-module $E^\iy$ via $\tPh^1$. Now 
$\QQ(v)\ot_\QQ E^\iy$ is naturally a $\tH^\iy_{\QQ(v)}$-module and this can be viewed as an $\tH^v$-module $E^v$ 
via $\tPh^v$. By restriction, $E$ can be viewed as a simple $\QQ[\WW]=H^1$-module $E_0$. From the definitions we 
see that $E_0^v$ is the restriction of the $\tH^v$-module $E^v$ to $H^v$.

Let $E'$ be the $\QQ[\tWW]$-module with the same underlying $\QQ[\WW]$-module structure as $E$ but with action of
$\vp$ equal to $-1$ times the action of $\vp$ on $E$. Then $E'{}^v$ is defined. Clearly, $E'{}^v,E^v$ have the
same underlying $H^v$-module and the action of $\tT_\vp$ on $E'{}^v$ is equal to $-1$ times the action of 
$\tT_\vp$ on $E^v$. 

Let $\sgn$ be the object of $\Irr(\tWW)$ with underlying vector space $\QQ$ on which $w\in\tWW$ acts as 
multiplication by $(-1)^{l(w)}$. We set $E^\da=E\ot\sgn\in\Irr(\tWW)$.

\subhead 43.4\endsubhead
Let $E\in\Irr(\tWW)$. From the definitions, for any $\x\in\tH$, $\z\in\tH^\iy$ we have:
$$\tr(\x,E^v)\in\ca,\qua\tr(\x,E^v)|_{v=1}=\tr(\x|_{v=1},E),\qua\tr(\z,E^\iy)\in\ZZ.\tag a$$
Hence it makes sense to write
$$\tr(\x,E^v)=\sum_{i\in\ZZ}\tr(\x,E^v;i)v^i\text{ where }\tr(\x,E^v;i)\in\ZZ.$$
More generally for $\x\in\tH^v$ we write 
$\tr(\x,E^v)=\sum_{i\in\ZZ}\tr(\x,E^v;i)v^i$ (possibly infinite sum) where $\tr(\x,E^v;i)\in\QQ$ (here 
$\tr(\x,E^v)\in\QQ(v)$ is viewed as a power series in $\QQ((v))$).

For any $\x\in\tH$ we show: 
$$\tr(\x,(E^\da)^v)=\tr(\x^\da,E^v).\tag b$$
Let $E^{v\da}$ be the $\tH^v$-module whose underlying $\QQ(v)$-module is $E^v$ but with $\x\in\tH^v$ acting as 
$\x^\da$ in the $\tH^v$-module $E^v$. Note that the $\tH^v$-module $E^{v\da}$ is simple and its restriction to an
$H^v$-module is simple. Also, the assignment $E'\m E'{}^v$ defines a bijection between the set of isomorphism
classes of objects of $\Irr(\tWW)$ and the set of isomorphism classes of simple $\tH^v$-modules whose restriction
to $H^v$ is simple. Thus we have $E^{v\da}\cong E_1^v$ for some $E_1\in\Irr(\tWW)$. It is enough to show that 
$(E^\da)^v\cong E^{v\da}$ or that $(E^\da)^v\cong E_1^v$ as $\tH^v$-modules. Using (a) for $\x\in\tH$ we have 
$$\align&\tr(\x_{v=1},E_1)=\tr(\x,E_1^v)_{v=1}=\tr(\x,E^{v^\da})_{v=1}=\tr(\x^\da,E^v)_{v=1}\\&
=\tr(\x^\da|_{v=1},E)=\tr(\x|_{v=1},E\ot\sgn).\endalign$$
Thus, $\tr(w,E_1)=\tr(w,E^\da)$ for any $w\in\tWW$ so that $E_1\cong E^\da$ in $\Irr(\tWW)$ and 
$(E^\da)^v\cong E_1^v$, as required.

For any $w\in\tWW$ we have:
$$\tr(\tT_w\i,E^v)=\tr(\tT_w,E^v).\tag c$$
The proof is the same as that of 34.17 (we use also (a)).

For any $\x\in\tH$ we show:
$$\tr(\ov{\x},E^v)=\ov{\tr(\x,E^v)}.\tag d$$
We may assume that $\x=c_{x\vp^j}^\da$ with $x\in\WW,j\in\ZZ$. Since $\ov{\x}=\x$, it is enough to verify:
$$\sum_{z\in\WW,d\in\cd,\aa(d)=\aa(z)}r_{x,d}^z\tr(t_z\vp^j,E^\iy)=
\sum_{z\in\WW,d\in\cd,\aa(d)=\aa(z)}\ov{r_{x,d}^z}\tr(t_z\vp^j,E^\iy).$$
This follows from the obvious identity $r_{x,y}^z=\ov{r_{x,y}^z}$ for any $x,y,z\in\WW$.

For any $w\in\tWW$ we show:
$$\tr(\tT_w,(E^\da)^v)=(-1)^{l(w)}\ov{\tr(\tT_w,E^v)}.\tag e$$
Using (b),(d), we see that the left hand side of (e) equals
$$(-1)^{l(w)}\tr(\tT_{w\i}\i),E^v)=(-1)^{l(w)}\tr(\ov{\tT_w},E^v)=(-1)^{l(w)}\ov{\tr(\tT_w,E^v)}.$$
This proves (e).

\subhead 43.5\endsubhead
For $E\in\Irr(\tWW)$ we define $f_E^v\in\QQ[v,v\i]$, $f_E^\iy\in\QQ$ by
$$\sum_{x\in\WW}\tr(\tT_x,E^v)^2=f_E^v\dim E, \sum_{x\in\WW}\tr(t_x,E^\iy)^2=f_E^\iy\dim E.\tag a$$
Note that $f_E^v,f_E^\iy$ depend only in $E|_{\QQ[\WW]}$. Now $f_E^v$ is $\ne0$; it specializes to $|\WW|/\dim E$
for $v=1$. Since $E_0^\iy$ is a simple $H^\iy_\QQ$-module, the integer $\tr(t_x,E_0^\iy)$ is $\ne0$ for some 
$x\in\WW$. Hence $f_E^\iy\ne0$. For $E,E'$ in $\Irr(\tWW)$, the following holds:

(b) $\sum_{x\in\WW}\tr(\tT_{x\vp},E^v)\tr(\tT_{x\vp},E'{}^v)$ equals $f_E^v\dim E$ if $E,E'$ are isomorphic and 
equals $0$ if $E|_{\QQ[\WW]}\not\cong E'|_{\QQ[W]}$.
\nl
This can be deduced from 34.15(c) using the commutative diagram in 43.2 (we use also 43.4(a)). Similarly,

(c) $\sum_{x\in\WW}\tr(x\vp,E)\tr(x\vp,E')$ equals $|\WW|$ if $E,E'$ are isomorphic and 
equals $0$ if $E|_{\QQ[\WW]}\not\cong E'|_{\QQ[W]}$.

\subhead 43.6\endsubhead
Let $E_0\in\Irr(\WW)$. Let $E_0^\iy$ be the irreducible $H^\iy_\QQ$-module corresponding to $E_0$ as in 43.3. 
Since $H^\iy_\QQ=\op_\boc\QQ\ot H^\iy_\boc$ as $\QQ$-algebras, there is a unique two-sided cell $\boc=\boc_{E_0}$
such that $E_0^\iy$ restricts to a simple module of the summand $\QQ\ot H^\iy_\boc$ (and all other summands act as
$0$ on $E_0^\iy$). Let $a_{E_0}$ be the value of $\aa$ on $\boc_{E_0}$.

Let $E\in\Irr(\tWW)$. We set $\boc_E=\boc_{E_0}$, $a_E=a_{E_0}$ where $E_0=E|_{\QQ[\WW]}\in\Irr(\WW)$. We show:

(a) {\it if $x\in\WW$, then $\tr(c_{x\vp}^\da,E^v)=\tr(t_x\vp,E^\iy)v^{-a_E}\mod v^{-a_E+1}\ZZ[v]$; equivalently,
$\tr(c_{x\vp}^\da,E^v;-a_E)=\tr(t_x\vp,E^\iy)$ and $\tr(c_{x\vp}^\da,E^v;\ta)=0$ for all $\ta<-a_E$;}

(b) {\it if $x\in\WW$ and the action of $c_{x\vp}^\da$ on $E^v$ is $\ne0$, then $z\prq x$ for some $z\in\boc_E$.}
\nl
From the definition, 
$$\tr(c_{x\vp}^\da,E^v)=\sum_{z\in\WW,d\in\cd,\aa(d)=\aa(z)}r_{x,d}^z\tr(t_z\vp,E^\iy).$$
In the last sum we have $\tr(t_z\vp,E^\iy)=0$ unless $z\in\boc_E$ in which case $\aa(z)=a_E$. For such $z$ we have
$r_{x,d}^z=\g_{x,d,z\i}v^{a_E}\mod v^{a_E-1}\ZZ[v\i]$ hence
$r_{x,d}^z=\ov{r_{x,d}^z}=\g_{x,d,z\i}v^{-a_E}\mod v^{-a_E+1}\ZZ[v]$ and
$$\tr(c_{x\vp}^\da,E^v)=\sum_{z\in\WW}\d_{x,z}\tr(t_z\vp,E^\iy)v^{-a_E}\mod v^{-a_E+1}\ZZ[v]$$
and (a) follows.

In the setup of (b), the action of $\sum_{z\in\WW,d\in\cd,\aa(d)=\aa(z)}r_{x,d}^zt_z\vp$ on $E^\iy$ is $\ne0$. 
Hence there exist $z\in\boc_E,d\in\cd$ such that $r_{x,d}^z\ne0$ (so that $z\prq x$). This proves (b).

We show:

(c) {\it if $x\in\WW$, then $\tr(\tT_{x\vp},E^v;-a_E)=\sgn(x)\tr(t_x\vp,E^\iy)$ and $\tr(\tT_{x\vp},E^v;\ta)=0$ 
for all $\ta<-a_E$.}
\nl
We argue by induction on $l(x)$. If $l(x)=0$ we have $x=1$ and $\tT_{x\vp}=c_{x\vp}^\da$ and the result follows 
from (a). Assume now that $l(x)>0$. From the definition we have $c_{x\vp}^\da=\sgn(x)\tT(x\vp)+\x$ where
$\x\in\sum_{x';l(x')<l(x)}v\ZZ[v]\tT_{x'\vp}$. The induction hypothesis shows that $\tr(\x,E^v;\ta)=0$ for all
$\ta\le-a_E$. Hence $\sgn(x)\tr(\tT_{x\vp},E^v;\ta)=\tr(c_{x\vp}^\da,E^v;\ta)$ for all $\ta\le-a_E$; now (c) for 
$x$ follows from (a).

Using (c) and 43.5(b) we see that
$$f_E^v\dim E=\sum_{x\in\WW}\tr(t_x\vp,E^\iy)^2v^{-2a_E}+\text{strictly higher powers of $v$}.$$
Using now 43.5(a) we obtain
$$f_E^v=f_E^\iy v^{-2a_E}+\text{strictly higher powers of $v$}.\tag d$$

Now let $E'$ be another object of $\Irr(\tWW)$. We show:

(e) {\it $\sum_{x\in\WW}\tr(t_x\vp,E^\iy)\tr(t_x\vp,E'{}^\iy)$ is equal to $f_E^\iy\dim E$ (if $E,E'$ are isomorphic) 
and is equal to $0$ if $E'|_{\QQ[\WW]}\not\cong E_0$.}
\nl
We can assume that $\boc_{E'}=\boc_E$ (otherwise, the sum in (e) is $0$). Combining 43.5(b) with (c) for $E$
and $E'$ and with (d) we see that 

$v^{-2a_E}\sum_{x\in\WW}\tr(t_x\vp,E^\iy)\tr(t_x\vp,E'{}^\iy)$
\nl
plus a $\ZZ$-linear combination of strictly higher powers of $v$ is equal to $f_E^\iy\dim Ev^{-2a_E}$ plus a 
$\ZZ$-linear combination of strictly higher powers of $v$ (if $E,E'$ are isomorphic) and is equal to $0$ if 
$E'|_{\QQ[\WW]}\not\cong E_0$. Taking coefficients of $v^{-2a_E}$ we obtain (e).

We show:  
$$\e(\boc_E)=\boc_E.\tag f$$
For any $x\in\WW$ we have $\tr(\e(x),E_0)=\tr(x,E_0)$. It follows that for any $x\in\WW$ we have
$\tr(t_{\e(x)},E^\iy_0)=\tr(t_x,E_0^\iy)$. We can find $x\in\boc_E$ such that $\tr(t_x,E_0^\iy)\ne0$. Then 
$\tr(t_{\e(x)},E^\iy_0)\ne0$ hence $\e(x)\in\boc_E$ and (f) follows.

\subhead 43.7\endsubhead
Let $(W,S)$ be a Weyl group ($S$ is the set of simple reflections). Let $\s:W@>>>W$ be an automorphism of $W$ such
that $\s(I)=I$ and such that whenever $s\ne s'$ in $S$ are in the same orbit of $\s$, the product $ss'$ has order
$2$ or $3$. Let $b\in\ZZ_{>0}$ be such that $\s^b=1$. Let $\tW$ be the semidirect product of $W$ with the cyclic 
group $C$ of order $b$ with generator $\s$ so that in $\tW$ we have the identity $\s x\s\i=\s(x)$ for 
any $x\in W$. Let $I$ be a $\s$-stable subset of $S$ and let $W_I$ be the subgroup of $W$ generated by $I$. Let 
$E$ be a simple $\QQ[\tW]$-module such that $E|_{\QQ[W]}$ is simple. Let $\tW_I=W_IC$, a subgroup of 
$\tW$. Let $E_{\bbq}=\bbq\ot E$. We show:

(a) {\it The $\bbq[\tW_I]$-module $E_{\bbq}|_{\tW_I}$ is isomorphic to $\op_jE'_j$ where each $E'_j$ is a 
$\bbq[\tW_I]$-module and either 

$E'_j$ is induced from a $\bbq[W_IC']$-module where $C'$ is a proper subgroup of $C$, or 

$E'_j|_{W_I}$ is simple and $E'_j$ is defined over $\QQ$.}
\nl
The general case reduces immediately to the case where $\s$ permutes transitively the irreducible components of 
$W$. In this case we may identify $W$ with $W_1\T W_1\T\do\T W_1$ and $S=S_1\T S_1\T\do\T S_1$ ($t$ factors) where
$(W_1,S_1)$ is an irreducible Weyl group; the automorphism $\s$ may be written as
$\s(w_1,w_2,\do,w_t)=(\s'(w_t),w_1,w_2,\do,w_{t-1})$, $w_i\in W_1$ where $\s'$ is an automorphism of $(W_1,S_1)$.
We have $I=I_1\T I_1\T\do I_1$ where $I_1\sub I$ is $\s'$-stable. Hence $W_I=W_{I_1}\T W_{I_1}\T\do\T W_{I_1}$. 
Note that $b/t\in\ZZ_{>0}$. Let $\tW_1$ be the semidirect product of $W_1$ with the cyclic group $C_1$ of order 
$b/t$ with generator $\s'$ so that in $\tW_1$ we have the identity $\s'x_1\s'{}\i=\s'(x_1)$ for any 
$x_1\in W_1$. We can find a simple $\QQ[\tW_1]$-module $E_1$ such that $E_1|_{W_1}$ is simple and such that 
$E=E_1\bxt E_1\bxt\do\bxt E_1$ ($t$ factors) as a $\QQ[W_1]$-module and $\s$ acts on $E$ as
$e_1\bxt e_2\bxt\do\bxt e_t\m\s'(e_t)\bxt e_1\bxt e_2\bxt\do\bxt e_{t-1}$, $(e_i\in E_i)$. Let 
$\tW_{I_1}=W_{I_1}C_1$, a subgroup of $\tW_1$.

Assume that (a) holds when $W,S,\s,b,I,E$ are replaced by $W_1,S_1,\s',b/t,I_1,E_1$. Let $E_{1,\bbq}=\bbq\ot E_1$.
Then we can identify $E_{1,\bbq}|_{\tW_{I_1}}=\op_{j\in\cj}E'_{1,j}$ where each $E'_{1,j}$ is a 
$\bbq[\tW_{I_1}]$-module with properties like those of $E'_j$ in (a). We have 
$E_{\bbq}=\op_{j_1,j_2,\do,j_t\text{ in }\cj}E'_{1,j_1}\bxt E'_{1,j_2}\bxt\do\bxt E'_{1,j_t}$ as a $W_I$-module. 
If we take the sum of all summands where $(j_1,j_2,\do,j_t)$ is fixed up to a cyclic permutation then we obtain a
$\tW_I$-submodule $\ce$ of $E_{\bbq}$. If $j_1,j_2,\do,j_t$ are not all equal then $\ce|_{W_I}$ is induced from a 
$\bbq[W_IC']$-module where $C'$ is a proper subgroup of $C$. If $j_1=j_2=\do=j_t$, then 
$\ce=E'_{1,j_1}\bxt E'_{1,j_1}\bxt\do\bxt E'_{1,j_1}$. If in addition $E'_{1,j_1}$ is induced from a 
$\bbq[W_{I_1}C'_1]$-module where $C'_1$ is a proper subgroup of $C_1$ then $\ce$ is a direct sum of
$\bbq[W_I]$-modules induced from $\bbq[W_IC']$-modules where $C'$ are proper subgroups of $C$. If on the other 
hand $E'_{1,j_1}|_{W_{I_1}}$ is simple and $E'_{1,j_1}$ is defined over $\QQ$ then $\ce|_{W_I}$ is simple and 
$\ce$ is defined over $\QQ$. Thus (a) holds for $W,S,\s,b,I,E$. We can therefore assume that $(W,S)$ is an 
irreducible Weyl group. Let $b'$ be the order of $\s:W@>>>W$. We have $b/b'\in\ZZ_{>0}$. By the proof of 
\cite{\ORA, 3.2} we can find a $\QQ$-linear isomorphism $\s':E@>>>E$ such that $\s'{}^{b'}=1$ and 
$\s'x\s'{}\i=\s(x):E@>>>E$ for any $x\in W$. Since $E|_W$ is absolutely simple we must have 
$\s'=\pm\s:E@>>>E$. Hence if (a) holds when $E$ is modified so that the action of $\s$ is replaced by that of 
$\s'$ (and $b$ is replaed by $b'$) then (a) also holds for the original $E$ and $b$. Thus we may assume that 
$b=b'$. In this case we have $b\le3$. Assume first that $b\le2$. We write $E_{\bbq}|_{\tW_I}=\op_jE'_j$ where each
$E'_j$ is a simple $\bbq[\tW_I]$-module. If $j$ is such that $E'_j|_{W_I}$ is not simple then $E'_j$ is induced by
a $\bbq[W_IC']$-module where $C'$ is a proper subgroup of $C$. If $j$ is such that $E'_j|_{W_I}$ is simple then 
there exists a $\QQ[W_I]$-module $E_0$ of finite dimension over $\QQ$ such that $E'_j|_{W_I}=\bbq\ot E_0$ as 
$\bbq[W_I]$-modules; moreover, by the proof of \cite{\ORA, 3.2}, there exists a $\QQ$-linear isomorphism 
$\ti\s:E_0@>>>E_0$ such that $\ti\s^2=1$ and $\ti\s x\ti\s\i=\s(x)$ for any $x\in W_I$. We extend 
$\ti\s$ to a $\bbq$-linear isomorphism $\bbq\ot E_0@>>>\bbq\ot E_0$ denoted again by $\ti\s$. Since $E_0$ is an 
absolutely simple $W_I$-module we have $\s=a\ti\s:\bbq\ot E_0@>>>\bbq\ot E_0$ where $a\in\bbq^*$. Since 
$\s^2=\ti\s^2=1$ on $\bbq\ot E_0$, we have $a=\pm1$. Hence $\s:\bbq\ot E_0@>>>\bbq\ot E_0$ is defined over $\QQ$.
We see that (a) holds for $E$. Next we assume that $b=3$ so that $W$ is of type $D_4$. In this case (a) is 
verified by examining the known explicit $W$-graph realization of $E$. This completes the proof of (a).

\subhead 43.8\endsubhead
We now return to the setup in 43.1, 43.2. Let $I$ be a subset of $\II$ such that $\e(I)=I$. Let $P\in\cp_I$ (see 
26.1). Then $N_DP\ne\em$ so that $D':=N_DP/U_P$ is a connected component of the reductive group $G':=N_GP/U_P$; 
note that $G'{}^0=P/U_P$. Let $\tWW_I$ be the subgroup of $\tWW$ generated by $\WW_I$ (see 26.1) and $\G$; now 
$\WW_I,I,\tWW_I$ play the same role for $G',D'$ as $\WW,\II,\tWW$ for $G,D$. Let $\tH_I^v$ be the algebra defined
like $\tH^v$ (with $\WW,\II$ replaced by $\WW_I,I$). We have naturally $\tH_I^v\sub\tH^v$ as algebras with $1$. 
For any subgroup $\G'$ of $\G$ let $\tH_I^{v,\G'}$ be the subspace of $\tH_I^v$ spanned by the elements 
$T_{x\vp^i}$ with $x\in\WW_I$ and $i\in\ZZ$ such that $\vp^i\in\G'$; this is a subalgebra of $\tH_I^v$. Let 
$\tH_{I,\bbq}^v,\tH_{\bbq}^v,\tH_{I,\bbq}^{v,\G'}$ be the $\bbq(v)$-algebras obtained by applying 
$\bbq(v)\ot_{\QQ(v)}()$ to $\tH_I^v,\tH^v,\tH_I^{v,\G'}$.

Let $E\in\Irr(\tWW)$. Let $E^v$ be the $\tH^v$-module corresponding to $E$, see 43.3. We have the following 
result:

(a) {\it The restriction to $\tH_{I,\bbq}^v$ of the $\tH_{\bbq}^v$-module $\bbq\ot E^v$ is isomorphic to 
$\op_j\EE'_j$ where each $\EE'_j$ is a $\tH_{I,\bbq}^v$-module and either

(i) $\EE'_j$ is of the form $\tH_{I,\bbq}^v\ot_{\tH_{I,\bbq}^{v,\G'}}\EE''_j$ for some proper subgroup $\G'$ of
$\G$ and some $\tH_{I,\bbq}^{v,\G'}$-module $\EE''_j$, or

(ii) $\EE'_j$ is of the form $\bbq\ot M_j{}^v$ where $M_j\in\Irr(\tWW_I)$;}
\nl
here $M_j^v$ is defined like $E^v$ in terms of $\tWW_I$ instead of $\tWW$.

Note that (a) is a $v$-analogue of 43.7(a). It can be proved by the same method as 43.7(a) or it can be reduced to
43.7(a) with $W=\WW,b=c$.

\subhead 43.9\endsubhead
In the setup of 43.8 let $x\in\WW_I$. We show:
$$\align&\tr(\tT_{x\vp},E^v)=\sum_{E'\in\un\Irr(\tWW_I)}\la E',E\ra\tr(\tT_{x\vp},E'{}^v),\\&
\tr(x\vp,E)=\sum_{E'\in\un\Irr(\tWW_I)}\la E',E\ra\tr(x\vp,E');\tag a\endalign$$
here for any $E'$ in the sum,
$$\la E',E\ra=\dim_{\QQ(v)}\Hom_{\tH_I^v}(E'{}^v,E^v)=\dim_\QQ\Hom_{\tWW_I}(E',E).$$
Using 43.8(a) we can write the left hand side of the first equality in (a) as 
$$\sum_j\tr_{\bbq(v)}(\tT_{x\vp},\EE'_j).$$
 Here $\EE'_j$ is as in 43.8(a); if it is as in 43.8(i) then 
$\tr_{\bbq(v)}(\tT_{x\vp},\EE'_j)=0$ since $\G'\ne\G$. The contribution of the $j$ as in 43.8(ii) yields the right
hand side of the first equality in (a). The proof of the second equality in (a) is entirely similar.

We show:

(b) {\it If $E'$ in (a) satisfies $\la E',E\ra\ne0$ then $a_{E'}\le a_E$;}
\nl
(here $a_E$ is as in 43.6 and $a_{E'}$ is defined similarly in terms of $E',\tWW_I$). Indeed the simple
$\WW_I$-module $E'|_{\WW_I}$ appears in the $\WW_I$-module $E|_{\WW_I}$ hence (b) follows from 
\cite{\UNE, 20.14(a)}.

Let $\tH_I^\iy$ be defined like $\tH^\iy$ but for $\tWW_I$ instead of $\tWW$. For $x\in\WW_I$ we show:
$$\tr(t_{x\vp},E^\iy)=\sum_{E'\in\un\Irr(\tWW_I);a_{E'}=a_E}\la E',E\ra\tr(t_{x\vp},E'{}^\iy).\tag c$$
(The simple $\QQ\ot H_I^\iy$-module $E'{}^\iy$ is defined like $E^\iy$ but for $\tWW_I$ instead of $\tWW$.)  
We take the coefficient of $v^{-a_E}$ in both sides of the first equality in (a) (they are in
$\ca$; using 43.6(c) we
obtain 
$$\sgn(x)\tr(t_x\vp,E^\iy)=\sum_{E'}\la E',E\ra\tr(\tT_{x\vp},E'{}^v;-a_E)$$
where the sum over $E'$ is as in (a). By (b) the previous sum can be restricted to the $E'$ such that
$a_{E'}\le a_E$. The contribution of $E'$ with $a_{E'}<a_E$ is $0$ by 43.6(c) (for $\tWW_I$). Thus the sum can be 
restricted to the $E'$ such that $a_{E'}=a_E$. For such $E'$ we have, using again 43.6 (for $\tWW_I$):
$$\tr(\tT_{x\vp},E'{}^v;-a_E)=\tr(\tT_{x\vp},E'{}^v;-a_{E'})=\sgn(x)\tr(t_x\vp,E'{}^\iy)$$
and (c) follows.

\subhead 43.10\endsubhead
For any $E\in\Irr(\tWW)$ we define $\ph_E:\WW\vp@>>>\ZZ$ by $\ph_E(x\vp)=\tr(x\vp,E)$. Note that
$\ph_{E\ot\io}=-\ph_E$ ($\io$ as in 43.1). The functions $\ph_E$ with $E\in\Irr(\tWW)$ generate a subgroup 
$\car(\tWW)$ of the group of all functions $\WW\vp@>>>\ZZ$ which are constant on the orbits of the conjugation 
$\WW$-action on $\WW\vp$. From 43.5(c) we see that $\{\ph_E;E\in\fE\}$ is a $\ZZ$-basis of $\car(\tWW)$. For any 
$x\in\WW$ we set
$$\ale_{x\vp}=\sum_{E\in\un\Irr(\tWW)}
\fra{1}{2}\tr(t_x\vp,E^\iy)\ph_E=\sum_{E\in\fE}\tr(t_x\vp,E^\iy)\ph_E\in\car(\tWW).\tag a$$
From 43.6(e) we see that for any $E\in\fE$ we have:
$$\sum_{x\in\WW}\tr(t_x\vp,E^\iy)\ale_{x\vp}=f_E^\iy\dim(E)\ph_E\in\car(\tWW).\tag b$$  
Now let $I$ be a subset of $\II$ such that $\e(I)=I$. We define a homomorphism
$J_{\tWW_I}^{\tWW}:\car(\tWW_I)@>>>\car(\tWW)$ by 
$$J_{\tWW_I}^{\tWW}(\ph_{E'})=\sum_{E\in\un\Irr(\tWW);a_{E'}=a_E}\la E',E\ra\ph_E$$
for any $E'\in\Irr(\tWW_I)$. This is clearly well defined. For $x\in\WW_I$ we define 
$\ale^I_{x\vp}\in\car(\tWW_I)$ in the same way as $\ale_{x\vp}\in\car(\tWW)$ but in terms of $\tWW_I$ instead of 
$\tWW$. From 43.9(c) we see that
$$J_{\tWW_I}^{\tWW}(\ale^I_{x\vp})=\ale^I_{x\vp}.\tag c$$

\subhead 43.11\endsubhead
Let $I$ be a subset of $\II$ such that $\e(I)=I$. We fix a two-sided cell $\boc'$ of $\WW_I$ such that
$\e(\boc')=\boc'$. There is a unique two-sided cell $\boc$ of $\WW$ such that $\boc'\sub\boc$; we must have
$\e(\boc)=\boc$. We show:

(a) {\it if $E'\in\Irr(\tWW_I),E\in\Irr(\tWW)$ satisfy $\boc'=\boc_{E'}$ (see 43.6 with $\tWW$ replaced by 
$\tWW_I$) and $\la E',E\ra\ne0$, then $\boc_E\prq\boc$.}
\nl
To prove this we may replace $E,E'$ by their restrictions to $\WW,\WW_I$. Thus we may assume that
$\tWW=\WW,\tWW_I=\WW_I,\vp=1$. Since $\boc'=\boc_{E'}$, there exists $x\in\boc'$ such that the action of $t_x$ in
the $\QQ\ot\tH_I^\iy$-module $E'{}^\iy$ is $\ne0$. Using 43.6(a) we see that the action of $c_x^\da$ in the 
$H_I^v$-module $E'{}^v$ is $\ne0$. Since $\la E',E\ra\ne0$, $E'_v$ may be regarded as a $H_I^v$-submodule of 
$E^v$. Hence the action of $c_x^\da$ in the $H^v$-module $E^v$ is $\ne0$. Using 43.6(b) we see that $z\prq x$ for
some $z\in\boc_E$. By definition we have $x\in\boc$. This proves (a).

We show:

(b) {\it if $E'\in\Irr(\tWW_I),E\in\Irr(\tWW)$ satisfy $\boc'=\boc_{E'}$ (see 43.6 with $\tWW$ replaced by 
$\tWW_I$) and $a_{e'}=a_E$, $\la E',E\ra\ne0$, then $\boc=\boc_E$.}
\nl
Since the $\aa$-function of $\WW_I$ is known to be the restriction of the $\aa$-function of $\WW$, we see that
the value of the $\aa$-function on $\boc$ and $\boc_E$ coincide. Since $\boc_E\prq\boc$ (see (a)) it follows that
$\boc=\boc_E$.

\subhead 43.12\endsubhead
Let $x\in\WW$. Let $\boc$ be the two-sided cell containing $x$. According to \cite{\ORA, (5.3.1)} there exists 
uniquely defined elements $a_{y,x}\in\QQ(v)$ (for $y\in\WW, y\prec x$) such that

$(-1)^{l(x)}c_x^\da-\sum_{y;y\prec x}(-1)^{l(y)}a_{y,x}c_y^\da$ acts as zero on $E_0^v$ for any $E_0\in\Irr(\WW)$
with $\boc_{E_0}\ne\boc$.
\nl
Moreover for $y\prec x$ we have 
$$a_{y,x}=\sum_{j\in\ZZ_{>0}}a_{y,x;j}v^j$$
where  $a_{y,x;j}\in\ZZ$ for all $j$ and $a_{y,x;j}=0$ unless $j=l(x)+l(y)\mod2$, see \cite{\ORA, (5.3.6)}. It 
follows that the sum
$$\sum_{E\in\un\Irr(\tWW)}\fra{1}{2}
\tr(c_{x\vp}^\da-\sum_{y;y\prec x}(-1)^{-l(x)+l(y)}a_{y,x}c_{y\vp}^\da,E^v)\ph_E\in\car(\tWW)\tag a$$
is equal to the same sum restricted to those $E$ such that $\boc_E=\boc$. For such $E$ we have $a_E=\aa(x)$ and
for any $y$ such that $y\prec x$, $a_{y,x}\tr(c_{y\vp}^\da,E^v)$ is of the form $v^{-\aa(x)+1}$ times a rational 
function in $v$ which is regular at $v=0$; moreover, $\tr(c_{x\vp}^\da,E^v)$ is of the form 
$v^{-\aa(x)}\tr(t_{x\vp},E^\iy)$ plus higher powers of $v$. Thus (a) is of the form 
$$\sum_{E\in\un\Irr(\tWW);\boc_E=\boc}\fra{1}{2}v^{-\aa(x)}\tr(t_x\vp,E^\iy)\ph_E+\s$$
where $\s$ is a linear combination of elements $\ph_E$ with coefficients of the form $v^{-\aa(x)+1}$ times a 
rational function in $v$ which is regular at $v=0$. In the previous sum the condition $\boc_E=\boc$  can be 
dropped and the sum is unchanged. We see that (a) is equal to $v^{-\aa(x)}\ale_{x\vp}+\s$ with $\s$ as above. 
Taking in this identity coefficients of $v^{-\aa(x)}$ in the expansions at $v=0$ we obtain
$$\align&\ale_{x\vp}=\sum_{E\in\un\Irr(\tWW)}\fra{1}{2}(\tr(c_{x\vp}^\da,E^v;-\aa(x))
\\&-\sum_{y,j;y\prec x,j>0}(-1)^{-l(x)+l(y)}a_{y,x;j}\tr(c_{y\vp}^\da,E^v;-\aa(x)-j))\ph_E.\tag b\endalign$$

\head 44. Unipotent character sheaves and two-sided cells\endhead
\subhead 44.1\endsubhead
In this section we study the unipotent character sheaves in connection with Weyl group representations and 
two-sided cells. A number of results in this section are conditional (they depend on a cleanness property and/or
on a parity property); they will become unconditional in \S46.

The following convention will be used in this section. In parts of 44.3-44.7, marked as $\sp...\sp$, we assume
that the ground field $\kk$ is an algebraic closure of $\FF_q$ and we fix an $\FF_q$-structure on $G$ with 
Frobenius map $F:G@>>>G$ which leaves $B^*,T$ (see 28.5) stable and induces the identity map on $\WW$ and on 
$G/G^0$; we will view the various varieties which appear with the natural $\FF_q$-structure induced by that of 
$G$. The results in other parts of this section are valid for a general $\kk$ (by a standard reduction to the case
$\kk=\bar\FF_q$).

If $X$ is an algebraic variety with a given $\FF_q$-structure we write $\cd_m(X)$ for the corresponding mixed 
derived category of $\bbq$-sheaves. If $A\in\cd_m(X)$ is perverse and $j\in\ZZ$, we denote by $A_j$ the canonical
subquotient of $A$ wich is pure of weight $j$. 

\subhead 44.2\endsubhead
For any $w\in\WW$ let 

$Z^w_{\em,\II,D}=\{(B,B',x)\in\cb\T\cb\T D;xBx\i=B',\po(B,B')=w\}$
\nl
(see 28.8),

$\bZ^w_{\em,\II,D}=\{(B,B',x)\in\cb\T\cb\T D;xBx\i=B',\po(B,B')\le w\}$;
\nl
note that $\bZ^w_{\em,\II,D}=\sqc_{w'\in\WW;w'\le w}Z^{w'}_{\em,\II,D}$. Let 

$\cb^w=\{(B,B')\in\cb\T\cb;\po(B,B')=w\}$,
 
$\bar\cb^w=\{(B,B')\in\cb\T\cb;\po(B,B')\le w\}$.
\nl
Define $\mu:\bZ^w_{\em,\II,D}@>>>\bar\cb^w$ by
$\mu(B,B',x)=(B,B')$. Note that $\mu$ is a fibration with connected smooth fibres and 
$Z^{w'}_{\em,\II,D}=\mu\i(\cb^{w'})$ for any $w'\le w$. Hence $Z^w_{\em,\II,D}$ is an irreducible smooth open 
dense subvariety of $\bZ^w_{\em,\II,D}$. Let $\bbq^w$ be the local system $\bbq$ on $\cb^w$ and let 
$\bbq^{w\sh}=IC(\bar\cb^w,\bbq^w)\in\cd(\bar\cb^w)$. Let $\dbbq^w$ be the local system $\bbq$ on 
$Z^w_{\em,\II,D}$ and let 
$$\dbbq^{w\sh}=IC(\bZ^w_{\em,\II,D},\dbbq^w)=\mu^*\bbq^{w\ch}\in\cd(\bZ^w_{\em,\II,D}).$$

\subhead 44.3\endsubhead
$\sp$ For $y,w\in\WW$, $y\le w$ and $i\in\ZZ$ let $n_{y,w,i}$ be as in 43.2; by \cite{\KLL},

(a) {\it $\ch^i(\bbq^{w\sh})|_{\cb^y}$ is a local system isomorphic to $(\bbq^y)^{\op n_{y,w,i}}$; moreover it 
admits a filtration (over $\FF_q$) with $n_{y,w,i}$ steps and each subquotient isomorphic over $\FF_q$ to 
$\bbq(-i/2)$.}
\nl
Using the fibration $\mu$ we deduce that 

(b) {\it $\ch^i(\dbbq^{w\sh})|_{Z^y_{\em,\II,D}}$ is a local system isomorphic to $(\dbbq^y)^{\op n_{y,w,i}}$; 
moreover, it admits a filtration (over $\FF_q$) with $n_{y,w,i}$ steps and each subquotient isomorphic over 
$\FF_q$ to $\bbq(-i/2)$.}
\nl
Define $\p_w:Z^w_{\em,\II,D}@>>>D$, $\bpi_w:\bZ^w_{\em,\II,D}@>>>D$ by $(B,B',x)\m x$. Let 
$$K^w_D=\p_{w!}\dbbq^w\in\cd(D),\qua\bK^w_D=\bpi_{!w}\dbbq^{w\sh}\in\cd(D).$$
(With notation of 28.12 we have $K^w_D=K^{w,\bbq}_{\II,D}$.) We view $\dbbq^w$ and $\dbbq^{w\sh}$ as objects of 
$\cd_m(Z^w_{\em,\II,D})$ and $\cd_m(\bZ^w_{\em,\II,D})$ such that Frobenius acts trivially on the stalk at any 
$\FF_q$-rational point of $Z^w_{\em,\II,D}$. Applying to them $\p_{w!}$ and $\bpi_{!w}$ we obtain objects 
$\uK^w_D\in\cd_m(D),\ubK^w_D\in\cd_m(D)$.

The following equality in the Grothendieck group of mixed perverse sheaves on $D$ is verified (using (b)) along 
the lines of \cite{\CS, 12.6}:
$$\sum_{i\in\ZZ}(-1)^iH^i(\ubK^w_D)=\sum_{y\in\WW;y\le w}\sum_{i,h\in\ZZ}(-1)^in_{y,w,h}H^i(\uK^y_D)(-h/2).\tag c
$$
We now take the part of weight $j$ in (c); note that $H^j(\ubK^w_D)$ is pure of weight $j$ since $\bpi_{w!}$ 
preserve weights and $\dbbq^{w\sh}$ is pure of weight $0$.) We see that for any $j\in\ZZ$, the following 
equality holds in the Grothendieck group of perverse sheaves on $D$:
$$(-1)^jH^j(\bK^w_D)=\sum_{y\in\WW;y\le w}\sum_{i,h\in\ZZ}(-1)^in_{y,w,h}H^i(\uK^y_D)_{j-h}. \sp\tag d$$

\subhead 44.4\endsubhead
We shall often write $\hD^{un}$ instead of $\hD^{\bbq}$ (see 28.14). 

{\it Definition.} We say that a character sheaf $A$ on $D$ is {\it unipotent} if $A\in\hD^{un}$.
\nl
Let $\uhD^{un}$ be the set of isomorphism classes of unipotent character sheaves on $D$. The following two 
conditions on a simple perverse sheaf $A$ on $D$ are equivalent:

(i) $A\in\hD^{un}$;

(ii) $A\dsv\bK^w_D$ for some $w\in\WW$.
\nl
This follows from (a) below which is verified along the lines of \cite{\CS,III, (12.7.1)}. 

(a) Let $w\in\WW$ be such that $A\not\dsv K^y_D$ for any $y\in\WW,y<w$. Then $(A:H^i(\bK^w_D))=(A:H^i(K^w_D))$ for
any $i\in\ZZ$. 
\nl
Let $\Xi$ be a set of representatives for the isomorphism classes of 
objects in $\hD^{un}$; note that $\Xi$ is a finite set. 

\subhead 44.5\endsubhead
Let $A\in\hD^{un}$. We regard $H\tT_\vp$ as an ideal in $\tH$. Let $\z^A_0:H\tT_\vp@>>>\ca$ be the composition 
of the map $H\tT_\vp@>>>H_1\tT_{\uD}$ (restriction of the natural surjection $\tH@>>>H^D_1$) with the map 
$\z_A:H_1\tT_{\uD}@>>>\ca$ in 31.7 (with $n=1$). From the definitions, $\z_A^0$ is an $\ca$-linear map $\sp$ and 
for any $x\in\WW$ we have
$$\z^A_0(v^{l(x)}\tT_{x\vp})=v^{-\dim G}\sum_{i,j}(-1)^i(A:H^i(\uK^x_D)_j)v^j.\sp\tag a$$
For $x\in\WW$ we show:
$$\z^A_0(c_x\tT_\vp)=v^{-\dim G-l(x)}\sum_{j\in\ZZ}(A:H^j(\bK^x_D))(-v)^j.\tag b$$
$\sp$ By 44.3(d) we have for any $j$:
$$(-1)^j(A:H^j(\bK^x_D))=\sum_{y\in\WW;y\le x}\sum_{i,h\in\ZZ}(-1)^in_{y,x,h}(A:H^i(\uK^y_D)_{j-h}).$$
We deduce
$$\align&v^{-\dim G-l(x)}\sum_{j\in\ZZ}(A:H^j(\bK^x_D))(-v)^j\\&=
v^{-\dim G-l(x)}\sum_{y\in\WW;y\le x}\sum_{i,j,h\in\ZZ}(-1)^in_{y,x,h}(A:H^i(K^y_D)_{j-h})v^j\\&
=v^{-\dim G-l(x)}\sum_{y\in\WW;y\le x}\sum_{i,j',h\in\ZZ}(-1)^in_{y,x,h}(A:H^i(K^y_D)_{j'})v^{j'+h}.\sp\tag c
\endalign$$
We can rewrite this as
$$\align&v^{-l(x)}\sum_{y\in\WW;y\le x}\sum_{h\in\ZZ}n_{y,x,h}v^h\z^A(v^{l(y)}\tT_{y\vp})\\&=
v^{-l(x)}\sum_{y\in\WW;y\le x}P_{y,x}(v^2)\z^A(v^{l(y)}\tT_{y\vp})=\z^A_0(c_x\tT_\vp).\endalign$$
This proves (b).

\subhead 44.6\endsubhead
Let $\ck^{un}(D)$ be the subgroup of the Grothendieck group of the category of perverse sheaves on $D$ generated 
by the objects in $\hD^{un}$. Let $\ck^{un}_\QQ(D)=\QQ\ot\ck^{un}(D)$. Let $(:)$ be the symmetric $\QQ$-bilinear 
form on $\ck^{un}_\QQ(D)$ with values in $\QQ$ such that $(A:A)=1$ if $A\in\hD^{un}$ and $(A:A')=0$ if 
$A,A'\in\hD^{un}$ are not isomorphic. Note that if $P$ is a perverse sheaf on $D$ all of whose simple 
subquotients are in $\hD^{un}$ then the present meaning of $(A:P)$ agrees with the earlier meaning, see 31.6. 

For any $x\in\WW$ we show:
$$gr_1(\bK^x_D)=\sum_{y\in\WW;y\le x}P_{y,x}(1)gr_1(K^y_D)\in\ck^{un}(D).\tag a$$
$\sp$ Specializing 44.5(c) for $v=1$ $\sp$ we deduce
$$gr_1(\bK^x_D)=\sum_{y\in\WW;y\le x}\sum_{j',h\in\ZZ}n_{y,x,h}gr_1((\uK^y_D)_{j'})\in\ck^{un}(D)$$
and (a) follows.

For any $E\in\Mod(\tWW)$ we set
$$R_E=|\WW|\i\sum_{x\in\WW}(-1)^{\dim G}\tr(x\vp,E)gr_1(K^x_D)\tag b$$
(an element of $\ck^{un}_\QQ(D)$). We show:
$$R_E=|\WW|\i\sum_{x\in\WW}(-1)^{\dim G}\tr(\tc_{x\vp}|_{v=1},E)gr_1(\bK^x_D)\tag c$$
where $\tc_{x\vp}$ is as in 43.2. We shall use the known inversion formula
$$\sum_{z\in\WW;y\le z\le x}(-1)^{l(y)-l(z)}P_{y,z}(\qq)P_{w_0x,w_0z}(\qq)=\d_{y,x}\tag d$$
for any $y\le x$ in $\WW$. Using (a),(d) and the definition of $\tc_{x\vp}$, we see that the right hand side of 
(c) is
$$\align&|\WW|\i\sum_{x,y,z\in\WW;y\le x\le z}(-1)^{\dim G}(-1)^{l(z)-l(x)}P_{y,x}(1)P_{w_0z,w_0x}(1)
\tr(z\vp,E)gr_1(K^y_D)\\&=|\WW|\i\sum_{y\in\WW}(-1)^{\dim G}\tr(y\vp,E)gr_1(K^y_D)=R_E,\endalign$$
as required. 

Let $\Mod_{\bbq}(\tWW)$ be the category of $\bbq[\tWW]$-modules of finite dimension over $\bbq$. For 
$E\in\Mod_{\bbq}(\tWW)$ we define $R_E\in\bbq\ot\ck^{un}(D)$ by the same formula as (b).

For any $\ph\in\car(\tWW)$ (see 43.10) we define $R_\ph\in\ck^{un}_\QQ(D)$ by $R_\ph=\sum_{E\in\fE}p_ER_E$ where 
$\ph=\sum_{E\in\fE}p_E\ph_E$ ($p_E\in\ZZ$). This is independent of the choice of $\fE$ since $R_{E\ot\io}=-R_E$ 
for $E\in\Irr(\tWW)$. Note that for $E\in\Irr(\tWW)$ we have $R_{\ph_E}=R_E$.

\subhead 44.7\endsubhead
Let $A\in\hD^{un}$. For any $E\in\Irr(\tWW)$ we set
$$b_{A,E}^v=\fra{1}{f_E^v\dim E}\sum_{x\in\WW}\z^A_0(\tT_{x\vp})\tr(\tT_{x\vp},E^v)\in\QQ(v).\tag a$$
Note that this definition is compatible with that in 34.19(b). Using 34.19(a) we see that for any $\x\in H$ we 
have
$$\z^A_0(\x\tT_\vp)=\sum_{E\in\fE}b^v_{A,E}\tr(\x\tT_\vp,E^v).\tag b$$
Taking here $\x=c_x,x\in\WW$ and using 44.5(b), we deduce:
$$\sum_{j\in\ZZ}(A:H^j(\bK^x_D))(-v)^j=v^{\dim G+l(x)}\sum_{E\in\fE}b^v_{A,E}\tr(c_x\tT_\vp,E^v).\tag c$$

Let $\hD^{unc}$ be the subcategory of $\hD^{un}$ whose objects are the unipotent character sheaves on $D$ which 
are cuspidal. 

An object $A\in\hD^{unc}$ is said to be {\it clean} if the following condition is satisfied: $A|_{\bS-S}=0$ where
$S$ is the isolated stratum of $D$ such that $\supp(A)$ is the closure $\bS$ of $S$.

We say that $D$ has property $\fA_0$ if any $A\in\hD^{unc}$ is clean. We say that $D$ has property $\fA$ if for 
any parabolic subgroup $P$ of $G^0$ such that $N_DP\ne\em$, the connected component $N_DP/U_P$ of $N_GP/U_P$ has 
property $\fA_0$. (Compare 33.4(b).) 

We say that $D$ has property $\tfA$ if for any $A\in\hD^{un}$ and any $w\in\WW$, $i\in\ZZ$ such that 
$(A:H^i(\bK^w_D))\ne0$ we have $i=\dim\supp(A)\mod2$.

{\it In the remainder of this section we assume that $D$ has property $\fA$.}

Using 35.18(g) we see that for any $E,E'$ in $\fE$ we have
$$\sum_{A'\in\Xi}b_{A',E}^vb_{A',E'}^v=\d_{E,E'}.\tag d$$
Let $A\in\hD^{un}$. Using 35.22 we see that for any $E\in\Irr(\tWW)$ we have
$$b_{A,E}^v\in\QQ.\tag e$$
(The quasi-rationality assumption in 35.22 is automatically satisfied in our case; see 43.4(a).)
In view of (e) we shall write $b_{A,E}$ instead of $b_{A,E}^v$. We show:
$$b_{A,E}=(-1)^{\dim G}(A:R_E).\tag f$$
Let $x\in\WW$. $\sp$ Setting $v=1$ in 44.5(a) we obtain
$$\z^A_0(\tT_{x\vp})|_{v=1}=\sum_{i,j}(-1)^i(A:H^i(\uK^x_D)_j)=(A:gr_1(K^x_D)).\sp\tag g$$
Setting $v=1$ in (b) with $\x=\tT_x$ and using (e) we obtain
$$\z^A_0(\tT_{x\vp})|_{v=1}=\sum_{E\in\fE}b_{A,E}\tr(x\vp,E).$$
Combining with (g) we obtain
$$(A:gr_1(K^x_D))=\sum_{E\in\fE}b_{A,E}\tr(x\vp,E).\tag h$$
Using the orthogonality relations 43.5(b) specialized for $v=1$ we obtain
$$b_{A,E}=|\WW|\i\sum_{x\in\WW}\tr(x\vp,E)(A:gr_1(K^x_D))$$
for any $E\in\fE$. This proves (f) in the case where $E\in\fE$. This clearly implies (f) in the general case.

We can now rewrite (h) as
$$gr_1(K^x_D)=(-1)^{\dim G}\sum_{E\in\fE}\tr(x\vp,E)R_E\tag i$$
in $\ck^{un}_\QQ(D)$ and (c) as:
$$\sum_{j\in\ZZ}(A:H^j(\bK^x_D))(-v)^j=(-1)^{\dim G}v^{\dim G+l(x)}\sum_{E\in\fE}(A:R_E)\tr(c_x\tT_\vp,E^v).\tag j
$$
We show:

(k) {\it there exists $E\in\fE$ such that $(A:R_E)\ne0$.}
\nl
We can find $x\in\WW$ and $j\in\ZZ$ such that $(A:H^j(\bK^x_D))\ne0$. Then the left hand side of (j) is $\ne0$
hence so is the right side. Thus (k) holds.

We show:

(l) {\it For $E,E'\in\Mod_{\bbq}(\tWW)$ we have 
$$(R_E:R_{E'})=|\WW|\i\sum_{x\in\WW}\tr(x\vp,E)\tr(x\vp,E').$$
Moreover, if $E,E'\in\fE$ then we have $(R_E:R_{E'})=\d_{E,E'}$.}
\nl
Here $(:)$ is the bilinear form $\bbq\ot\ck^{un}(D)\T\bbq\ot\ck^{un}(D)@>>>\bbq$ extending $(:)$ in 44.6.

Assume first that $E,E'\in\fE$. Clearly, $R_E=\sum_{A'\in\Xi}(A':R_E)A'$, 
$R_{E'}=\sum_{A'\in\Xi}(A':R_{E'})A'$. It follows that 
$$(R_E:R_{E'})=\sum_{A'\in\Xi}(A':R_E)(A':R_{E'})=\sum_{A'\in\Xi}b_{A',E}b_{A',E'}=\d_{E,E'}$$
where the last two equalities come from (f),(d). This proves the second equality in (l). To prove the first
equality in (l) we may assume that $E,E'$ are simple objects of $\Mod_{\bbq}(\tWW)$. If the restriction of $E$ to
$\bbq[\WW]$ is not simple then $\tr(x\vp,E)=0$ for any $x\in\WW$ hence both sides of the first equality in (l) are
$0$. Thus we may assume in addition that $E|_{\bbq[\WW]}$ is simple; similarly we may assume that
$E'|_{\bbq[\WW]}$ is simple. Replacing $E,E'$ by their tensor products with one dimensional representations of 
$\tWW$ which are trivial on $\WW$ reduces us to the case where $E,E'$ come from objects of $\fE$ by extension of
scalars. Using then the second identity in (l) we see that it is enough to show that for $E,E'\in\fE$ we have
$|\WW|\i\sum_{x\in\WW}\tr(x\vp,E)\tr(x\vp,E')=\d_{E,E'}$. But this is known from 43.5(c). This completes the proof
of (l).

For any $x\in\WW,i\in\ZZ$ we take the coefficient of $v^{i+l(x)+\dim G}$ in the two sides of (j); we obtain
$$(-1)^{i+l(x)}(A:H^{i+l(x)+\dim G}(\bK^x_D))=\sum_{E\in\un\Irr(\tWW)}
\fra{1}{2}\tr(c_x\tT_\vp,E^v;i)(A:R_E).\tag m$$
For any $y,z$ in $\WW$ we show
$$gr_1(K^{y\i z\vp y\vp\i}_D)=gr_1(K^z_D).\tag n$$
Using (i) this is the same as
$$\sum_{E\in\fE}\tr(y\i z\vp y,E)R_E=\sum_{E\in\fE}\tr(z\vp,E)R_E$$
which is clear since $\tr(y\i z\vp y,E)=\tr(z\vp,E)$ for any $E\in\fE$.

We show:

(o) {\it if $E\in\Mod(\tWW)$ then $R_E$ is a $\ZZ$-linear combination of elements $R_{E_1}$ with 
$E_1\in\Irr(\tWW)$.}
\nl
We can write $\bbq\ot E=\op_h\EE_h$ where $\EE_h$ are simple $\bbq[\tWW]$-modules. Hence 
$R_E=R_{\bbq\ot E}=\sum_hR_{\EE_h}$. If $h$ is such that $\EE_h|_\WW$ is not a simple $\bbq[\WW]$-module then
$\tr(x\vp,\EE_h)=0$ for any $x\in\WW$ hence $R_{\EE_h}=0$. If $h$ is such that $\EE_h|_\WW$ is a simple 
$\bbq[\WW]$-module then by taking the tensor products of $\EE_h$ with a one dimensional representation of $\tWW$ 
which is trivial on $\WW$ we obtain a module which comes from an object of $\Irr(\tWW)$. It follows that 
$R_E=\sum_{E_1\in\fE}c_{E_1}R_{E_1}$ where $c_{E_1}$ are integer combination of roots of $1$. 
Using (l) we have $c_{E_1}=(R_E:R_{E_1})=|\WW|\i\sum_{x\in\WW}\tr(x\vp,E)\tr(x\vp,E_1)$. This is a rational
number. Being also an algebraic integer it is an integer. This proves (o).

We show:

(p) {\it for $E\in\fE$, $x\in\WW$ we have $(R_E:gr_1(K^x_D))=(-1)^{\dim G}\tr(x\vp,E)$.}
\nl
Using (i) we have $(R_E:gr_1(K^x_D))=(R_E:(-1)^{\dim G}\sum_{E'\in\fE}\tr(x\vp,E')R_{E'})$ so that (p) follows 
from (l).

\subhead 44.8\endsubhead
The $\ca$-linear involution $\dd:\fK(D)@>>>\fK(D)$ in 42.2 induces (by the specialization $v=1$) a $\ZZ$-linear 
involution $\dd:\ck(D)@>>>\ck(D)$ ($\ck(D)$ as in 38.9). By extension of scalars, $\dd$ gives rise to a 
$\QQ$-linear involution $\QQ\ot\ck(D)@>>>\QQ\ot\ck(D)$ denoted again by $\dd$. 

Let $A\in\hD$. We show that:
$$\dd(A)=(-1)^{\codim(\supp(A))}A^\circ\tag a$$
where $A^\circ\in\hD$. We can find a parabolic $P_0$ of $G^0$ such that $N_DP_0\ne\em$ and a cuspidal character
sheaf $A_0$ on $D_0:=N_DP_0/U_{P_0}$ such that $A$ is a direct summand of $\ind_{D_0}^D(A_0)$. We have 
$P_0\in\cp_J$ where $J\sub\II,\e(J)=J$. By 38.11(a) we have $\dd(A)=(-1)^{|J_\e|}A^\circ$ where $A^\circ\in\hD$ 
and $J_\e$ is the set of orbits of $\e:J@>>>J$. It remains to show that $\codim(\supp(A))=|J_\e|\mod2$. From the 
theory of admissible complexes (6.7) and from 3.13(b) we see that 
$\dim\supp(A)=\dim G^0-\dim(P_0/U_{P_0})+\dim\supp(A_0)$ that is, $\codim(\supp(A))=\codim(\supp(A_0))$. Also the
analogue of $J_\e$ for $A_0$ is $J_\e$ itself. Thus we are reduced to the case where $A=A_0$ that is, we may 
assume that $A$ is cuspidal. Let $G'={}^D\cz_{G^0}^0\bsl G,D'={}^D\cz_{G^0}^0\bsl D$. Then the support of $A$ is 
the closure of a subset of $D$ which is the inverse image of a single $G'{}^0$-conjugacy class $C$ in $D'$ under 
the obvious map $D@>>>D'$. Moreover, ${}^{D'}\cz_{G'{}^0}^0=\{1\}$. The set $\II$ for $G'$ can be identified with
that for $G$. Since $\codim(\supp(A))=\codim_{D'}C$, it is enough to show that $\codim_{D'}C=|\II_\e|\mod2$ for 
any $G'{}^0$-conjugacy class $C$ in $D'$. According to Spaltenstein \cite{\SP} we have $\codim_{D'}C=2\b+r$ where
$\b$ is the dimension of the variety of Borel subgroups of $G'{}^0$ that are normalized by some fixed element of 
$C$ and $r$ is the rank of the connected centralizer in $G'$ of any quasisemisimple element of $D'$. Thus,
$\codim_{D'}C=r\mod2$. It remains to note that $r=|\II_\e|$. 

By 42.9 (specialized with $v=1$) we see that for any $x\in\WW$ we have
$$\dd(\sum_{i\in\ZZ}H^i(K^x_D))=(-1)^{l(x)}\sum_{i\in\ZZ}H^i(K^x_D)\tag b$$
in $\ck(D)$. Here $H^i(K^x_D)$ is identified with the element $\sum_{A'\in\Xi}(A':H^i(K^x_D))A'$ of $\ck(D)$. We
show that for any $E\in\Irr(\tWW)$ we have
$$\dd(R_E)=R_{E\ot\sgn}.\tag c$$
Indeed, by (b), this is the same as the obvious equality
$$\align&|\WW|\i\sum_{i\in\ZZ}\sum_{x\in\WW}(-1)^{i+\dim G+l(x)}\tr(x\vp,E)H^i(K^x_D)\\&
=|\WW|\i\sum_{i\in\ZZ}\sum_{x\in\WW}(-1)^{i+\dim G}\tr(x\vp,E\ot\sgn)H^i(K^x_D).\endalign$$
If $A\in\hD^{un}$ then, by 44.7(k), there exists $E\in\Irr(\tWW)$ such that the coefficient of $A$ in $R_E$ is 
$\ne0$. Applying $\dd$ to $R_E$ we see that the coefficient of $A^\circ$ in $\dd(R_E)$ is $\ne0$ that is, the 
coefficient of $A^\circ$ in $R_{E\ot\sgn}$ is $\ne0$. In particular, $A^\circ\in\hD^{un}$. In the same way we see
that for any $E\in\Irr(\tWW)$ we have 
$$(A:R_E)=\pm(A^\circ:R_{E\ot\sgn}).\tag d$$
Using (a) and the equality $\dd\dd=1$ we obtain
$$A=(-1)^{\codim(\supp(A))}\dd(A^\circ)=(-1)^{\codim(\supp(A))}(-1)^{\codim(\supp(A^\circ))}(A^\circ)^\circ.$$
It follows that $(A^\circ)^\circ\cong A$ and 
$$\codim(\supp(A))=\codim(\supp(A^\circ))\mod2.\tag e$$

\subhead 44.9\endsubhead
For any sequence $\ss=(s_1,s_2,\do,s_r)$ in $\II$ we write $K^\ss_D,\bK^\ss_D$ instead of $K^{\ss,\bbq}_{\II,D}$,
$\bK^{\ss,\bbq}_{\II,D}$, see 28.12.

Let $A\in\hD^{un}$. Then $(A:H^i(\bK^w_D))\ne0$ for some $w\in\WW,i\in\ZZ$. We set 
$$\ee^A=(-1)^{i+\dim G}.\tag a$$
We show that $\ee^A$ is well defined. Assume that we have also $(A:H^{i'}(\bK^{w'}_D))\ne0$ with $w'\in\WW$,
$i'\in\ZZ$. We must show that $i=i'\mod2$. Let $\ss=(s_1,s_2,\do,s_r)$, $\ss'=(s'_1,s'_2,\do,s'_{r'})$ be 
sequences in $\II$ such that $s_1s_2\do s_r=w$, $s'_1s'_2\do s'_{r'}=w'$, $r=l(w)$, $r'=l(w')$. We will show that 

(b) $\bK^w_D$ is a direct summand of $\bK^\ss_D$.
\nl
Assuming this and the similar statement for $w',\ss'$ instead of $w,\ss$ we see that $(A:H^i(\bK^\ss_D))\ne0$ and
$(A:H^{i'}(\bK^{\ss'}_D))\ne0$ and the congruence $i=i'\mod2$ follows from 35.17(a). (Although in 35.17 it is 
assumed that $D$ is clean, in the present application it is enough to use the weaker hypothesis that $\fA$ holds
for $D$.)

Recall that $\bK^\ss_D=\bpi_\ss!\bbq$ where 
$$\align&\bZ^\ss_{\em,\II,D}=\{(B_0,B_1,\do,B_r,g)\in\cb^{r+1}\T D;gB_0g\i=B_r,\po(B_{i-1},B_i)\in\{1,s_i\}
\\&\text{ for }i\in[1,r]\}\endalign$$
and $\bpi_\ss:\bZ^\ss_{\em,\II,D}@>>>D$ is given by $(B_0,B_1,\do,B_r,g)\m g$. Recall from 44.2 that 
$\bK^w_D=\bpi_{w!}\dbbq^{w\sh}$. We have $\bpi_\ss=\bpi_w\r$ where $\r:\bZ^\ss_{\em,\II,D}@>>>\bZ^w_{\em,\II,D}$ 
is given by $(B_0,B_1,\do,B_r,g)\m(B_0,B_r,g)$. Hence $\bK^\ss_D=\bp_{w!}(\r_!\bbq)$ so that to prove (b) it is 
enough to show that $\dbbq^{w\sh}$ is a direct summand of $\r_!\bbq$. This follows from the fact that $\r$ is 
proper and is an isomorphism over an open dense subset of $\bZ^w_{\em,\II,D}$. This proves (b).

\subhead 44.10\endsubhead
We now fix a subset $I\sub\II$ such that $\e(I)=I$. Let $P\in\cp_I$ (see 26.1). Then $N_DP\ne\em$ so that 
$D':=N_DP/U_P$ is a connected component of the reductive group $G':=N_GP/U_P$; note that $G'{}^0=P/U_P$. Let 
$\p':N_DP@>>>D'$ be the obvious map. As in 27.1 we consider the diagram $D'@<\ua<<V_1@>a'>>V_2@>a''>>D$ where 
$V_1=\{(g,x)\in D\T G^0;x\i gx\in N_DP\}$, $V_2=\{(g,xP)\in D\T G^0/P;x\i gx\in N_DP\}$, $\ua(g,x)=\p'(x\i gx)$, 
$a'(g,x)=(g,xP)$, $a''(g,xP)=g$. As in 27.1 for any $G'{}^0$-equivariant perverse sheaf $A'$ we define a complex
of sheaves $A=\ind_{D'}^D(A')\in\cd(D)$ by $A=a''_!A'_1[2\dim U_P]$ where $A'_1\in\cd(V_2)$ is such that 
$\ua^*A'=a'{}^*A'_1$. We show:

(a) {\it if $A'\in\hD'{}^{un}$ then $\ind_{D'}^D(A')$ is isomorphic to a direct sum of objects of $\hD^{un}$.}
\nl
The proof is similar to that of \cite{\CS, I, 4.8}. Before giving it we need some preliminaries. Let $\cb'$ be the
flag manifold of $G'{}^0=P/U_P$. For $\b\in\cb'$ let $\ti\b\in\cb$ be the inverse image of $\b$ under the obvious
map $P@>>>G'{}^0$. Let $w\in\WW_I$ (see 26.1). Recall that 
$$\bZ^w_{\em,\II,D}=\{(B,B',x)\in\cb\T\cb\T D;xBx\i=B',\po(B,B')\le w\}.$$
Replacing here $D,\II$ by $D',I$ we have
$$\bZ^w_{\em,I,D'}=\{(\b,\b',y)\in\cb'\T\cb'\T D';y\b y\i=\b',\po(\b,\b')\le w\}.$$
We have a commutative diagram with cartesian squares
$$\CD
\bZ^w_{\em,I,D'}@<\tua<<\tV_1@>\ta'>>\bZ^w_{\em,\II,D}@.{}\\
@V\d VV         @V\d'VV    @V\d''VV     @.\\
D'@<\ua<<V_1@>a'>>V_2@>a''>>D
\endCD$$
where 
$$\align&\tV_1=\{(\b,b',y,g,x)\in\cb'\T\cb'\T D';y\b y\i=\b',x\i gx\in N_DP,\\&
y=\p'(x\i gx),\po(\b,\b')\le w\},\endalign$$
$\tua(\b,\b',y,g,x)=(\b,\b',y)$, $\ta'(\b,\b',y,g,x)=(x\ti\b x\i,x\ti\b'x\i,g)$,    

$\d(\b,\b',y)=y$, $\d'(\b,\b',y,g,x)=(g,x)$, $\d''(B,B',x)=(x,zP)$ with $z\in G^0$ such that $z\i Bz\sub P$.
\nl
Note that $\ua,\tua$ are smooth with connected fibres and $\ta',a'$ are principal $P$-bundles. It follows that
$$IC(\tV_1,\bbq)=\tua^*IC(\bZ^w_{\em,I,D'},\bbq)=\ta'{}^*IC(\bZ^w_{\em,\II,D},\bbq)$$
where the first $\bbq$ lives on 
$$\{(\b,b',y,g,x)\in\tV_1;\po(\b,\b')=w\}=\tua\i(Z^w_{\em,I,D'})=\ta'{}\i(Z^w_{\em,\II,D}),$$
the second $\bbq$ lives on $Z^w_{\em,I,D'}$ and the third $\bbq$ lives on $Z^w_{\em,\II,D}$. Hence
$$\d'_!IC(\tV_1,\bbq)=\ua^*\d_!IC(\bZ^w_{\em,I,D'},\bbq)=\a'{}^*\d''_!IC(\bZ^w_{\em,\II,D},\bbq)$$
that is, $\d'_!IC(\tV_1,\bbq)=\ua^*\bK^w_{D'}=\a'{}^*K'$ where $K'=\d''_!IC(\bZ^w_{\em,\II,D},\bbq)\in\cd(V_2)$.
Since $\ua$, $a'$ are smooth with connected fibres of dimension $\dim D+\dim U_P$, $\dim D-\dim U_P$ respectively
we see that for any $i$ we have 
$$\align&\ua^*(H^{i-\dim D-\dim U_P}\bK^w_{D'})[\dim D+\dim U_P]=H^i(\ua^*\bK^w_{D'})\\&
=H^i(a'{}^*\bK^w_D)=a'{}^*(H^{i-\dim D+\dim U_P}K')[\dim D-\dim U_P]\endalign$$
hence (seting $j=i-\dim D-\dim U_P$):
$$\ua^*(H^j\bK^w_{D'})=a'{}^*(H^{j+2\dim U_P}K')[-2\dim U_P].$$
We see that
$$\ind_{D'}^D(H^j\bK^w_{D'})=a''_!(H^{j+2\dim U_P}K').$$
We have 
$$\op_j\ind_{D'}^D(H^j\bK^w_{D'})[-j]=\op_j(H^{j+2\dim U_P}\bK^w_D)[-j]\text{ in }\cd(D).\tag b$$
Indeed the left hand side is
$$\align&\op_ja''_!(H^{j+2\dim U_P}K')=a''_!K'[2\dim U_P]\\&=a''_!\d''_!IC(\bZ^w_{\em,\II,D},\bbq)[2\dim U_P]
=\bK^w_D[2\dim U_P];\endalign$$
(we have used that $K'\cong\op_jH^jK'[-j]$ which follows from the decomposition theorem \cite{\BBD} applied to the
proper map $\d''$). This is equal to the right hand side of (b) since $\bK^w_D\cong\op_jH^j(\bK^w_d)[-j]$, by the
decomposition theorem applied to the proper map $a''\d''$. Now $H^j\bK^w_{D'}$ is a direct sum of character 
sheaves on $D'$; hence, by 30.6(a), $\ind_{D'}^D(H^j\bK^w_{D'})$ is a perverse sheaf on $D$ for any $j$.
Taking $H^i$ for both sides of (b) we obtain for any $i\in\ZZ$:
$$\ind_{D'}^D(H^i\bK^w_{D'})=H^{i+2\dim U_P}\bK^w_D.\tag c$$
Now let $A'\in\hD'{}^{un}$. We can find $w\in\WW_I$ and $i\in\ZZ$ such that $A'$ appears in $H^i\bK^w_{D'}$. Since
$H^i\bK^w_{D'}$ is semisimple, $A'$ is a direct summand of $H^i\bK^w_{D'}$. Using (c) we see that 
$\ind_{D'}^D(A')$ is a direct summand of $H^{i+2\dim U_P}\bK^w_D$. Hence (a) holds.

From (a) we see that $A'\m\ind_{D'}^D(A')$ (with $A'\in\hD'{}^{un}$) defines a group homomorphism
$\ck^{un}(D')@>>>\ck^{un}(D)$ and a $\QQ$-linear map $\ck^{un}_\QQ(D')@>>>\ck^{un}_\QQ(D)$ denoted again by 
$\ind_{D'}^D$.

Applying this homomorphism to both sides of 44.6(a) for $D'$ instead of $D$ and for $x\in\WW_I$ and using (c) we 
obtain
$$gr_1(\bK^x_D)=\sum_{y\in\WW_I;y\le x}P_{y,x}(1)\ind_{D'}^D(gr_1(K^y_{D'})).$$ 
Here $P_{y,x}$ are as in 43.2 for $\WW_I$ or equivalently for $\WW$. The left hand side can be evaluated using 
44.3(d) for $D$; we obtain
$$\sum_{y\in\WW_I;y\le x}P_{y,x}(1)gr_1(K^y_D)=\sum_{y\in\WW_I;y\le x}P_{y,x}(1)\ind_{D'}^D(gr_1(K^y_{D'})).$$ 
Since the matrix $(P_{y,x})_{x,y\in\WW_I}$ is invertible, we deduce for any $y\in\WW_I$:
$$\ind_{D'}^D(gr_1(K^y_{D'}))=gr_1(K^y_D).\tag d$$

\subhead 44.11\endsubhead
We preserve the setup of 44.10. Let $\G,\tWW$ be as in 43.1 and let $\tWW_I$ be the subgroup of $\tWW$ generated 
by $\WW_I$ and $\G$; now $\tWW_I$ plays the same role for $\WW_I$ as $\tWW$ for $\WW$. For any 
$E'\in\Mod(\tWW_I)$, the element $R_{E'}\in\ck^{un}_\QQ(D')$ is defined as in 44.6(b). Let 
$\ind_{\tWW_I}^{\tWW}E'\in\Mod(\tWW)$ be the induced module. We show:
$$\ind_{D'}^D(R_{E'})=R_{\ind_{\tWW_I}^{\tWW}E'}\in\ck^{un}_\QQ(D).\tag a$$
Applying $\ind_{D'}^D$ to 44.6(b) with $E,D$ replaced by $E',D'$ and using 44.10(d) we obtain
$$\ind_{D'}^D(R_{E'})=|\WW_I|\i\sum_{i\in\ZZ}\sum_{x\in\WW_I}(-1)^{i+\dim G'}\tr(x\vp,E')H^i(K^x_D).$$
Using the definitions and 44.7(n) we have
$$\align&R_{\ind_{\tWW_I}^{\tWW}E'}=|\WW|\i\sum_{i\in\ZZ}\sum_{x\in\WW}(-1)^{i+\dim G}\tr(x\vp,ind E')H^i(K^x_D)
\\&=|\WW|\i|\WW_I|\i\sum_{i\in\ZZ}\sum_{x\in\WW,y\in\WW;yx\vp y\i\in\WW_I\vp}
(-1)^{i+\dim G}\tr(yx\vp y\i,E')H^i(K^x_D)\\&
=|\WW|\i|\WW_I|\i\sum_{i\in\ZZ}\sum_{z\in\WW_I,y\in\WW}(-1)^{i+\dim G}\tr(z\vp,E')H^i(K^{y\i z\vp y\vp\i}_D)\\&
=|\WW|\i|\WW_I|\i\sum_{i\in\ZZ}\sum_{z\in\WW_I,y\in\WW}(-1)^{i+\dim G}\tr(z\vp,E')H^i(K^z_D)\\&
=|\WW_I|\i\sum_{i\in\ZZ}\sum_{z\in\WW_I}(-1)^{i+\dim G}\tr(z\vp,E')H^i(K^z_D).\endalign$$
Now (a) follows since $\dim G=\dim G'\mod2$.

\subhead 44.12\endsubhead
We preserve the setup of 44.10. Let $\ss$ be a sequence in $\II$. From 29.14 we see that 
$\res_D^{D'}(\bK^\ss_D)\cong\op_{\tt\in\ct}\bK^\tt_{D'}[-d_\tt]$ where $\ct$ is a certain finite collection of 
sequences in $I$ and $d_\tt$ are integers. Since $\bK^\ss_D\cong\op_iH^i(\bK^\ss_D)[-i]$, 
$\bK^\tt_{D'}\cong\op_iH^i(\bK^\tt_{D'})[-i]$, we have 
$$\op_i\res_D^{D'}(H^i(\bK^\ss_D))[-i]\cong\op_{\tt\in\ct,i}H^i(\bK^\tt_{D'})[-i-d_\tt].\tag a$$
By 31.14, $\res_D^{D'}(H^i(\bK^\ss_D))$ is a perverse sheaf on $D'$. Hence taking $H^i$ for both sides of (a) we
obtain
$$\res_D^{D'}(H^i(\bK^\ss_D))\cong\op_{\tt\in\ct}H^{i-d_\tt}(\bK^\tt_{D'}).\tag b$$
In particular, if $A\in\hD^{un}$ then, $\res_D^{D'}(A)$ is a direct sum of objects in $\hD'{}^{un}$. Hence 
$A\m\res^{D'}_D(A)$ (with $A\in\hD^{un}$) defines a group homomorphism $\ck^{un}(D)@>>>\ck^{un}(D')$ and a 
$\QQ$-linear map $\ck^{un}(D)_\QQ@>>>\ck^{un}_\QQ(D')$ denoted again by $\res^{D'}_D$. Taking alternating sum over
$i$ in (b) we obtain
$$\res^{D'}_D(gr_1(\bK^\ss_D))=\sum_{\tt\in\ct}(-1)^{d_\tt}gr_1(\bK^\tt_{D'}).\tag c$$
For any $\x\in\ck^{un}_\QQ(D)$, $\x'\in\ck^{un}_\QQ(D')$ we have
$$(\res_D^{D'}(\x):\x')=(\x:\ind_{D'}^D(\x'))\tag d$$
where the first $(:)$ refers to $D'$ and the second $(:)$ refers to $D$. Indeed, we can assume that 
$\x=A\in\hD^{un}$, $\x'=A'\in\hD'{}^{un}$; in this case (d) follows from the equalities in 30.9 and the 
semisimplicity of the perverse sheaves $\res_D^{D'}(A)$, $\ind_{D'}^D(A')$. 

The following subspaces of $\ck^{un}_\QQ(D)$ coincide:

-the subspace $(1)$ spanned by the $R_E$ (with $E\in\Mod(\tWW)$);

-the subspace $(2)$ spanned by the $R_E$ (with $E\in\Irr(\tWW)$);

-the subspace $(3)$ spanned by the elements $gr_1(K^x_D)$ (with $x\in\WW$);

-the subspace $(4)$ spanned by the elements $gr_1(K^\ss_D)$ for various sequences $\ss$ in $\II$.
\nl
Indeed $(1)\sub(3)$ by 44.6(b); $(3)\sub(2)$ by 44.7(i); $(2)\sub(1)$ obviously; moreover, $(3)=(4)$ by the 
arguments in 31.7. We denote any of the four subspaces above by $V_D$. We define similarly a subspace $V_{D'}$ of
$\ck^{un}_\QQ(D')$. We show:
$$\res_D^{D'}(R_E)=R_{E|_{\tWW_I}}\tag e$$
where $E|_{\tWW_I}\in\Mod(\tWW_I)$ is the restriction of $E$. From (c) we see that $\res^{D'}_D$ maps $V_D$ into 
$V_{D'}$. Thus both sides of (e) are in $V_{D'}$. Now the restriction of $(:)$ (for $D'$) to $V_{D'}$ is 
nondegenerate (we use the analogue of 44.7(l) for $D'$). Hence to prove (e) it is enough to show that
$$(\res_D^{D'}(R_E):R_{E'})=(R_{E|_{\tWW_I}}:R_{E'})\tag f$$
for any $E'\in\Mod(\tWW_I)$. By (d) and 44.11(a), the left hand side of (f) is equal to 
$$(R_E:\ind_{D'}^D(R_{E'})=(R_E:R_{\ind_{\tWW_I}^{\tWW}E'}).$$
Using 44.7(l) for $D$ and for $D'$ we see that it is enough to use the equality
$$|\WW|\i\sum_{x\in\WW}\tr(x\vp,E)\tr(x\vp,\ind_{\tWW_I}^{\tWW}E')
=|\WW_I|\i\sum_{x\in\WW_I}\tr(x\vp,E)\tr(x\vp,E').$$
which follows from the standard character formula for an induced representation. This proves (f) and hence (e).

\subhead 44.13\endsubhead
Let $x\in\WW$ be such that for any $y\in\WW$ we have $yx\vp y\i\vp\i\n\WW_I$. We show:
$$\res_D^{D'}(gr_1(K^x_D))=0.\tag a$$
Using 44.7(i), we see that it is enough to show:
$$(-1)^{\dim G}\sum_{E\in\fE}\tr(x\vp,E)\res_D^{D'}(R_E)=0.$$
Using 44.12(e) and 44.6(b) for $D'$, we see that left hand side is
$$\sum_{E\in\fE}\tr(x\vp,E)R_{E|_{\tWW_I}}=
|\WW_I|\i\sum_{E\in\fE}\sum_{z\in\WW_I}(-1)^{\dim G'}\tr(z\vp,E)\tr(x\vp,E)gr_1(K^z_{D'}).$$
To show that this is zero it is enough to show that for any $z\in\WW_I$ we have
$$\sum_{E\in\fE}\tr(z\vp,E)\tr(x\vp,E)=0.$$
The left hand side is equal to $|\tWW|\i|\WW|$ times $\sum_E\tr(z\vp,E)\tr((x\vp)\i,E)$ where $E$ runs over the 
simple $\bbq[\tWW]$-modules up to isomorphism. (A module $E$ whose restriction to $\WW$ is not simple contributes
$0$ to the last sum.) It is enough to show that the last sum is $0$. It is also enough to show that
$z\vp$ and $x\vp$ are not conjugate in $\tWW$. But this follows from our assumption on $x$. This proves (a).

\subhead 44.14\endsubhead
An element $w\in\WW$ is said to be {\it $D$-anisotropic} if the following condition holds:
for any $x\in\WW$, $I\subsetneqq\II$ such that $\e(I)=I$ we have $xw\e(x)\i\notin\WW_I$. Let $A\in\hD^{un}$. 

We show:

(a) {\it $A$ is cuspidal if and only if any $w\in\WW$ such that $(A:gr_1(K^w_D))\ne0$ is $D$-anisotropic.}
\nl
Assume first that $A$ is not cuspidal. By 31.15 there exists $I\subsetneqq\II$, $\e(I)=I$ and $P\in\cp_I$ (so that
$N_DP\ne\em$) such that setting $D'=N_DP/U_P$, $G'=N_GP/U_P$ we have $\res_D^{D'}(A)\ne0$. By 31.14 and 44.12, 
$\res_D^{D'}(A)$ is an $\NN$-linear combinations of objects in $\hD'{}^{un}$. Hence there exists $x\in\WW_I$ and 
$i\in\ZZ$ such that $(\res_D^{D'}(A):H^i(\bK^x_D))\ne0$. Using 44.7(m) for $D'$  we see that there exists 
$E'\in\Irr(\tWW_I)$ such that $(\res_D^{D'}(A):R_{E'})\ne0$. Hence there exists $y\in\WW_I$ such that 
$(\res_D^{D'}(A):gr_1(K^y_{D'}))\ne0$. Using 44.12(d) we deduce $(A:\ind_{D'}^D(gr_1(K^y_{D'}))\ne0$ and using 
44.10(d) we see that $(A:gr_1(K^y_D))\ne0$. Since $y\in\WW_I$, $y$ is not $D$-anisotropic.

Conversely, assume that there exist $w\in\WW,x\in\WW$, $I\subsetneqq\II$ such that $(A:gr_1(K^w_D))\ne0$, 
$\e(I)=I$ and $xw\e(x)\i\in\WW_I$. Using 44.7(n) we see that we can assume that $x=1,w\in\WW_I$. Choose 
$P\in\cp_I$ (so that $N_DP\ne\em$) and set $D'=N_DP/U_P$, $G'=N_GP/U_P$. Using 44.10(d) we see that 
$(A:\ind_{D'}^D(gr_1(K^w_D)))\ne0$. Using 44.12(d) we see that $(\res_D^{D'}(A):gr_1(K^w_{D'}))\ne0$ so that 
$\res_D^{D'}(A)\ne0$. Thus $A$ is not cuspidal. This proves (a).

We show: 

(b) {\it Let $w\in\WW$ be such that $w$ is $D$-anisotropic. Then $l(w)=|\II_\e|\mod2$ where $\II_\e$ is the set of
orbits of $\e:\II@>>>\II$.}
\nl
We use the notation in 42.7. We consider the equality 
$$(-1)^{|\II|}H_c^{|\II|}(\cv_\RR)=\sum_\et(-1)^{r_\et}H^{r_\et}(\cv^\et_\RR)$$
(see 42.7) in the Grothendieck group of $\WW^D$-modules. Taking the trace of $w\uD\in\WW^D$ we obtain
$$(-1)^{|\II|}\det(w\uD,\cv_\RR)=\sum_\et t_\et$$
where 
$$t_\et=(-1)^{r_\et}\tr(w\uD,\op_{J\in\et}\op_{F\in\cf_J}\L^{r_\et}(|F|).$$
Since $w\uD$ permutes the summands in the last direct sum, we have $t_\et=0$ unless there exist $J\in\et$ and 
$F\in\cf_J$ such that $\uD(J)=J$ and $w\uD(F)=F$. For such $J,F$ we can find $F_J\in\cf_J$ such that
$\uD(F_J)=F_J$ and $\{y\in\WW;y(F_J)=F_J\}=\WW_J$; moreover, $F=x\i(F_J)$ for some $x\in\WW$ and
$w\e(x)\i(F_J)=x\i(F_J)$ so that $xw\e(x)\i(F_J)=F_J$ and $xw\e(x)\i\in\WW_J$. Since $w$ is $D$-anisotropic we see
that $J=\II$. Thus $t_\et=0$ unless $\et=\{\II\}$. On the other hand, if $\et=\{\II\}$ then $\cf_J=\{0\}$, 
$r_\et=0$ and $t_\et=1$. Thus we have $(-1)^{|\II|}\det(w\uD,\cv_\RR)=1$. Note that $\det(w,\cv_\RR)=(-1)^{l(w)}$.
Since $\uD$ permutes a basis of $\cv_\RR$ indexed by $\II$ (according to $\e$) we have 
$\det(\uD,\cv_\RR)=(-1)^{|\II|-|\II_\e|}$. We see that $(-1)^{l(w)}(-1)^{|\II_\e|}=1$. This proves (b).

\subhead 44.15\endsubhead
Let $P$ be a parabolic subgroup of $G^0$ such that $N_DP\ne\em$. Let $D'=N_DP/U_P$ (a connected component of 
$N_GP/U_P$). We show:

(a) {\it If $A'\in\hD'{}^{un}$, $A\in\hD^{un}$, are such that $A$ appears with non-zero coefficient in 
$\ind_{D'}^D(A')$ (or equivalently $A'$ appears with non-zero coefficient in $\res_D^{D'}(A)$) then 
$\ee^A=\ee^{A'}$. Moreover, $\codim(\supp(A))=\codim(\supp(A'))$.}
\nl
We can find $I\sub\II$, $\e(I)=I$ such that $P\in\cp_I$ and $w\in\WW_I$, $i\in\ZZ$ such that $A'$ is a direct 
summand of $H^i(\bK^w_{D'})$. Then $\ind_{D'}^D(A')$ is a direct summand of $\ind_{D'}^D(H^i\bK^w_{D'})$ hence a 
direct summand of $H^{i+2\dim U_P}\bK^w_D$ (see 44.10(c)). It follows that $A$ is a direct summand of 
$H^{i+2\dim U_P}\bK^w_D$. By definition we have
$\ee^{A'}=(-1)^{i+\dim(P/U_P)}$, $\ee^A=(-1)^{i+2\dim U_P+\dim G^0}$. Thus, $\ee^A=\ee^{A'}$. This proves the 
first statement of (a). We can find a parabolic subgroup $P_1$ of $G^0$ such that $N_DP_1\ne\em$, $P_1\sub P$ and
$A_1\in\hD_1^{unc}$ (where $D_1=NDP_1/U_{P_1}$) such that $A'$ is a component of $\ind_{D_1}^{D'}(A_1)$ hence $A$
is a component of $\ind_{D_1}^{D}(A_1)$. To prove the second statement of (a) it is enough to show that 
$(-1)^{\codim(\supp(A))}=(-1)^{\codim(\supp(A_1))}$, $(-1)^{\codim(\supp(A'))}=(-1)^{\codim(\supp(A_1))}$. Thus we
are reduced to the case where $A'$ is cuspidal. In this case, by 3.13(b) we have 
$\dim\supp(A)=\dim(G^0)-\dim(P/U_P)+\dim\supp(A')$. Thus, $\codim(\supp(A))=\codim(\supp(A'))$ and (a) is proved.

We show:

(b) {\it If $A\in\hD^{un}$ and $A^\circ\in\hD^{un}$ is defined by $\dd(A)=(-1)^{\codim(\supp(A))}A^\circ$ (see 
44.8(a)) then $\ee^{A^\circ}=\ee^{A}$.}
\nl
If $P,D'$ are as in (a) then, by (a), $\ind_{D'}^D\res_D^{D'}(A)$ is a linear combination of objects
$A_1\in\hD^{un}$ with $\ee^{A_1}=\ee^A$. Since $\dd(A)$ is an alternating sum of elements of the form 
$\ind_{D'}^D\res_D^{D'}(A)$, we see that $\dd(A)$ is a linear combination of objects $A_1\in\hD^{un}$ with 
$\ee^{A_1}=\ee^A$. Now (b) follows.

Let $x\in\WW$. We show:

(c) {\it The element $R_{\ale_{x\vp}}\in\ck^{un}_\QQ(D)$ is a $\ZZ$-linear combination of objects $A\in\hD^{un}$
such that $\ee^A=(-1)^{l(x)-\aa(x)}$.}
\nl
Let $\boc$ be the two-sided cell containing $x$. Using 43.12(b), for any $A\in\hD^{un}$ we have (with notation in
43.12):    
$$\align&(A:R_{\ale_{x\vp}})=\sum_{E\in\un\Irr(\tWW)}\fra{1}{2}(\tr(c_{x\vp}^\da,E^v;-\aa(x))\\&
-\sum_{y,j;y\prec x,j>0}(-1)^{-l(x)+l(y)}a_{y,x;j}\tr(c_{y\vp}^\da,E^v;-\aa(x)-j))(A:R_E).\tag d\endalign$$
From 44.7(j) we have for any $A\in\hD^{un}$ and $z\in\WW$:
$$\align&\sum_{j\in\ZZ}(\dd(A):H^j(\bK^z_D))(-v)^j\\&
=(-1)^{\dim G}v^{\dim G+l(z)}\sum_{E\in\un\Irr(\tWW)}\fra{1}{2}(\dd(A):R_E)\tr(c_{z\vp},E^v)\\&=
(-1)^{\dim G}v^{\dim G+l(z)}\sum_{E\in\un\Irr(\tWW)}\fra{1}{2}(A:R_{E\ot\sgn})\tr(c_{z\vp},E^v)\\&=
(-1)^{\dim G}v^{\dim G+l(z)}\sum_{E\in\un\Irr(\tWW)}\fra{1}{2}(A:R_E)\tr(c_{z\vp},(E\ot\sgn)^v)\\&=
(-1)^{\dim G}v^{\dim G+l(z)}\sum_{E\in\un\Irr(\tWW)}\fra{1}{2}(A:R_E)\tr(c_{z\vp}^\da,E^v).\endalign$$
(We have used 44.8(c), 43.4(c).) Hence for any $N\in\ZZ$ we have
$$\sum_{E\in\un\Irr(\tWW)}\fra{1}{2}(A:R_E)\tr(c_{z\vp}^\da,E^v;N)
=(\dd(A):H^{N+\dim G+l(z)}(\bK^z_D))(-1)^{N+l(z)}.$$
Introducing this in (c) we obtain
$$\align&(A:R_{\ale_{x\vp}})
=(-1)^{l(x)-\aa(x)}(\dd(A):H^{\dim G+l(x)-\aa(x)}(\bK^x_D))\\&-
\sum_{y,j;y\prec x,j>0}a_{y,x;j}(-1)^{l(x)-\aa(x)-j}(\dd(A):H^{\dim G+l(y)-\aa(x)-j}(\bK^y_D)).\tag e\endalign$$
Since $a_{y,x;j}$ are integers (see 43.12) we see that $(A:R_{\ale_{x\vp}})\in\ZZ$. Assume now that 
$(A:R_{\ale_{x\vp}})\ne0$. Using (e) and 43.12 we see that either 

$(A^\circ:H^{\dim G+l(x)-\aa(x)}(\bK^x_D))\ne0$
\nl
or there exist $y,j$ such that $j=l(x)+l(y)\mod2$,

$(A^\circ:H^{\dim G+l(y)-\aa(x)-j}(\bK^y_D))\ne0$
\nl
(here $A^\circ$ is as in 44.8(a)). In the first case we have 
$\ee^{A^\circ}=(-1)^{l(x)-\aa(x)}$. In the second case we have 
$\ee^{A^\circ}=(-1)^{l(y)-\aa(x)-j}=(-1)^{l(x)-\aa(x)}$ since $j=l(x)+l(y)\mod2$. This implies (c) in view of (b).

Note that $D$ has property $\tfA$ (see 44.7) if and only if for any $A\in\hD^{un}$ we have 
$\ee^A=(-1)^{\codim(\supp(A))}$). 

\subhead 44.16\endsubhead
We show that if $D$ has property $\tfA$ then for any $A\in\hD^{un}$, $w\in\WW$, $i\in\ZZ$ we have
$$(-1)^{i+\dim G}(A:\dd(H^i(\bK^w_D)))\in\NN.\tag a$$
Indeed the expression (a) is equal to $(-1)^{i+\dim G}(\dd(A):H^i(\bK^w_D))$ (see 38.10(e)). If this is $\ne0$
then it is equal to $(-1)^{\codim(\supp(A))}\ee^{A^\circ}(A^\circ:H^i(\bK^w_D))$. By property $\tfA$ for $A^\circ$
and 44.8(e), this is equal to
$$(-1)^{\codim(\supp(A))}(-1)^{\codim(\supp(A^\circ))}(A^\circ:H^i(\bK^w_D))=(A^\circ:H^i(\bK^w_D))\in\NN.$$
This proves (a).

\subhead 44.17\endsubhead
Let $x\in\WW$ and let $\boc$ be the two-sided cell of $\WW$ that contains $x$. Let $a$ be the value of 
$\aa:\WW@>>>\NN$ on $\boc$. We show that in $\ck^{un}_\QQ(D)$ we have:
$$\align&(-1)^{-a+l(x)}H^{-a+l(x)+\dim G}(\bK^x_D)\\&=
R_{\ale_{x\vp}\ot\sgn}+\QQ\text{-linear combination of elements $R_{\ale_{x'\vp}\ot\sgn}$ with }x'\prec x,\tag a
\endalign$$
$$\align&(-1)^{-a+l(x)}\dd(H^{-a+l(x)+\dim G}(\bK^x_D))\\&=
R_{\ale_{x\vp}}+\QQ\text{-linear combination of elements $R_{\ale_{x'\vp}}$ with }x'\prec x.\tag b\endalign$$
By 44.7(m), the left hand side of (b) is equal to $\sum_E\fra{1}{2}\tr(c_x\tT_\vp,E^v;-a)\dd(R_E)$. By 44.8(c) and
43.4(b), 43.6(b), this equals 
$$\align&\sum_E\fra{1}{2}\tr(c_x\tT_\vp,E^v;-a)R_{E\ot\sgn}=\sum_E\fra{1}{2}\tr(c_x\tT_\vp,(E\ot\sgn)^v;-a)R_E\\&=
\sum_E\fra{1}{2}\tr(c_{x\vp}^\da,E^v;-a)R_E=\sum_{E;\boc_E\prq\boc}\fra{1}{2}\tr(c_{x\vp}^\da,E^v;-a)R_E=b'+b''
\endalign$$
where 
$$\align&b'=\sum_{E;\boc_E=\boc}\fra{1}{2}\tr(c_{x\vp}^\da,E^v;-a)R_E
=\sum_{E;\boc_E=\boc}\fra{1}{2}\tr(t_x\vp,E^\iy)R_E\\&=\sum_E\fra{1}{2}\tr(t_x\vp,E^\iy)R_E=R_{\ale_{x\vp}},
\endalign$$
$$b''=\sum_{E;\boc_E\prec\boc}\fra{1}{2}\tr(c_{x\vp}^\da,E^v;-a)R_E.$$
Now $b''$ is a $\ZZ$-linear combination of elements of the form $R_E$ where $E$ is such that $\boc_E\prec\boc$ and
these elements are $\QQ$-linear combinations of elements of the form $R_{\ale_{x'\vp}}$ for various $x'\in\WW$ 
such that $x'\prec x$, by 43.10(b). This proves (b). Now (a) is obtained by applying $\dd$ to both sides of (b) 
and using the equality $\dd(R_\ph)=R_{\ph\ot\sgn}$ for any $\ph\in\car(\tWW)$ (see 44.8(c)).

Now let $a'$ be the value of $\aa:\WW@>>>\NN$ on the two-sided cell $w_0\boc=\boc w_0$. We show:
$$\align&(-1)^{-a'+l(w_0x)}H^{-a'+l(w_0x)+\dim G}(\bK^{w_0x}_D)\\&=R_{\ale_{w_0x\vp}\ot\sgn}+\QQ-\text{linear 
combination of elements $R_{\ale_{w_0x'\vp}\ot\sgn}$ with }x\prec x'.\tag c\endalign$$
This is obtained by replacing $x$ by $w_0x$ in (a) and noting that for $y\in\WW$ we have $w_0y\prec w_0x$ if and 
only if $x\prec y$.

{\it In the remainder of this section we assume that $D$ satisfies property $\tfA$ (in addition to property 
$\fA$).} 
\nl
For any $x\in\WW$ we set $r_x=R_{\ale_{x\vp}}$, $\tir_x=(-1)^{-\aa(w_0x)+l(w_0x)}R_{\ale_{w_0x\vp}\ot\sgn}$. We 
note the following properties of the elements $r_x,\tir_x$.

(i) $(r_x:r_{x'})=0$ whenever $x\not\si x'$;

(ii) for any two-sided cell $\boc$, the $\QQ$-vector space spanned by $\{r_x;x\in\boc\}$ coincides with the 
$\QQ$-vector space spanned by $\{\tir_x;x\in\boc\}$;

(iii) for any $x\in\WW$ there exist $d_{x,x'}\in\QQ$ defined for $x'\prec x$ such that \lb
$(A:r_x+\sum_{x';x'\prec x}d_{x,x'}r_{x'})\in\NN$ for any $A\in\hD^{un}$; 

(iv) for any $x\in\WW$ there exist $\td_{x,x'}\in\QQ$ defined for $x\prec x'$ such that \lb
$(A:\tir_x+\sum_{x';x\prec x'}\td_{x,x'}\tir_{x'})\in\NN$ for any $A\in\hD^{un}$.
\nl
In the setup of (ii), let $V_\boc$ be the $\QQ$-vector space spanned by $R_E$ with $E\in\Irr(\tWW)$ such that
$\boc_E=\boc$. From the definitions, for any $x\in\boc$, $r_x$ belongs to $V_\boc$. Conversely, for any 
$E\in\Irr(\tWW)$ such that $\boc_E=\boc$, $R_E$ belongs to the first vector space in (ii), by 43.10(b). Thus the 
first vector space in (ii) is equal to $V_\boc$. Let $V'_\boc$ be the $\QQ$-vector space spanned by 
$R_{E'\ot\sgn}$ with $E'\in\Irr(\tWW)$ such that $\boc_{E'}=w_0\boc$. From the definitions, for any $x\in\boc$, 
$\tir_x$ belongs to $V'_\boc$. Conversely for any $E'\in\Irr(\tWW)$ such that $\boc_{E'}=w_0\boc$, $R_{E'\ot\sgn}$
belongs to the second vector space in (ii), by 43.10(b). Thus the second vector space in (ii) is equal to 
$V'_\boc$. If $E'\in\Irr(\tWW)$ then we have $\boc_{E'}=w_0\boc$ if and only $\boc_{E'\ot\sgn}=\boc$ (a known 
property of two-sided cells). It follows that $V_\boc=V'_\boc$ and (ii) is proved.

We prove (i). Let $\boc$, $\boc'$ be the two-sided cells that contain $x,x'$ respectively. Assume that 
$\boc\ne\boc'$. It is enough to show that $(h:h')=0$ for any $h\in V_\boc,h'\in V_{\boc'}$. Hence it is enough to
show that if $E,E'\in\Irr(\tWW)$ are such that $\boc_E=\boc,\boc_{E'}=\boc'$ then $(R_E:R_{E'})=0$. This follows 
from 44.7(l) since $E,E'$ have nonisomorphic restrictions to $\QQ[\WW]$.

Now (iv) follows from (c) and (iii) follows from (b) in view of 44.16(a).

From (i)-(iv) we deduce, by a general result in \cite{\CS, III, 16.8}, that:
$$(A:r_x)\in\NN,\qua (A:\tir_x)\in\NN\text{ for any }A\in\hD^{un},x\in\WW.\tag d$$
We show:

(e) {\it Let $A\in\hD^{un}$ and let $E,E'\in\Irr(\tWW)$ be such that $(A:R_E)\ne0$, $(A:R_{E'})\ne0$. Then
$\boc_E=\boc_{E'}$.}
\nl
By the proof of (ii) we see that there exists $x\in\boc_E$ such that $(A:r_x)\ne0$; similarly there exists 
$x'\in\boc_{E'}$ such that $(A:r_{x'})\ne0$. Using this and (d) we deduce $(A:r_x)>0$, $(A:r_{x'})>0$. It follows
that $(r_x:r_{x'})>0$. (By (d), $(r_x:r_{x'})$ is a sum of terms in $\NN$, at least one of which is $>0$.) Again 
by the proof of (ii) we have 
$$r_x=\sum_{E_1;\boc_{E_1}=\boc_E}s_{E_1}R_{E_1}, r_{x'}=\sum_{E_2;\boc_{E_2}=\boc_{E'}}s'_{E_2}R_{E_2},$$
where $s_{E_1}\in\QQ$, $s'_{E_2}\in\QQ$. From $(r_x:r_{x'})\ne0$ it follows that there exist $E_1,E_2$ such that 
$\boc_{E_1}=\boc_E$, $\boc_{E_2}=\boc_{E'}$, $(R_{E_1}:R_{E_2})\ne0$. From 44.7(l) we deduce that $E_1,E_2$ have 
isomorphic restrictions to $\QQ[\WW]$ hence $\boc_{E_1}=\boc_{E_2}$. It follows that $\boc_E=\boc_{E'}$. This 
proves (e). 

\proclaim{Proposition 44.18} Recall that $D$ is assumed to have property $\fA$ and property $\tfA$. Let 
$A\in\hD^{un}$.

(a) There exists a well defined two-sided cell $\boc'_A$ in $\WW$  such that whenever $E\in\Irr(\tWW)$ and 
$(A:R_E)\ne0$, we have $\boc_E=\boc'_A$. Moreover we have $\e(\boc'_A)=\boc'_A$.

(b) We have $w_0\boc'_A=\boc_A$ where $\boc_A$ is as in 41.4.
\endproclaim
(a) follows from 44.17(e) and 43.6(f). We prove (b). Recall (41.8) that

(c) {\it $A\dsv\bK^x_D$ for some $x\in\boc_A$; if $x'\in\WW$ and $A\dsv\bK^{x'}_D$ then $x\prq x'$.}
\nl
We show:

(d) {\it if $E\in\Irr(\tWW)$ is such that $(A:R_E)\ne0$ then $\boc_A\prq w_0\boc_E$.}
\nl
Using 44.6(c) we see that
$$|\WW|\i\sum_{i\in\ZZ}\sum_{x\in\WW}(-1)^{i+\dim G}\tr(\tc_{x\vp}|_{v=1},E)(A:H^i(\bK^x_D))\ne0.$$
Hence there exist $x\in\WW,i\in\ZZ$ such that $\tr(\tc_{x\vp}|_{v=1},E)\ne0$ and $(A:H^i(\bK^x_D))\ne0$. Using (c)
we deduce that $y\prq x$ for some $y\in\boc_A$. From the definitions we have
$$\tc_{x\vp}=(-1)^{l(w_0x)}\tT_{w_0}c_{w_0x\vp}^\da.$$
It follows that $\tr(w_0c_{w_0x\vp}^\da|_{v=1},E)\ne0$. Thus the action of $c_{w_0x\vp}^\da|_{v=1}$ on $E$ is 
$\ne0$. Using 43.6(b) we see that $z\prq w_0x$ for some $z\in\boc_E$. Hence $x\prq w_0z$. Since $y\prq x$, we have
$y\prq w_0z$. Since $y\in\boc_A$ we have $\boc_A\prq w_0\boc_E$. This proves (d).

We show:

(e) {\it There exists $E\in\Irr(\tWW)$ such that $(A:R_E)\ne0$ and $w_0\boc_E=\boc_A$.}
\nl
Let $x$ be as in (c). We have $\sum_{j\in\ZZ}(A:H^j(\bK^x_D))(-v)^j\ne0$. Using 6.7(c) we deduce that
$$v^{\dim G+l(x)}\sum_{E\in\fE}b_{A,E}\tr(c_{x\vp},E^v)\ne0.$$
Hence there exists $E\in\Irr(\tWW)$ such that $(A:R_E)\ne0$ and $\tr(c_{x\vp},E^v)\ne0$ that is,
$\tr(c_{x\vp}^\da,(E^\da)^v)\ne0$. The last condition implies, in view of 43.6(b) that $z\prq x$ for some 
$z\in\boc_{E^\da}=w_0\boc_E$. Thus, $w_0\boc_E\prq\boc_A$. Since
$\boc_A\prq w_0\boc_E$ by (d), it follows that $\boc_A=w_0\boc_E$. This proves (e).

From (e) we see that $w_0\boc'_A=\boc_A$. The proposition is proved.

\subhead 44.19\endsubhead
For any $\e$-stable two-sided cell $\boc$ of $\WW$ let $\hD^{un}_\boc$ be the category whose objects are those
$A\in\hD^{un}$ such that $\boc'_A=\boc$ (see 44.18) and let $\ck^{\boc}(D)$ be the subgroup of $\ck^{un}(D)$ 
generated by the various $A\in\hD^{un}_\boc$ up to isomorphism. We have $\ck^{un}(D)=\op_\boc\ck^{\boc}(D)$ where
$\boc$ runs over the $\e$-stable two-sided cells of $\WW$. We show:

(a) {\it $A\m A^\circ$ (see 44.8(a)) induces a bijection between the set of isomorphism classes in $\hD^{un}_\boc$
and the set of isomorphism classes in $\hD^{un}_{w_0\boc}$; it also induces an isomorphism 
$\ck^{\boc}(D)@>\si>>\ck^{w_0\boc}(D)$.}
\nl
Let $A\in\hD^{un}_\boc$. Then $(A:R_E)\ne0$ for some $E\in\Irr(\tWW)$ such that $\boc_E=\boc$. We have
$(\dd(A):\dd(R_E))\ne0$ and $(A^\circ:R_{E\ot\sgn})\ne0$ (see 44.8(d)). Thus 
$A^\circ\in\hD^{un}_{\boc_{E\ot\sgn}}=\hD^{un}_{w_0\boc}$. The remaining statements of (a) are immediate.

\subhead 44.20\endsubhead
Let $I$ be a subset of $\II$ such that $\e(I)=I$. We fix a two-sided cell $\boc'$ of $\WW_I$ (see 26.1) such that
$\e(\boc')=\boc'$. There is a unique two-sided cell $\boc$ of $\WW$ such that $\boc'\sub\boc$; we must have 
$\e(\boc)=\boc$. 

Let $\Irr_\boc(\tWW)=\{E\in\Irr(\tWW);\boc_E=\boc\}$, 
$\Irr_{\boc'}(\tWW_I)=\{E'\in\Irr(\tWW_I);\boc_{E'}=\boc'\}$.

Let $\car_\boc(\tWW)$ be the subgroup of $\car(\tWW)$ generated by the elements $\ph_E$ with 
$E\in\Irr_\boc(\tWW)$. Let $\car_{\boc'}(\tWW_I)$ be the subgroup of $\car(\tWW_I)$ generated by the elements 
$\ph_{E'}$ with $E'\in\Irr_{\boc'}(\tWW_I)$. From 43.11(b) we see that 
$$J_{\tWW_I}^{\tWW}:\car(\tWW_I)@>>>\car(\tWW)\text{ satisfies }
J_{\tWW_I}^{\tWW}(\car_{\boc'}(\tWW_I))\sub\car_\boc(\tWW).\tag a$$
Let $\ck^{\boc}(D)$ be as in 44.19. We define similarly $\ck^{\boc'}(D')$. 
Define a $\QQ$-linear map $p_\boc:\QQ\ot\ck^{un}(D)@>>>\QQ\ot\ck^\boc(D)$ by $A\m A$ if 
$A\in\hD^{un}_\boc$ and $A\m 0$ if $A\in\hD^{un},\boc'_A\ne\boc$; this restricts to a homomorphism
$\ck^{un}(D)@>>>\ck^{\boc}(D)$. Note that for $E_1\in\Irr(\tWW)$ we have $R_{E_1}\in\QQ\ot\ck^{\boc_{E_1}}(D)$ 
hence 

(b) {\it $p_\boc(R_{E_1})=R_{E_1}$ if $\boc_{E_1}=\boc$ and $p_\boc(R_{E_1})=0$ if $\boc_{E_1}\ne\boc$.}
\nl
Let $E'\in\Irr_{\boc'}(\tWW_I)$. We show:
$$R_{J_{\tWW_I}^{\tWW}(\ph_{E'})}=p_\boc(R_{\ind_{\tWW_I}^{\tWW}E'}).\tag c$$
By 44.7(o) and (b), both sides of (c) are integer combinations of elements of form $R_{E_1}$ with 
$E_1\in\fE$. Hence (using 44.7(l)) it is enough to show that for any $E_1\in\fE$ we have
$$(R_{J_{\tWW_I}^{\tWW}(\ph_{E'})}:R_{E_1})=(p_\boc(R_{\ind_{\tWW_I}^{\tWW}E'}):R_{E_1}).\tag d$$
If $\boc_{E_1}\ne\boc$ then from (b) we see that the right hand side of (d) is zero; moreover, since
$\ph_{E'}\in\car_{\boc'}(\tWW_I)$ we have $J_{\tWW_I}^{\tWW}(\ph_{E'})\sub\car_\boc(\tWW)$ (see (a)) hence 
$R_{J_{\tWW_I}^{\tWW}(\ph_{E'})}\in\ck^{\boc}(D)$ so that the left hand side of (d) is also zero. Thus, we 
may assume that $\boc_{E_1}=\boc$. In this case (d) can be rewritten in the form
$$(R_{J_{\tWW_I}^{\tWW}(\ph_{E'})}:R_{E_1})=(R_{\ind_{\tWW_I}^{\tWW}E'}:R_{E_1})$$
or equivalently (using 44.7(l)) in the form
$$\align&\sum_{E\in\un\Irr(\tWW);a_{E'}=a_E}\la E',E\ra|\WW|\i\sum_{u\in\WW}\tr(u\vp,E)\tr(u\vp,E_1)\\&
=|\WW|\i\sum_{x\in\WW}\tr(x\vp,\ind_{\tWW_I}^{\tWW}E')\tr(x\vp,E_1).\tag e\endalign$$ 
The right hand side of (e) can be rewritten as $|\WW_I|\i\sum_{z\in\WW_I}\tr(z\vp,E')\tr(z\vp,E_1)$; substituting 
$\tr(z\vp,E_1)=\sum_{E'_1\in\un\Irr(\tWW_I)}\la E'_1,E_1\ra\tr(z\vp,E')$ (see 43.9(a)) this becomes
$$\align&|\WW_I|\i\sum_{z\in\WW_I}\tr(z\vp,E')\sum_{E'_1\in\un\Irr(\tWW_I)}\la E'_1,E_1\ra\tr(z\vp,E'_1)\\&=
\sum_{E'_1\in\un\Irr(\tWW_I)}\la E'_1,E_1\ra\a(E',E'_1)=\la E',E_1\ra-\la E'\ot\io,E_1\ra\endalign$$
where $\a(E',E'_1)$ is $1$ if $E'\cong E'_1$, is $-1$ if $E'\cong E'_1\ot\io$ and is $0$ otherwise. Now in the 
left hand side of (e) the second sum is zero unless $E$ is isomorphic to $E_1$ or to $E_1\ot\io$ in which case we
have automatically $a_{E'}=a_E$ (since $a_E=a_{E_1}$). Thus the left hand side of (e) is equal to
$$\align&\sum_{E\in\un\Irr(\tWW)}\la E',E\ra|\WW|\i\sum_{u\in\WW}\tr(u\vp,E)\tr(u\vp,E_1)\\&=
\sum_{E\in\un\Irr(\tWW)}\la E',E\ra \a(E,E_1)=\la E',E_1\ra-\la E',E_1\ot E\ra.\endalign$$
This proves (e) and hence (c).

For any $A'\in\hD'{}^{un}_{\boc'}$ we set $\tind_{D'}^D(A')=p_\boc(\ind_{D'}^D(A))$, (see 44.13). Now 
$A'\m\tind_{D'}^D(A')$ defines a group homomorphism $\ck^{\boc'}(D')@>>>\ck^{\boc}(D)$ and a 
$\QQ$-linear map $\QQ\ot\ck^{\boc'}(D')@>>>\QQ\ot\ck_{\boc}(D)$; these are denoted again by $\tind_{D'}^D$.

Let $\ph'\in\car_{\boc'}(\tWW_I)$. We show:
$$\tind_{D'}^D(R_{\ph'})=R_{J_{\tWW_I}^{\tWW}(\ph')}.\tag f$$
We may assume that $\ph'=\ph_{E'}$ where $E'\in\Irr_{\boc'}(\tWW_I)$. From the definitions we have 
$R_{\ph_{E'}}\in\QQ\ot\ck^{\boc'}(D')$ and $\tind_{D'}^D(R_{\ph_{E'}})\in\QQ\ot\ck^{\boc}(D)$. Applying $p_\boc$
to the identity
$$\ind_{D'}^D(R_{\ph_{E'}})=R_{\ind_{\tWW_I}^{\tWW}E'}\in\ck^{un}_\QQ(D)$$
(see 44.14(a)) we obtain
$$\tind_{D'}^D(R_{\ph_{E'}})=p_\boc(R_{\ind_{\tWW_I}^{\tWW}E'}).$$
Now (f) follows from (c).

For any $x\in\boc$ we have $\ale_{x\vp}\in\car_\boc(\tWW)$. Similarly for any $x\in\boc'$ we have 
$\ale^I_{x\vp}\in\car_{\boc'}(\tWW_I)$. Combining (f) (with $\ph'=\ale^I_{x\vp}$, $x\in\boc'$) with 43.10(c) we 
see that
$$\tind_{D'}^D(R_{\ale^I_{x\vp}})=R_{\ale^I_{x\vp}}.\tag g$$
We define a homomorphism ${}'J^{\tWW_I}_{\tWW}:\car(\tWW)@>>>\car(\tWW_I)$ by 
$${}'J^{\tWW_I}_{\tWW}(\ph_E)=\sum_{E'\in\un\Irr(\tWW_I);a_{E'}=a_E}\la E',E\ra\ph_{E'}$$
for any $E\in\Irr(\tWW)$.

Let $\ph\in\car_\boc(\tWW)$ and let $A'\in\hD'{}^{un}_{\boc'}$. We show:
$$(\tind_{D'}^D(A'):R_\ph)=(A':R_{{}'J^{\tWW_I}_{\tWW}(\ph)}).\tag h$$ 
We may assume that $\ph=\ph_E$ where $E\in\Irr_\boc(\tWW)$. By the definition of $\tind_{D'}^D(A')$, the left hand
side of (h) is equal to $(\ind_{D'}^D(A'):R_E)$. From the second equality in 43.9(a) we see that
$$R_{E|_{\tWW_I}}=\sum_{E'\in\un\Irr(\tWW_I)}\la E',E\ra R_{E'}.$$
By 43.9(b) we may restrict the previous sum to those $E'$ such that $a_{E'}\le a_E$; moreover for $E'$ such that 
$a_{E'}<a_E$ we have $\boc_{E'}\ne\boc'$. Thus we have $R_{E|_{\tWW_I}}=R_{{}'J^{\tWW_I}_{\tWW}(\ph)}$ plus a 
linear combination of $A''\in\hD'{}^{un}$ with $\boc'_{A''}\ne\boc'$. We see that the right hand side of (h) is 
equal to $(A':R_{E|_{\tWW_I}})$ hence to $(A':\res_D^{D'}(R_E))$ (see 44.12(e)) and (h) is equivalent to
$(\ind_{D'}^D(A'):R_E)=(A':\res_D^{D'}(R_E))$; but this follows from 44.12(d). This proves (h).

\subhead 44.21\endsubhead
We preserve the setup of 44.20. We assume that 

(i) for any $E'\in\Irr_{\boc'}(\tWW_I)$ there exists a unique $E\in\Irr_\boc(\tWW)$ (up to isomorphism) such that
$\la E',E\ra\ne0$; moreover we then have $\la E',E\ra=1$;

(ii) for any $E\in\Irr_\boc(\tWW)$ there exists a unique $E'\in\Irr_{\boc'}(\tWW_I)$ (up to isomorphism) such that
$\la E',E\ra\ne0$; moreover we then have $\la E',E\ra=1$;
\nl
Note that the $E'\m E$ in (i) and $E\m E'$ in (ii) defined inverse bijections $E'\lra E$ between the sets of 
isomorphism classes of objects in $\Irr_{\boc'}(\tWW_I)$ and $\Irr_\boc(\tWW)$. If $E'\lra E$ then 
$$\align&J_{\tWW_I}^{\tWW}(\ph_{E'})=\ph_E,\\&
{}'J^{\tWW_I}_{\tWW}(\ph_E)=\ph_{E'}+\text{ linear combination of elements $\ph_{E''}$ with}\\&
E''\in\Irr(\tWW_I)-\Irr_{\boc'}(\tWW_I).\tag a\endalign$$
The second equality in (a) is obvious. To prove the first equality in (a) we consider $\tE\in\Irr(\WW)$ such that
$a_{E'}=a_{\tE}$ and $\la E',\tE\ra\ne0$. It is enough to show that $\tE=E$. By 43.11(b) we have 
$\boc_{\tE}=\boc$. Using (i) we see that $\tE=E$, as required.

We show:
$$\text{ if }A'\in\hD'{}^{un}_{\boc'}\text{ then }\tind_{\tWW_I}^{\tWW}(A')\ne0.\tag b$$
Assume that $\tind_{\tWW_I}^{\tWW}(A')=0$. From 44.20(h) we deduce $(A':R_{{}'J^{\tWW_I}_{\tWW}(\ph_E)})=0$ for 
any $E\in\Irr_\boc(\tWW)$. Thus, for any $E'\in\Irr_{\boc'}(\tWW_I)$ we have $(A':R_{\ph_{E'}})=0$ (see (a)). This
contradicts 44.7(k) for $D'$. This proves (b).

We show:
$$\text{ if }A'\in\hD'{}^{un}_{\boc'}\text{ then }A:=\tind_{\tWW_I}^{\tWW}(A')\text{ is a single object of }
\hD^{un}_{\boc}.\tag c$$
By 44.7(k) we can find $E'\in\Irr(\tWW_I)$ such that $(A':R_{E'})\ne0$. We have necessarily 
$E'\in\Irr_{\boc'}(\tWW_I)$. By 43.10(b), $R_{E'}$ is a $\QQ$-linear combination of elements $R_{\ale^I_{x\vp}}$ 
such that $\tr(t_x\vp,E^\iy)\ne0$ (and in particular $x\in\boc'$). Hence there exists $x\in\boc'$ such that 
$(A':R_{\ale^I_{x\vp}})\ne0$. By 44.20(d) we have
$$R_{\ale^I_{x\vp}}=n_1A_1+n_2A_2+\do+n_rA_r\tag d$$ 
where $A_i\in\hD'{}^{un}$ are nonisomorphic to each other and $n_i\in\ZZ_{>0}$; we can assume that $A_1=A'$. We 
have:
$$\align&{}'J_{\tWW}^{\tWW_I}(\ale_{x\vp})=\ale^I_{x\vp}+
\text{ linear combination of elements $\ph_{E''}$ with}\\&E''\in\Irr(\tWW_I)-\Irr_{\boc'}(\tWW_I).\tag e\endalign
$$
Using (a) we see that this is equivalent to the identity $\tr(t_x\vp,E^\iy)=\tr(t_x\vp,E'{}^\iy)$ (for any 
$E'\lra E$ as above) which follows from 43.10(c). For $i,j$ in $[1,r]$ we set 
$x_{i,j}=(\tind_{\tWW_I}^{\tWW}(A_i):\tind_{\tWW_I}^{\tWW}(A_j))$. We have
$$\align&\sum_{i,j\in[1,r]}n_in_jx_{i,j}=
(\tind_{\tWW_I}^{\tWW}(R_{\ale^I_{x\vp}}):\tind_{\tWW_I}^{\tWW}(R_{\ale^I_{x\vp}}))=
(\tind_{\tWW_I}^{\tWW}(R_{\ale^I_{x\vp}}):R_{\ale_{x\vp}})\\&=
(R_{\ale^I_{x\vp}}:R_{'J_{\tWW}^{\tWW_I}(\ale_{x\vp})})=
(R_{\ale^I_{x\vp}}:R_{\ale^I_{x\vp}})=\sum_{i\in[1,r]}n_i^2.\endalign$$
(The first equality comes from (d); the second equality comes from 44.20(g); the third equality comes from 
44.20(h); the fourth equality comes from (e); the fifth equality comes from (d).) Since 
$\tind_{\tWW_I}^{\tWW}(A_i)$ is an $\NN$-linear combination of objects in $\hD^{un}$ and is $\ne0$ by (b), we see
that $(\tind_{\tWW_I}^{\tWW}(A_i):\tind_{\tWW_I}^{\tWW}(A_j))$ is $\ge1$ if $i=j$ and is $\ge0$ if $i\ne j$. Hence
from the equality $\sum_{i,j\in[1,r]}n_in_jx_{i,j}=\sum_{i\in[1,r]}n_i^2$ it follows that $x_{i,j}=1$ if $i=j$ and
$x_{i,j}=0$ if $i\ne j$. Since $A'=A_1$ we see that (c) holds.

We show:

(f) {\it If $A_1,A_2$ are objects of $\hD'{}^{un}_{\boc'}$ and
$A:=\tind_{\tWW_I}^{\tWW}(A_1)=\tind_{\tWW_I}^{\tWW}(A_2)$ then $A_1\cong A_2$.}
\nl
Assume that $A_1\not\cong A_2$. Let $E'\in\Irr_{\boc'}(\tWW_I)$. We can find $E\in\Irr_\boc(\tWW)$ such that 
$\la E',E\ra=1$. For $i=1,2$ we have
$$(A:R_E)=(\tind_{\tWW_I}^{\tWW}(A_i):R_E)=
(A_i:R_{{}'J^{\tWW_I}_{\tWW}(\ph_E)})=(A_i:R_{E'}).$$
(The second equality holds by 44.20(h); the third equality holds by (a).) Thus we have $(A_1:R_{E'})=(A_2:R_{E'})$
for any $E'\in\Irr_{\boc'}(\tWW_I)$. This implies that $(A_1:R_{\ale^I_{x\vp}})=(A_2:R_{\ale^I_{x\vp}})$ for any 
$x\in\WW_I$. We can choose $x\in\boc'$ such that $(A_1:R_{\ale^I_{x\vp}})\ne0$. Then we have also
$(A_2:R_{\ale^I_{x\vp}})\ne0$. We can assume that $A_1,A_2$ are the first two terms in the right hand side of (d).
But in the proof of (c) we have seen that $(\tind_{\tWW_I}^{\tWW}(A_1):\tind_{\tWW_I}^{\tWW}(A_2))=0$. This 
contradicts the assumption that $\tind_{\tWW_I}^{\tWW}(A_1)=\tind_{\tWW_I}^{\tWW}(A_2)$ which is $\ne0$ by (b).
This proves (f).

We show:

(g) {\it If $A\in\hD^{un}_{\boc}$ then there exists $A'\in\hD'{}^{un}_{\boc'}$ such that
$A=\tind_{\tWW_I}^{\tWW}(A')$.}
\nl
By 44.7(k) we can find $E\in\Irr(\tWW)$ such that $(A:R_E)\ne0$. We have necessarily $E\in\Irr_\boc(\tWW)$. Let 
$E'\in\Irr_{\boc'}(\tWW_I)$ be such that $E'\lra E$. By 44.20(f) we have 
$0\ne(A:R_E)=(A:R_{J_{\tWW_I}^{\tWW}(E')})=(A:\tind_{\tWW_I}^{\tWW}(R_{E'}))$. Hence there exists 
$A'\in\hD'{}^{un}_{\boc'}$ such that $(A':R_{E'})\ne0$ and $(A:\tind_{\tWW_I}^{\tWW}(A'))\ne0$. This implies that
$A=\tind_{\tWW_I}^{\tWW}(A')$. This proves (g).

Combining (c),(f),(g) and using 44.20(h) and (a), we obtain the following result:

(h) {\it $A'\m\tind_{\tWW_I}^{\tWW}(A')$ defines a bijection between the set of isomorphism classes in 
$\hD'_{\boc'}{}^{un}$ and the set of isomorphism classes in $\hD_\boc^{un}$; this bijection has the following
property: for any $E\in\Irr(\tWW)$ and any $A'\in\hD'_{\boc'}{}^{un}$ we have 
$(\tind_{\tWW_I}^{\tWW}(A'):R_E)=0$ if $E\n\Irr_\boc(\tWW)$ and $(\tind_{\tWW_I}^{\tWW}(A'):R_E)=(A':R_{E'})$ 
where $E'\in\Irr_{\boc'}(\tWW_I)$ is defined uniquely by $\la E',E\ra=1$.}

\head 45. Reductions\endhead
\subhead 45.1\endsubhead
In this section we show that the problem of classifying the unipotent character sheaves on $D$ can be reduced to
the analogous problem in the case where $G^0$ is simple and $\cz_G=\{1\}$.

Let $\t:G^0_{sc}@>>>G^0$ be a simply connected covering of the derived group of $G^0$. Let 
$\wt{G^0}=\cz_{G^0}^0\T G^0_{sc}$. The homomorphism $\ps:\wt{G^0}@>>>G^0$, $(z,g)\m z\t(g)$ is surjective with 
finite kernel which may be identified with $\{z\in\cz_{G^0_{sc}};\t(z)\in\cz_{G^0}^0\}$. Let $\fs(G^0)$ be the 
category whose objects are the local systems $\ce$ of rank $1$ on $G^0$ such that for some 
$\ce_0\in\fs(\cz_{G^0}^0)$ we have $\ps^*\ce\cong\ce_0\bxt\bbq$ or equivalently $\ce$ is a direct summand of the 
local system $\ps_!(\ce_0\bxt\bbq)$. (When $G^0$ is a torus this definition of $\fs(G^0)$ agrees with that in 
28.1.) Let $\ce\in\fs(G^0)$. We show:

(a) {\it $\ce$ is $G^0$-equivariant for the conjugation action of $G^0$ on $G^0$;}

(b) {\it $\ce$ is $G^0_{sc}\T G^0_{sc}$-equivariant for the $G^0_{sc}\T G^0_{sc}$-action on $G^0$ given by 
$(x_1,x_2):g\m\t(x_1)g\t(x_2\i)$;}

(c) {\it for any $x\in G^0$ we have $L_x^*\ce\cong\ce$ where $L_x:G^0@>>>G^0$ is given by $g\m xg$.}
\nl
Let $\ce_0\in\fs_n(\cz_{G^0}^0)$ be such that $\ce$ is a direct summand of $\ps_!(\ce_0\bxt\bbq)$. The 
$G^0$-action on $\wt{G^0}$ given by $y:(z,x)\m\ty(z,x)\ty\i$ (where $\ty\in\ps\i(y)$) is well defined and is 
compatible under $\ps$ with the conjugation action of $G^0$ on $G^0$; moreover, $\ce_0\bxt\bbq$ is 
$G^0$-equivariant. Hence $\ps_!(\ce_0\bxt\bbq)$ is $G^0$-equivariant and (a) holds. The
$G^0_{sc}\T G^0_{sc}$-action on $\wt{G^0}$ given by $(x_1,x_2):(z,x)\m(z,x_1xx_2\i)$ is compatible under $\ps$ 
with the $G^0_{sc}\T G^0_{sc}$-action on $G^0$ (as in (b)) and $\ce_0\bxt\bbq$ is 
$G^0_{sc}\T G^0_{sc}$-equivariant. Hence $\ps_!(\ce_0\bxt\bbq)$ is $G^0_{sc}\T G^0_{sc}$-equivariant and (b) 
holds. We prove (c). The $\wt{G^0}$-action on $\wt{G^0}$ given by $(z,x):(z',x')\m(z^nz',xx')$ is compatible under
$\ps$ with the $\wt{G^0}$-action on $G^0$ given by $(z,x):g\m z^n\t(x)g$ and $\ce_0\bxt\bbq$ is 
$\wt{G^0}$-equivariant. Hence $\ps_!(\ce_0\bxt\bbq)$ is $\wt{G^0}$-equivariant. Since the map $\wt{G^0}@>>>G^0$,
$(z,x)\m z^n\t(x)$ is surjective, we see that (c) holds.

Let $B^*,T$ be as in 28.5. Let $h:T@>>>G^0$ be the inclusion; let $\tT=\t\i(T)$ (a maximal torus of $G^0_{sc}$).
Let $\t_T:\tT@>>>T,\ps_T:\cz_{G^0}^0\T\tT@>>>T$ be the restrictions of $\t,\ps$. Let $\fs(T)^1$ be the category 
whose objects are the local systems $\ce'$ in $\fs(T)$ which satisfy one of the following four equivalent
conditions:

(i) for some $\ce_0\in\fs(\cz_{G^0}^0)$ we have $\ps_T^*\ce'\cong\ce_0\ot\bbq$;

(ii) $\ce'$ is a direct summand of the local system $\ps_{T!}(\ce_0\ot\bbq)$;

(iii) $\t_T^*\ce'\cong\bbq$;

(iv) for any coroot $f:\kk^*@>>>T$ we have $f^*\ce'\cong\bbq$.
\nl
From the definitions we see that 

(d) $\ce\m\ce_T:=h^*\ce$ is an equivalence of categories $\fs(G^0)@>>>\fs(T)^1$. 
\nl
Let $\fs(\TT)^1$ be the category whose objects are the local systems $\ce'$ in $\fs(\TT)$ such that 
$\cha^*\ce'\cong\bbq$ for any $\a\in R$ (see 28.3). We identify $T=\TT$ as in 28.5. Then $\fs(T)^1$ becomes 
$\fs(\TT)^1$.

\subhead 45.2\endsubhead
Let $d\in N_D(B^*)\cap N_D(T)$. There is a unique automorphism $\d_0:G^0_{sc}@>>>G^0_{sc}$ such that 
$\t(\d_0(g))=d\i\t(g)d$ for all $g\in G^0_{sc}$. Define an automorphism $\d:\wt{G^0}@>>>\wt{G^0}$ by 
$\d(z,g)=(d\i zd,\d_0(g))$. Then $\ps(\d(y))=d\i\ps(y)d$ for all $y\in\wt{G^0}$. 

Let $\ce\in\fs(G^0)$. Note that $\Ad(d\i)^*\ce\in\fs(G^0)$. Define $L_{d\i}:D@>>>G^0$ by $g\m d\i g$.
We set $\ce_D=L_{d\i}^*\ce$, a local system of rank 
$1$ on $D$. We show that the following three conditions are equivalent:

(i) $\Ad(d\i)^*\ce\cong\ce$;

(ii) $\Ad(d\i)^*\ce_T\cong\ce_T$;

(iii) the local system $\ce_D$ on $D$ is $G^0$-equivariant for the conjugation action of $G^0$ on $D$.
\nl
Now (i), (ii) are equivalent by 45.1(d); moreover if (i) or (iii) holds for some $d\in N_D(B^*)\cap N_D(T)$ then 
it holds for any $d\in N_D(B^*)\cap N_D(T)$ (by the $G^0$-equivariance of $\ce$, see 45.1(a)).

Assume first that (i) holds. Let $\tD=\{(y,x')\in D\T\wt{G^0};d\i y=\ps(x')\}$. Let $L':\tD@>>>\wt{G^0}$, 
$\ps':\tD@>>>D$ be the obvious projections. Let $\ce_0\in\fs(\cz_{G^0}^0)$ be such that 
$\ps^*\ce\cong\ce_0\bxt\bbq$. Then 
$$\Ad(d\i)^*\ce_0\bxt\bbq=\d^*(\ce_0\bxt\bbq)\cong\d^*\ps^*\ce\cong\ps^*\Ad(d\i)^*\ce\cong\ps^*\ce\cong
\ce_0\bxt\bbq,$$
hence $\Ad(d\i)^*\ce_0\cong\ce_0$. By 28.2, $\ce_0$ is $\cz_{G^0}^0$-equivariant for the $\cz_{G^0}^0$-action on 
$\cz_{G^0}^0$ given by $z_0:z\m d\i z_0dzz_0\i$. Hence $\ce_0\bxt\bbq$ is $\wt{G^0}$-equivariant for the 
$\wt{G^0}$-action on $\wt{G^0}$ given by $x:x'\m\d(x)x'x\i$. Define a $\wt{G^0}$-action on $\tD$ by 
$x:(y,x')\m(\ps(x)y\ps(x)\i,\d(x)x'x\i)$. This action is compatible under $\ps'$ with the $\wt{G^0}$-action on $D$
given by $x:y\m\ps(x)y\ps(x)\i$ and is compatible under $L'$ with the $\wt{G^0}$-action on $\wt{G^0}$ given by 
$x:x'\m\d(x)x'x\i$. It follows that $L'{}^*(\ce_0\bxt\bbq)$ is $\wt{G^0}$-equivariant and 
$\ps'_!L'{}^*(\ce_0\bxt\bbq)=L_{d\i}^*\ps_!(\ce_0\bxt\bbq)$ is $\wt{G^0}$-equivariant. Since $L_{d\i}^*\ce$ is a 
direct summand of $L_{d\i}^*\ps_!(\ce_0\bxt\bbq)$, we see that $L_{d\i}^*\ce$ is $\wt{G^0}$-equivariant. Since 
$\wt{G^0}$ acts on $D$ through its quotient $G^0$, we see that $\ker\ps$ acts naturally on the stalk of 
$L_{d\i}^*\ce$ at $y\in D$ through a character $\c$ which is independent of $y$. To show that $L_{d\i}^*\ce$ is 
$G^0$-equivariant it is enough to show that $\c=1$. Let $\tT=\ps\i(T)$, a maximal torus of $\wt{G^0}$. Then 
$L_{\d\i}^*\ce|_{dT}$ is $\tT$-equivariant (for the restriction of the $\wt{G^0}$-action to $\tT$). Since 
$\ker\ps\sub\tT$, $\c$ is determined by the $\tT$-equivariant structure of $L_{\d\i}^*\ce|_{dT}$. To show that 
$\c=1$ it is then enough to show that $L_{\d\i}^*\ce|_{dT}$ is $T$-equivariant for the conjugation $T$-action on 
$dT$. From (i) we deduce $\Ad(d\i)^*\ce_T\cong\ce_T$. By 28.2, $\ce_T$ is $T$-equivariant for the $T$-action on 
$T$ given by $t_0:t\m d\i t_0dtt_0\i$. Also $\l:dT@>>>T,dt\m t$ is compatible with the $T$-action on $T$ (as 
above) and the $T$-action on $dT$ given by conjugation. Hence $\l^*\ce_T$ is $T$-equivariant. Hence 
$L_{\d\i}^*\ce|_{dT}$ is $T$-equivariant. We see that (iii) holds.

Conversely, assume that (iii) holds. Then $m^*L_{d\i}^*\ce\cong m'{}^*L_{d\i}^*\ce$ where $m,m':G^0\T D@>>>D$ are 
given by $m(g,y)=gyg\i,m'(g,y)=y$. Define $j:G^0@>>>G^0\T D$ by $j(g)=(g,dg)$. Then 
$L_{d\i}mj=\Ad(d\i),L_{d\i}m'j=1$ hence $\Ad(d\i)^*\ce=j^*m^*L_{d\i}^*\ce\cong j^*m'{}^*L_{d\i}^*\ce=\ce$. We see
that (i) holds.

\subhead 45.3\endsubhead
Let $\ce\in\fs(G^0)$ and let $\cl=\ce_T\in\fs(T)^1$. Then $\uD\in\Wb_\cl$. Moreover, for any $w\in\WW$ we have
$w\in\Wb_\cl$ (see 45.1(iv) and 28.3(a)); hence $w\uD\in\Wb_\cl$. Hence the local system $\tcl$ on 
$Z_{\em,\II,D}^w$ is defined as in 28.7. From the definitions we see that $\tcl=\p_w^*\ce_D$ where
$\p_w:Z_{\em,\II,D}^w@>>>D$ is the map $(B,B',g)\m g$. Hence 
$$K^{w,\cl}_D=\p_{w!}\p_w^*\ce_D=\ce_D\ot\p_{w!}\p_w^*\bbq=\ce_D\ot\p_{w!}\bbq=\ce_D\ot K^w_D\in\cd(D),$$
(notation of 28.19). 

\subhead 45.4\endsubhead
Now let $\G$ be a closed normal subgroup of $G$ contained in $\cz_{G^0}$. Then $G'=G/\G$ is a reductive group and
the image $D'$ of $D$ under the obvious homomorphism $\o:G@>>>G'$ is a connected component $D'$ of $G'$ that 
generates $G'$. We may regard naturally $\G$ as a subgroup of the canonical torus $\TT$ of $G^0$ and we may 
identify naturally $\TT/\G$ with $\TT'$, the canonical torus of $G'$. Let $\WW'$ be the Weyl group of $G'{}^0$ and
let $\II'$ be its set of simple reflections (see 26.1). We identify $\WW'=\WW$, $\II'=\II$ in an obvious way. Then
$\WW$ acts on $\TT,\TT'$ compatibly with the canonical map $\TT@>>>\TT'$. Let $\o_D:D@>>>D'$ be the restriction of
$\o$.

Let $w\in\WW$. Then $K^w_D,\bK^w_D\in\cd(D)$, $K^w_{D'},\bK^w_{D'}\in\cd(D')$ are defined. We show
$$K^w_D\cong\o_D^*K^w_{D'}\in\cd(D),\bK^w_D\cong\o_D^*\bK^w_{D'}\in\cd(D).\tag a$$
Define $Z^w_{\em,\II,D'}$ in terms of $G'$ in the same way that $Z^w_{\em,\II,D}$ is defined in terms of $G$. Let
$\p_w:Z^w_{\em,\II,D}@>>>D$ be as in 45.3 and let $\p'_w:Z^w_{\em,\II,D'}@>>>D'$ be the analogous map defined in 
terms of $G'$. Define $\o':Z^w_{\em,\II,D}@>>>Z^w_{\em,\II,D'}$ by $(B,B',g)\m(\o(B),\o(B'),\o(g))$. We have a 
cartesian diagram
$$\CD
Z^w_{\em,\II,D}@>\o'>>Z^w_{\em,\II,D'}\\
@V\p_wVV                   @V\p'_wVV \\
D                  @>\o_D>>                 D'          \endCD$$
Hence 
$$\o_D^*K^w_{D'}=\o_D^*\p'_{w!}\bbq=\p_{w!}\bbq=K^w_D,$$
as required. The second statement in (a) is proved similarly. We set $r=\dim(\G)$. From (a) we deduce for any
$i\in\ZZ$:

(b) $H^i(K^w_D)\cong\o_D^*(H^{i-r}(K^w_{D'}))[r]$,
$H^i(\bK^w_D)\cong\o_D^*(H^{i-r}(\bK^w_{D'}))[r]$;

(c) {\it if $A'\in\hD'{}^{un}$ then the perverse sheaf $\o_D^*(A')[r]$ is a direct sum of finitely many objects of
$\hD^{un}$.}

\subhead 45.5\endsubhead
In the setup of 45.4 we assume that $\G=\cz_{G^0}^0$. Then $\o_D:D@>>>D'$ is a a fibration with smooth, connected 
fibres. Using this and 45.4(c) we see that if $A'\in\hD'{}^{un}$ then $\o_D^*(A')[r]\in\hD^{un}$ and (in the setup
of 45.4(b)): 
$$\align&(A':H^{i-r}(K^w_{D'}))=(\o_D^*(A')[r]:H^i(K^w_D)),\\&
(A':H^{i-r}(\bK^w_{D'}))=(\o_D^*(A')[r]:H^i(\bK^w_D)).\tag a\endalign$$
Now let $A\in\hD^{un}$. We show that $A\cong\o_D^*(A')[r]$ for some $A'\in\hD'{}^{un}$. We can find $w\in\WW$
and $i\in\ZZ$ such that $(A:H^i(K^w_D))>0$. By 45.4(b) we then have 
$(A:\o_D^*(H^{i-r}(K^w_{D'}))[r])>0$.
Hence there exists $A'\in\hD'{}^{un}$ such that $(A:\o_D^*(A')[r])>0$, as required. Note that if $A',A''$ are
objects of $\hD'{}^{un}$ such that $\o_D^*(A')[r]\cong\o_D^*(A'')[r]$ then $A'\cong A''$ (a standard 
property of $\o_D^*$). We see that $A'\m\o_D^*(A')[r]$ defines a bijection $\uhD'{}^{un}@>\si>>\uhD^{un}$.

Let $E\in\Irr(\tWW)$. Let $R_E\in\ck^{un}_\QQ(D)$ be as in 44.6(b) and let $R'_E\in\ck^{un}_\QQ(D')$ be
the analogous object defined in terms of $G'$. From (a) we see that for $A'\in\hD'{}^{un}$ we have
$$(A':R'_E)=(\o_D^*(A')[r]:R_E).\tag b$$
Moreover, since $\dim\supp(\o_D^*A[r])=\dim\supp(A)+r$, we see from (a) that:

(c) {\it if $D'$ has property $\tfA$ then $D$ has property $\tfA$.}

\subhead 45.6\endsubhead
In the setup of 45.4 we assume that $\cz_{G^0}^0=\{1\}$ so that $\G$ is a finite abelian group. Then 
$\cz_{G'{}^0}^0=\{1\}$. Let $\G^*=\Hom(\G,\bbq^*)$. For $\c\in\G^*$ define ${}^D\c\in\G^*$ by $x\m\c(dxd\i)$ 
(with $d\in N_D(B^*)\cap N_D(T)$). Let ${}^D\G^*=\{\c\in\G^*;{}^D\c=\c\}$. Let $\o_0:G^0@>>>G'{}^0$ be the 
restriction of $\o$. Since $\G$ is abelian, we have 
$$\o_{0!}\bbq\cong\op_{\c\in\G^*}\ce^\c\tag a$$
where $\ce^\c$ is a local system of rank $1$ on $G'{}^0$, equivariant for the $G^0$-action $g:g'\m\o_0(g)g'$ of 
$G^0$ on $G'{}^0$, which induces an action of $\G$ on any stalk of $\ce^\c$ through $\c$. Let $\ce^\c_{T'}$ be the
restriction of $\ce^\c$ to $T'$. Let $\ps'$ be the composition $G^0_{sc}@>\ps>>G^0@>\o_0>>G'{}^0$ ($\ps$ as in 
45.1). For $\c\in\G^*$ we have $\o_0^*\ce^\c\cong\bbq$ hence $\ps'{}^*\ce^\c\cong\bbq$ and 
$\ce^\c\in\fs(G'{}^0)$. Let $d'=\o(d)\in D'$. Define $L'_{d'{}\i}:D'@>\si>>G'{}^0$ by $g'\m d'{}\i g'$. For 
$\c\in\G^*$ we set $\ce^\c_{D'}=L'_{d'{}\i}{}^*\ce^\c$, a local system of rank $1$ on $D'$. From (a) we deduce
$$\o_{D!}\bbq\cong\op_{\c\in\G^*}\ce^\c_{D'}.\tag b$$
It follows that $\op_{\c\in\G^*}\ce^\c_{D'}$ is $G^0$-equivariant for the $G^0$-action 
$$g:g'\m\o_0(g)g'\o_0(g)\i\tag c$$
on $D'$. Hence for any $\c$, $\ce^\c_{D'}$ is $G^0$-equivariant for the action (c). Since the restriction of the 
action (c) to $\G$ is trivial, we see that (c) induces an action of $\G$ on the stalk of $\ce^\c_{D'}$ at 
$y\in D'$ through a character $\ti\c$ which is independent of $y$. Moreover, we have $\ti\c=1$ if and only if 
$\ce^\c_{D'}$ is $G'{}^0$-equivariant for the conjugation action of $G'{}^0$ on $D'$. By 45.2 (for $G'$ instead of
$G$), this last condition is equivalent to the condition that $\Ad(d'{}\i)^*\ce^\c\cong\ce^\c$ that is, to the 
condition that ${}^D\c=\c$. Thus we have $\ti\c=1$ if and only if ${}^D\c=\c$. We show:

(d) {\it if $A'\in\hD'$ and $\c\in\G^*$ satisfies ${}^D\c\ne\c$ then the simple perverse sheaf 
$A'_1:=\ce^\c_{D'}\ot A'$ is not in $\hD'$.}
\nl
Indeed, $A'$ is a $G'{}^0$-equivariant simple perverse sheaf (for the conjugation action of $G'{}^0$) and $A'_1$ 
is $G^0$-equivariant for the action (c) in such a way that the induced action of $\G$ on stalks is via the 
non-trivial character $\ti\c$. We see that $A'_1$ is not $G'{}^0$-equivariant for the conjugation action of
$G'{}^0$; (d) follows.

Let $w\in\WW$. We show:
$$\o_{D!}K^w_D=\op_{\c\in\G^*;{}^D\c=\c}K^{w,\ce^\c_{T'}}_{D'}
\op\op_{\c\in\G^*;{}^D\c\ne\c}\ce^\c_{D'}\ot K^w_{D'}.\tag e$$
Using the cartesian diagram in 45.4 we have
$$\align&\o_{D!}K^w_D=\o_{D!}\p_{w!}\bbq=\p'_{w!}\o'_!\bbq=\p'_{w!}\p'_w{}^*\o_{D!}\bbq
=\o_{D!}\bbq\ot(\p'_{w!}\p'_w{}^*\bbq)\\&=\o_{D!}\bbq\ot\p'_{w!}\bbq
=\o_{D!}\bbq\ot K^w_{D'}=\op_{\c\in\G^*}\ce^\c_{D'}\ot K^w_{D'}.\endalign$$
It remains to use 45.3 (for $G',T'$ instead of $G,T$).

We show:

(f) {\it if $A'\in\hD'{}^{un}$ and $\c\in\G^*$, ${}^D\c=\c$, $\c\ne1$ then $A'\n\hD'{}^{\ce^\c_{T'}}$.}
\nl
Indeed, if $A'\in\hD'{}^{\ce^\c_{T'}}$ then by 32.24 there exists $a\in\WW$ such that 
$\bbq=a^*\bbq\cong\ce^\c_{T'}$ as local systems on $T'=\TT'$. Using 45.1(d) (for $G'$ instead of $G$) it follows 
that $\ce^\c=\ce^1$ hence $\c=1$, a contradiction.

We show that for any $A'\in\hD'{}^{un}$ and $i\in\ZZ$ we have
$$(A':\o_{D!}H^i(K^w_D))=(A':H^i(K^w_{D'})).\tag g$$
We use that $\o_{D!}H^i(K^w_D)=H^i(\o_{D!}K^w_D)$ which holds since $\o_D$ is a finite covering. Hence the left 
hand side of (g) can be rewritten using (e) as 
$$\sum\op_{\c\in\G^*;{}^D\c=\c}(A':H^i(K^{w,\ce^\c_{T'}}_{D'}))+
\sum_{\c\in\G^*;{}^D\c\ne\c}(A':\ce^\c_{D'}\ot H^i(K^w_{D'})).$$
The term corresponding to $\c$ such that ${}^D\c\ne\c$ is $0$ by (d); the term corresponding to $\c$ such that 
${}^D\c=\c$, $\c\ne1$ is $0$ by (f) and (g) follows. 

Using 45.4(b) we can reformulate (g) as follows:
$$(A':\o_{D!}\o_D^*(H^i(K^w_{D'}))=(A':H^i(K^w_{D'})).\tag h$$
In $\ck^{un}(D')$ we have $H^i(K^w_{D'})=\sum_{j=1}^sm_jA'_j$ where $A'_1,A'_2,\do,A'_s$ are mutually 
non-isomorphic objects in $\hD'{}^{un}$ and $m_j\in\ZZ_{>0}$. Applying (h) with $A'=A'_h$ we obtain
$\sum_{j=1}^sm_j(A'_h:\o_{D!}\o_D^*(A'_j))=m_h$ hence $\sum_{j=1}^sm_j(\o_D^*(A'_h):\o_D^*(A'_j))=m_h$ for 
$h\in[1,s]$. Since $(\o_D^*(A'_h):\o_D^*(A'_j))\ge\d_{h,j}$ it follows that $(\o_D^*(A'_h):\o_D^*(A'_j))=\d_{h,j}$
for $h,j\in[1,s]$. It follows that the perverse sheaf $\o_D^*A'_j$ is simple. Since any 
$A'\in\hD'{}^{un}$ appears in some $H^i(K^w_{D'})$ we see that in our case we have the following refinement of
45.4(c):

(i) {\it if $A'\in\hD'{}^{un}$ then $\o_D^*(A')\in\hD^{un}$.}
\nl
Now let $A\in\hD^{un}$. Let $\o_{D!}^0A$ be the sum of all simple summands of the semisimple perverse sheaf
$\o_{D!}A$ which are in $\hD'{}^{un}$. We show that:

(j) {\it $\o_{D!}^0A\in\hD'{}^{un}$.}
\nl
We can find $w\in\WW$ and $i\in\ZZ$ such that $A$ appears in $H^i(K^w_D)$. Using 45.4(b) we see that $A$ appears 
in $\o_D^*(H^i(K^w_{D'}))$. Hence there exists $C\in\hD'{}^{un}$ which appears in $H^i(K^w_{D'})$ such that 
$(A:\o_D^*C)>0$. By (i), $\o_D^*C$ is a simple perverse sheaf. It follows that $A\cong\o_D^*C$. Thus $C$
appears in $\o_{D!}A$. In particular, $\o_{D!}^0A\ne0$. Now assume that $C,C'$ are two objects in $\hD'{}^{un}$ 
such that both $C$ and $C'$ appear in $\o_{D!}A$. Then $A\cong\o_D^*C$; similarly, $A\cong\o_D^*C'$. Thus the 
simple objects $\o_D^*C,\o_D^*C'$ are isomorphic. It follows that $\dim\Hom(C',\o_{D!}\o_D^*C)=1$. We have 
$$\o_{D!}\o_D^*C=C\ot\o_{D!}\o_D^*\bbq=C\ot\o_{D!}\bbq=\op_{\c\in\G^*}C\ot\ce^\c_{D'}.$$
It follows that for some $\c\in\G^*$ we have $\dim\Hom(C',C\ot\ce^\c_D)=1$ hence $C'\cong C\ot\ce^\c_{D'}$. This 
forces ${}^D\c=\c$, by (d). Then $\ce^\c_{T'}$ is defined and from 45.3 we see that
$C\ot\ce^\c_{D'}\in\hD'{}^{\ce^\c_{T'}}$ so that $C'\in\hD'{}^{\ce^\c_{T'}}$. Using (f) we deduce that $\c=1$ and
$C'\cong C$. Thus, the semisimple perverse sheaf $\o_{D!}^0A$ is nonzero and isotypic. If $C\in\hD'{}^{un}$ 
appears in $\o_{D!}^0A$ then, as we have seen, we have $A\cong\o_D^*C$ hence $\dim\Hom(C,\o_{D!}A)=1$ so that
$\dim\Hom(C,\o_{D!}^0A)=1$. Thus $\o_{D!}^0A$ is simple. This proves (j).

From (i),(j) and the proof of (j) we see that:

(k) {\it $A'\m\o_D^*(A')$ defines a bijection $\uhD'{}^{un}@>\si>>\uhD^{un}$; the inverse bijection is induced by
$A\m\o_{D!}^0A$.}
\nl
We define $\tWW'$ in terms of $G',D'$ in the same way as $\tWW$ was defined in terms of $G,D$. We may assume that
$\tWW'=\tWW$. Let $E\in\Irr(\tWW)$. Let $R_E\in\ck^{un}_\QQ(D)$ be as in 44.6(b) and let 
$R'_E\in\ck^{un}_\QQ(D')$ be the analogous object defined in terms of $G'$. From (g) we see that for 
$A'\in\hD'{}^{un}$ we have
$$(A':R'_E)=(\o_D^*(A'):R_E).\tag l$$
If $A\in\hD^{un}$, $w\in\WW$, $i\in\ZZ$ then
$$(A:H^i(\bK^w_D))=(A:\o_D^*H^i(\bK^w_{D'}))=(\o_{D!}A:H^i(\bK^w_{D'}))=(\o_{D!}^0A:H^i(\bK^w_{D'})).$$
Since $A=\o_D^*(\o_{D!}^0A)$ we have $\dim\supp(A)=\dim\supp(\o_{D!}^0A)$. We see that 

(m) {\it if $D'$ has property $\tfA$ then $D$ has property $\tfA$.}

\subhead 45.7\endsubhead
In the setup of 45.4 assume that $\G=\cz_{G^0}$. Then $A'\m\o_D^*(A')[r]$ defines a bijection 
$\uhD'{}^{un}@>\si>>\uhD^{un}$. Moreover, for any $w\in\WW$, any $A'\in\hD'{}^{un}$ and any $i\in\ZZ$ we have
$$(A':H^{i-r}(K^w_{D'}))=(\o_D^*(A')[r]:H^i(K^w_D)).\tag a$$
Note that $G/\cz_{G^0}$ can be obtained from $G$ in two steps: we first form $G_1=G/\cz_{G^0}^0$ which has
$\cz_{G_1^0}^0=\{1\}$ and then we have $G/\cz_{G^0}=G_1/\cz_{G_1^0}$. We use 45.5 to compare $G$ to $G_1$ and
45.6(k),(h) to compare $G_1$ to $G/\cz_{G^0}$. The statements above follow.

We define $\tWW'$ in terms of $G',D'$ in the same way as $\tWW$ was defined in terms of $G,D$. We may assume that
$\tWW'=\tWW$. Let $E\in\Irr(\tWW)$. Let $R_E\in\ck^{un}_\QQ(D)$ be as in 44.6(b) and let 
$R'_E\in\ck^{un}_\QQ(D')$ be the analogous object defined in terms of $G'$. From (a) we see that for 
$A'\in\hD'{}^{un}$ we have
$$(A':R'_E)=(\o_D^*(A'):R_E).\tag b$$
Combining 45.5(c), 45.6(m) we see that

(c) {\it if $D'$ has property $\tfA$ then $D$ has property $\tfA$.}
\nl
Now if $A'\in\hD'{}^{un}$ then $A'$ is cuspidal if and only if $\o_D^*(A')[r]$ is cuspidal. It follows that 

(d) {\it if $D'$ has property $\fA_0$ then $D$ has property $\fA_0$.}

\subhead 45.8\endsubhead
Assume now that $\cz_{G^0}=\{1\}$. Let $\D=\cz_G$. Let $G'=G/\D$.

If $g\in G$ satisfies $gg_1=g_1g\mod\cz_G$ for any $g_1\in G$ then for any $g_1\in G$ we have
$gg_1g\i g_1\i\in G^0$ (since $G/G^0$ is abelian) hence $gg_1g\i g_1\i\in G^0\cap\cz_G\sub\cz_{G^0}=\{1\}$; thus,
$g\in\cz_G$. We see that $\cz_{G'}=\{1\}$.

Let $\p:G@>>>G'$ be the obvious map. Then $\p$ induces an isomorphism $G^0@>\si>>G'{}^0$ and an isomorphism of $D$
onto a connected component $D'$ of $G'$ which generates $G'$. We identify the canonical tori and Weyl groups of 
$G^0,G'{}^0$ in the obvious way. 

Let $w\in\WW$. From the definitions it is clear that
$$K^w_D=\p^*K^w_{D'}, \bK^w_D=\p^*\bK^w_{D'}.\tag a$$
It follows that 

(b) {\it $A'\m\p^*A'$ induces a bijection $\uhD'{}^{un}@>\si>>\uhD^{un}$;}
\nl
moreover, if $w\in\WW$, $A'\in\hD'{}^{un}$ and $i\in\ZZ$ then 
$$(A':H^i(K^w_{D'}))=(\p^*A':H^i(K^w_D)).\tag c$$
Let $E\in\Irr(\tWW)$. Let $R_E\in\ck^{un}_\QQ(D)$ be as in 44.6(b) and let $R'_E\in\ck^{un}_\QQ(D')$ be
the analogous object defined in terms of $G'$. From (c) we see that for $A'\in\hD'{}^{un}$ we have
$$(A':R'_E)=(\p^*A':R_E).\tag d$$
From the definitions we see that

(e) {\it if $D'$ has property $\tfA$ then $D$ has property $\tfA$;}

(f) {\it if $D'$ has property $\fA_0$ then $D$ has property $\fA_0$.}

\subhead 45.9\endsubhead
Assume now that $\cz_G=\{1\}$ with $G^0$ adjoint. We have $G^0=\prod_{f\in\fF}G_f$ where $\fF$ is a finite set and
$G_f$ ($f\in\fF$) are the maximal connected simple closed subgroups of $G^0$. There is a well defined permutation
$\io:\fF@>\si>>\fF$ such that $gG_fg\i=G_{\io(f)}$ for all $g\in D,f\in\fF$. Let $\bfF$ be the set of orbits of 
$\io$ on $\fF$. For any $\co\in\bfF$ we set $G_\co=\prod_{f\in\co}G_f$. Then $G_\co$ is a closed connected normal
subgroup of $G$; hence we have a well defined homomorphism $\th_\co:G@>>>\Aut(G_\co)$ given by $g:x\m gxg\i$. The
image of $\th_\co$ is denoted by $\tG_\co$. Since $G_\co$ is adjoint, $\tG_\co$ is a reductive grup with identity
component $G_\co$; it is generated by its connected component $D_\co:=\th_\co(D)$. 

Let $\bg\in\cz_{\tG_\co}$. We have $\bg=\th_\co(g)$ with $g\in G$ and $ygxg\i y\i=gyxy\i g\i$ (that is 
$y\i g\i ygx=xy\i g\i yg$) for any $y\in G_\co$. Thus $y\i g\i yg$ (an element of $G_\co$) is in the centre of 
$G_\co$ so that $y\i g\i y g=1$ for any $y\in G_\co$. We see that $\th_\co(g\i)=1$ that is $\bg\i=1$. Thus, 
$\cz_{\tG_\co}=\{1\}$.

Note that the homomorphism $G@>>>\prod_{\co\in\bfF}\tG_\co$ given by $(\th_\co)_{\co\in\bfF}$ is an imbedding of 
reductive groups by which we can identify the identity components $G^0=\prod_{\co\in\bfF}G_\co$ and the component
$D$ with the component $\prod_{\co\in\bfF}D_\co$.

We can identify $\WW=\prod_{\co\in\bfF}\WW_\co$ where $\WW_\co$ is the Weyl group of $G_\co$. Let $w\in\WW$ and 
let $w_\co$ be the $\WW_\co$-component of $w$. From the definitions we have
$$K^w_D=\bxt_{\co\in\bfF}K^{w_\co}_{D_\co},\qua\bK^w_D=\bxt_{\co\in\bfF}\bK^{w_\co}_{D_\co}.\tag a$$
Hence for $i\in\ZZ$ we have
$$\align&H^i(K^w_D)=\op_{(i_\co);\sum_\co i_\co=i}\bxt_{\co\in\bfF}H^{i_\co}(K^{w_\co}_{D_\co})\\&
H^i(\bK^w_D)=\op_{(i_\co);\sum_\co i_\co=i}\bxt_{\co\in\bfF}H^{i_\co}(\bK^{w_\co}_{D_\co}).\tag b\endalign$$
Assume that $A_\co\in\hD_\co^{un}$ is given for each $\co\in\bfF$. Let $A=\bxt_{\co\in\bfF}A_\co$, a simple
perverse sheaf on $D$. We can find $w=(w_\co)\in\WW$ and $(i_\co)\in\NN^{\bfF}$ such that
$(A_\co:H^{i_\co}(\bK^{w_\co}_{D_\co}))>0$ for all $\co$ hence
$(A:\bxt_{\co\in\bfF}H^{i_\co}(\bK^{w_\co}_{D_\co}))>0$.
Using (b) we deduce that $(A:H^i(\bK^w_D))>0$ where $i=\sum_\co i_\co$. Hence $A\in\hD^{un}$.

Conversely, let $A\in\hD^{un}$. We can find $w=(w_\co)\in\WW$ and $(i_\co)\in\NN^{\bfF}$ such that
$(A:H^i(\bK^w_D))>0$. Using (b) we deduce that $(A:\bxt_{\co\in\bfF}H^{i_\co}(\bK^{w_\co}_{D_\co}))>0$ for some 
$(i_\co)\in\NN^{\bfF}$ such that $i=\sum_\co i_\co$. Hence there exist $A_\co\in\hD_\co^{un}$ ($\co\in\bfF$) such
that $(A_\co:H^{i_\co}(\bK^{w_\co}_{D_\co}))>0$ and $A\cong\bxt_{\co\in\bfF}A_\co$. We see that
$$(A_\co)\m\bxt_{\co\in\bfF}A_\co\text{ induces a bijection }\prod_{\co\in\bfF}\uhD_\co^{un}@>\si>>\uhD^{un}.
\tag c$$
Moreover if $(A_\co)\lra A$ under this bijection then
$$(A:H^i(K^w_D))=\sum_{(i_\co);\sum_\co i_\co=i}\prod_{\co\in\bfF}(A_\co:H^{i_\co}(K^{w_\co}_{D_\co})).\tag d$$
For $\co\in\bfF$ we define $\tWW_\co$, $\Irr(\tWW_\co)$ in terms of $\tG_\co$ in the same way as $\tWW$, 
$\Irr(\tWW)$ were defined in terms of $G$ (see 43.1). For each $\co\in\bfF$ we assume given an object 
$E_\co\in\Irr(\tWW_\co)$. Then the vector space $E=\ot_\co E_\co$ can be naturally regarded as an object of 
$\Irr(\tWW)$. (Any object of $\Irr(\tWW)$ can be obtained in this way.) Define 
$R_{E_\co}\in\ck^{un}_\QQ(D_\co)$ in terms of $\tG_\co$ in the same way as $R_E$ was defined in terms of $G$. Let
$(A_\co)\lra A$ be as above. From (d) we see that for $A'\in\hD'{}^{un}$ we have
$$(A:R_E)=\prod_{\co\in\bfF}(A_\co:R_{E_\co}).\tag e$$
From the definitions we see that

(f) {\it if $D_\co$ has property $\tfA$ for any $\co$ then $D$ has property $\tfA$;}

(g) {\it if $D_\co$ has property $\fA_0$ for any $\co$ then $D$ has property $\fA_0$}.

\subhead 45.10\endsubhead
Let $x,x',y\in\WW$ be such that $x'=yx\e(y)\i$. We show
$$gr_1(K^x_D)=gr_1(K^{x'}_D)\in\ck^{un}(D).\tag a$$
The proof is similar to that in \cite{\DL, 1.6}. Arguing by induction on $l(y)$ we see that we may assume that 
$y=s\in\II$. 

Assume first that $l(x)=l(x')=l(sx)+1$. Define an isomorphism $Z^x_{\em,\II,D}@>>>Z^{x'}_{\em,\II,D}$ by 
$(B,B',g)\m(B_1,B'_1,g)$ where $B_1,B'_1\in\cb$ are given by $\po(B,B_1)=s$, $\po(B_1,B')=sx$, $B'_1=gB_1g\i$. (We
then have $\po(B_1,B'_1)=(sx)\e(s)=x'$.) It follows that $K^x_D=K^{x'}_D$. 

The case where $l(x)=l(x')=l(sx')+1$ can be reduced to the previous case by exchanging $x,x'$.

Assume next that $l(x')=l(x)+2$. If $(B,B',g)\in Z^{x'}_{\em,\II,D}$ then there are well defined $B_1,B'_1$ in 
$\cb$ such that $\po(B,B_1)=s$, $\po(B_1,B'_1)=x,\po(B'_1,B')=\e(s)$. We partition $Z^{x'}_{\em,\II,D}$ into two
pieces $Z',Z''$ (one closed, one open) defined respectively by the conditions $B'_1=gB_1g\i$, $B'_1\ne gB_1g\i$. 
Let $K',K''$ be the direct image with compact support of $\bbq$ under the maps $Z'@>>>D$, $Z''@>>>D$,
$(B,B',g)\m g$. Then $gr_1(K^{x'}_D)=gr_1(K')+gr_1(K'')$.
Now $(B,B',g)\m(B_1,B'_1,g)$ defines an affine line bundle $Z'@>>>Z^x_{\em,\II,D}$. Hence $gr_1(K')=gr_1(K^x_D)$.
It remains to show that $gr_1(K'')=0$. Let $\tZ$ be the set of all
$(B,B_0,B'_0,B',g)$ in $\cb^4\T D$ such that $\po(B,B_0)=s$, $\po(B_0,B'_0)=x\e(s)$, $gBg\i=B'$, $gB_0g\i=B'_0$. 
If $(B,B_0,B'_0,B',g)\in\tZ$ there is a unique $\tB\in\cb$ such that $\po(B_0,\tB)=x$, $\po(\tB,B'_0)=\e(s)$. We 
partition $\tZ$ into two subsets $\tZ_1,\tZ_2$ (one closed, one open) defined respectively by the conditions 
$\tB=B'$, $\tB\ne B'$. Let $\tK,K_1,K_2$ be the direct image with compact support of $\bbq$ under the maps 
$\tZ@>>>D$, $\tZ_1@>>>D$, $\tZ_2@>>>D$, $(B,B_0,B'_0,B',g)\m g$. We have $gr_1(\tK)=gr_1(K_1)+gr_1(K_2)$.
Now $(B,B_0,B'_0,B',g)\m(B_0,B'_0,g)$ is an isomorphism $\tZ_1@>>>Z^{x\e(s)}_{\em,\II,D}$ and an affine line 
bundle $\tZ@>>>Z^{x\e(s)}_{\em,\II,D}$; hence $\tK=K_1$ and $gr_1(K_2)=0$. Moreover, 
$(B,B_0,B'_0,B',g)\m(B,B',g)$ is an isomorphism $\tZ_2@>>>Z''$. Hence $K_2=K''$ and $gr_1(K'')=0$, as required.

The case where $l(x)=l(x')+2$ can be reduced to the previous case by exchanging $x,x'$.
It remains to consider the case where $l(x)=l(x')=l(sx)-1=l(sx')-1$. In this case we have $x=x'$ (see 
\cite{\DL, 1.6.4}) and there is nothing to prove.

\subhead 45.11\endsubhead
Assume now that $\cz_G=\{1\}$, that $G^0$ is adjoint $\ne\{1\}$ and that $G$ has no closed connected normal 
subgroups other than $G^0$ and $\{1\}$. Let $e$ be a pinning (or \'epinglage, see 1.6) of $G^0$ which projects to
$(B^*,T)$ under the map $p$ in 1.6. By the adjointness of $G^0$ there is a unique element $d\in D$ such that 
$\Ad(d):G^0@>>>G^0$ stabilizes $e$ under the action 1.6(i). We have $G^0=\prod_{f\in\fF}G_f$ as in 45.9. Let 
$\io:\fF@>>>\fF$, $\bfF$ be as in 45.9. If $\co\in\bfF$ then $G_\co$ (as in 45.9) is a closed connected normal 
subgroup of $G$ other than $\{1\}$ hence it is equal to $G^0$. Thus, we have $\co=\fF$ that is, $\io:\fF@>>>\fF$ 
has a single orbit. Let $k=|\fF|$. We can identify $\fF=\ZZ/k\ZZ$ in such a way that $\io(j)=j+1$ for any 
$j\in\ZZ/k\ZZ$. 

For $j\in\ZZ/k\ZZ$ let $\cb_j$ be the variety of Borel subgroups of $G_j$. We can identify
$\cb=\prod_{j\in\ZZ/k\ZZ}\cb_j$ by $B\lra(B_0,B_1,\do,B_{k-1})$ where $B\in\cb,B_j\in\cb_j$ satisfy
$B=\prod_{j\in\ZZ/k\ZZ}B_j$. In particular we have $B^*=\prod_{j\in\ZZ/k\ZZ}B^*_j$ where $B^*_j$ is a Borel
subgroup of $G_j$. We also have $T=\prod_{j\in\ZZ/k\ZZ}T_j$, where $T_j$ is a maximal torus of $B^*_j$. We can view
$e$ as a collection $(e_j)_{j\in\ZZ/k\ZZ}$ where $e_j$ is a pinning of $G_j$ which projects to $(B^*_j,T_j)$. Note
that $\Ad(d)$ carries $e_j$ to $e_{j+1}$ for any $j\in\ZZ/k\ZZ$.

We can identify $\WW=\prod_{j\in\ZZ/k\ZZ}\WW_j$, where $\WW_j$ is the Weyl group of $G_j$ and
$\II=\sqc_{j\in\ZZ/k\ZZ}\II_j$ where $\II_j$ is the set of simple reflections in $\WW_j$. Recall that 
$\e:\WW@>>>\WW$ is the automorphism induced by $\Ad(d):G^0@>>>G^0$. We have $\e(\WW_j)=\WW_{j+1}$ for 
$j\in\ZZ/k\ZZ$.

Now $d^k$ normalizes $G_0$ and $\Ad(d^k):G_0@>>>G_0$ stabilizes $e_0$. Let $G'$ be the subgroup of $G$ generated 
by $G_0$ and $d^k$. Since $d$ has finite order, $G'$ is closed, $G'{}^0=G_0$ and $D'=d^kG_0$ is a connected 
component of $G'$ that generates $G'$. 

We show that $\cz_{G'}=\{1\}$. If $g'\in\cz_{G'}$ then we have $g'=d^{kr}x$ for some $r\in\ZZ,x\in G_0$ and
$\Ad(g'):G_0@>>>G_0$ is the identity map hence $\Ad(g')$ stabilizes $e_0$. Since $\Ad(d^{kr})$ also stabilizes 
$e_0$ we see that $\Ad(x)$ stabilizes $e_0$. Since $G_0$ is adjoint we must have $x=1$ hence $g'=d^{kr}$. Thus 
$g'$ commutes with $d$. Since $g'$ also centralizes $G_0$ and $d,G_0$ generate $G$ we see that $g'$ centralizes 
$G$ hence $g'=1$ (by our assumption that $\cz_G=\{1\}$). This verifies our assertion.

Define $\b:D@>>>D'$ by $\b(dg_0g_1\do g_{k-1})=dg_{k-1}dg_{k-2}\do dg_0$ where $g_j\in G_j$ or equivalently by the
requirement that $\z^k\in\b(\z)G_1G_2\do G_{k-1}$ for $\z\in D$. This is a principal 
$\{1\}\T G_1\T G_2\T\do\T G_{k-1}$-bundle where this group acts on $D$ by restriction of the conjugation action of
$G^0$. Moreover, $\b$ is compatible with the conjugation action of $G^0$ on $D$ and the conjugation action of 
$G_0$ on $D'$ via the homomorphism $G^0@>>>G_0$ which takes $g_0$ to $g_0$ if $g_0\in G_0$ and $g_i$ to $1$ if 
$i\in[1,k-1]$. We see that (setting $t=(k-1)\dim G_0$):

(a) {\it $A'\m\b^*A'[t]$ is an equivalence between the category of $G_0$-equivariant perverse sheaves
on $D'$ and the category of $G^0$-equivariant perverse sheaves on $D$.}
\nl
Let $w\in\WW_0\sub\WW$. The variety $Z^w_{\em,\II,D}$ may be identified with
$$\align&\{((B_0,B_1,\do,B_{k-1}),(B'_0,B'_1,\do,B'_{k-1}),dg_0g_1\do g_{k-1});B_j,B'_j\in\cb_j,
g_j\in G_j,\\&B'_j=\Ad(dg_{j-1})B_{j-1}) (j\in\ZZ/k\ZZ),\\&\po(B_0,B'_0)=w,B_j=B'_j(j\ne0)\}\endalign$$
or with
$$\align&\{(B_0,B'_0,dg_0g_1\do g_{k-1});\\&B_0,B'_0\in\cb_0,g_j\in G_j,
B'_0=\Ad(dg_{k-1}dg_{k-2}\do dg_0)B_0,\po(B_0,B'_0)=w\}.\endalign$$
We see that we have a cartesian diagram
$$\CD
Z^w_{\em,\II,D} @>\ti\b>>Z^w_{\em,\II_0,D'}\\
@VVV                     @VVV           \\
D          @>\b>>          D' \endCD$$
where 
$$\align&\ti\b:(B_0,B_1,\do,B_{k-1}),(B'_0,B'_1,\do,B'_{k-1}),dg_0g_1\do g_{k-1})\\&
\m(B_0,B'_0,dg_{k-1}dg_{k-2}\do dg_0).\endalign$$
Using this cartesian diagram we see that $K^w_D=\b^*K^w_{D'}$. Similarly we have 
$\bK^w_D=\b^*\bK^w_{D'}$. Since $\b$ is smooth with connected fibres we see that for any $i\in\ZZ$ we have 
$$H^i(K^w_D)=\b^*H^{i-t}(K^w_{D'})[t],H^i(\bK^w_D)=\b^*H^{i-t}(\bK^w_{D'})[t]$$
and
$$\align&(\b^*A'[t]:H^i(K^w_D))=(A':H^{i-t}(K^w_{D'})),
\\&(\b^*A'[t]:H^i(\bK^w_D))=(A':H^{i-t}(\bK^w_{D'}))\tag b\endalign$$
for any simple perverse sheaf $A'$ on $D'$. From (b) we see that, if $A'\in\hD'{}^{un}$, then 
$\b^*A'[t]\in\hD^{un}$.

Conversely, assume that $A\in\hD^{un}$. Let $\cx$ be the set of sequences \lb
$\ss=(s_1,s_2,\do,s_r)$ in $\II$ such
that $(A:H^i(K^\ss_D))>0$ for some $i$. Let $\cx_0$ be the set of all $\ss=(s_1,s_2,\do,s_r)\in\cx$ such 
that $s_h\in\II_0$ for all $h$. Note that $\cx\ne\em$. Let $N$ be the minimum value of 
$N_\ss:=\sum_{j\in[0,k-1],h\in[1,r];s_h\in\II_j}j$ where $\ss=(s_1,s_2,\do,s_r)$ runs through $\cx$. 

Assume that $N>0$. We choose $\ss\in\cx$ such that $N_\ss=N$. We can find $h\in[1,r]$ such that $s_h\in\II_j$ for
some $j\in[1,k-1]$; moreover we can assume that $h$ is maximum possible with this property. Then $s_{h'}\in\II_0$
for $h'\in[h+1,r]$. Let $\ss'=(s_1,s_2,\do,s_{h-1},s_{h+1},\do,s_r,s_h)$. Since $s_hs_{h'}=s_{h'}s_h$ for 
$h'\in[h+1,r]$ we see using the definitions that $K^\ss_D=K^{\ss'}_D$. Thus $\ss'\in\cx$. Note that 
$N_{\ss'}=N$. Let $\ss''=(\e\i(s_h),s_1,s_2,\do,s_{h-1},s_{h+1},\do,s_r)$. By 28.16 we have 
$K^{\ss'}_D=K^{\ss''}_D$. Thus $\ss''\in\cx$. Since $s_h\in\II_j$ with $j\in[1,k-1]$ we have 
$\e\i(s_h)\in\II_{j-1}$. Thus $N_{\ss''}=N_{\ss'}-1=N-1$. This contradicts the minimality of $N$. We have shown 
that $N=0$. We choose $\ss\in\cx$ such that $N_\ss=0$. We then have $\ss\in\cx_0$. Thus we have $\cx_0\ne\em$.

By the proof of the implication (iii)$\imp$(i) in 28.13 we deduce that there exists $w\in\WW_0$ and $i\in\ZZ$
such that $(A:H^i(K^w_D))>0$. Using (a) we can write $A=\b^*A'[t]$ where $A'$ is a well defined simple 
$G_0$-equivariant perverse sheaf on $D'$. Using (b) we see that $(A':H^{i-t}(K^w_{D'}))>0$. Hence
$A'\in\hD'{}^{un}$. Thus:

(c) {\it $A'\m\b^*A'[t]$ induces a bijection $\hD'{}^{un}@>\si>>\hD^{un}$.}
\nl
We define $\tWW'$ in terms of $G',D'$ in the same way as $\tWW$ was defined in terms of $G,D$; let $\vp'$ be the
element of $\tWW'$ which plays the same role for $\tWW'$ as $\vp$ for $\tWW$. We can assume that
the order of $\vp'$ in $\tWW'$ is the same as
the order of $\vp$ in $\tWW$. Let $E'\in\Irr(\tWW')$. Then the vector space $E=E'\ot E'\ot\do\ot E'$ ($k$ factors)
can be regarded as an object of $\Irr(\tWW)$ with $x=(x_0,x_1,\do,x_{k-1})$ ($x_j\in\WW_j$) acting by
$$e'_0\ot e'_1\ot\do\ot e'_{k-1}\m x_0(e'_0)\ot\e\i(x_1)(e'_1)\ot\do\ot\e^{-k+1}(x_{k-1})(e'_{k-1})$$ and $\vp$ 
acting by $e'_0\ot e'_1\ot\do\ot e'_{k-1}\m\vp'(e'_{k-1})\ot e'_0\ot\do\ot e'_{k-2}$. (Note that any object of 
$\Irr(\tWW)$ can be obtained in this way.) Define $R_{E'}\in\ck^{un}_\QQ(D')$ in terms of $G'$ in the same way as
$R_E$ was defined in terms of $G$. We show that for $A'\in\hD'{}^{un}$ we have
$$(\b^*A'[t]:R_E)=(A':R_{E'}).\tag d$$
Let $A=\b^*A'[t]$. Using (b) we see that the right hand side of (d) equals
$$\align&|\WW_0|\i\sum_{x\in\WW_0,i\in\ZZ}(-1)^{\dim G'+i}\tr(x\vp',E')(A':H^i(K^x_{D'}))\\&
=|\WW_0|\i\sum_{x\in\WW_0,i\in\ZZ}(-1)^{\dim G+i}\tr(x\vp',E')(A:H^i(K^x_D))\\&
=|\WW_0|\i\sum_{x\in\WW_0,i\in\ZZ}(-1)^{\dim G+i}\tr(x\vp,E)(A:H^i(K^x_D)).\endalign$$
(We have used that $\tr(x\vp,E)=\tr(x\vp',E')$ for $x\in\WW_0$, which follows from definitions.) Let 
$\WW_*=\prod_{j\in\ZZ/k\ZZ;j\ne0}\WW_j$. We note that the map $\WW_*\T\WW_0@>>>\WW$, $(y,x)\m yx\e(y)\i$ is a
bijection. Using 45.10(a) we see that the left hand side of (d) equals
$$\align&|\WW|\i\sum_{y\in\WW_*,x\in\WW_0,i\in\ZZ}(-1)^{\dim G+i}\tr(yx\e(y)\i\vp,E)(A:H^i(K^{yx\e(y\i)}_D))\\&
=|\WW|\i\sum_{y\in\WW_*,x\in\WW_0,i\in\ZZ}(-1)^{\dim G+i}\tr(x\vp,E)(A:H^i(K^x_D)).\endalign$$
Thus the two sides of (d) are equal.

Using (b) and the definitions we see that

(e) {\it if $D'$ has property $\tfA$ then $D$ has property $\tfA$.}
\nl
Note that if $O$ is a $G^0$-orbit on $D$ then $\b(O)$ is a $G'{}^0$-orbit on $D'$. Moreover, if $O'$ is a 
$G'{}^0$-orbit on $D'$ then $\b\i(O')$ is a $G^0$-orbit on $D$. We see that

(f) {\it the map $O\m\b(O)$ is a bijection between the set of $G^0$-orbits on $D$ and the set of $G'{}^0$-orbits 
on $D'$; the inverse bijection takes a $G'{}^0$-orbit $O'$ on $D'$ to $\b\i(O')$.}
\nl
We show:

(g) {\it if $D'$ has property $\fA_0$ then $D$ has property $\fA_0$}.
\nl
Let $A\in\hD^{unc}$. Then $\supp(A)$ is the closure of a single $G^0$-orbit $O$ in $D$. We have 
$A=\b^*A'[t]$ where $A'\in\hD'{}^{un}$. Hence $\supp(A')=\b\i(\supp(A))$. From (f) we see that 
$\supp(A')$ is the closure of a single $G'{}^0$-orbit $O'$ in $D'$. Hence $A'$ is cuspidal. By the assumption of 
(g) we see that $A'$ is zero outside $O'$. Hence $A$ is zero outside $\b\i(O')$ which is a single $G^0$-orbit 
necessarily equal to $O$. Thus $D$ has property $\fA_0$.

\head 46. Classification of unipotent character sheaves\endhead
\subhead 46.1\endsubhead
Let $p\ge1$ be the characteristic exponent of $\kk$. In this section we extend the results of \cite{\CS, IV,V} on
the classification of unipotent character sheaves on $D$ from the case $G=G^0$ to the general case.

In the remainder of this subsection we assume that $D=G^0$ 
and that (a) below holds:

(a) {\it if $G^0$  has a factor of type $E_8$ or $F_4$ then $p\ne2$.}
\nl
We note that:

(b)  {\it any character sheaf on $D$ is clean;}

(c)  {\it any admissible complex (see 6.7) on $D$ is a character sheaf.}
\nl
This is reduced to the case where $G^0$ is almost simple as in \cite{\CS, V, 23.21}. In that case, (b) is proved 
in \cite{\CS, IV,V} assuming in addition that: if $G^0$ has a factor $E_8$ then $p\ne3, p\ne5$; if $G^0$ has a 
factor $E_7$ or $F_4$ then $p\ne3$; if $G^0$ has a factor $E_6$ then $p\ne2$; if $G^0$ has a factor $G_2$ then 
$p\ne2,p\ne3$. In the remaining cases an additional argument (given by Shoji \cite{\SH, Sec.5} and Ostrik
\cite{\OS}) is needed. The fact that (b) implies (c) is proved as in \cite{\CS, IV,V}. 
 
\subhead 46.2\endsubhead
Assume that $G^0$ is semisimple and that for any proper parabolic subgroup $P$ of $G^0$ such that $N_DP\ne\em$ the
following condition is satisfied: any irreducible cuspidal admissible complex on $N_DP/U_P$ whose support contains
some unipotent element is a character sheaf. Let $A\in\hD^{unc}$ be such that for some unipotent $G^0$-orbit $S$ 
in $D$ and some irreducible cuspidal local system $\ce$ on $S$ we have $A=IC(\bS,\ce)[\dim S]$ extended by $0$ on
$D-\bS$. We assume that for any $G^0$-orbit $C\sub\bS-S$ there is no irreducible cuspidal local system on $C$. We
show:

(a) {\it $A$ is clean.}
\nl
The proof is along the lines of that of \cite{\CS, II, 7.9}. Assume that $A$ is not clean. Let $C\sub\bS-S$ be a 
$G^0$-orbit of minimum possible dimension such that $\ch^i(A)$ is nonzero on $C$ for some $i$; let $i_0$ be the
largest $i$ such that $\ch^i(A)$ is nonzero on $C$. Let $\cl$ be an irreducible local system on $C$ which is a 
direct summand of $\ch^{i_0}(A)|_C$. By our assumption, $\cl$ is not a cuspidal local system on $C$. By 8.8, 8.3,
8.2(b) we can find $(L',S')\in\AA$ (see 3.5) such that $S'$ contains unipotent elements and an irreducible 
cuspidal local system $\ce'$ on $S'$ such that, setting $\fK'=IC(\bY_{L',S'},\p_!\tce')$ extended by $0$ outside 
$\bY_{L',S'}$ ($\tce'$ as in 5.6), there exists a direct summand $A_1$ of $\fK'$ whose restriction to the 
unipotent variety of $D$ is (up to shift) $IC(\bC,\che\cl)$ extended to the unipotent variety by zero outside 
$\bC$. Let $(L'',S'')=(G^0,S)$. Our assumption implies that $L'\ne G^0$ so that $L',L''$ are not $G^0$-conjugate.
Hence 23.7 is applicable and yields $H^j_c(D,\fK'\ot A)=0$ for any $j$. Hence $H^j_c(D,A_1\ot A)=0$ for any $j$. 
Since $\supp(A)\sub\bS$ we have $\supp(A_1\ot A)\sub\bS$ so that $H^j_c(D,A_1\ot A)=H^j_c(\bS,A_1\ot A)$. Since 
$\supp(A_1)\cap\bS\sub\supp(\fK')\cap\bS\sub\bC$ we have $H^j_c(\bS,A_1\ot A)=H^j_c(\bC,A_1\ot A)$. Since $A$ is 
zero on $\bC-C$ (by the minimality of $C$) we have $H^j_c(\bC,A_1\ot A)=H^j_c(C,A_1\ot A)$. We see that 
$H^j_c(C,A_1\ot A)=0$ for all $j$. Since $A_1|_C$ is $\che\cl$ up to shift, it follows that 
$H^j_c(C,\che\cl\ot A)=0$ for all $j$. In particular we have $H^{2b+i_0}_c(C,\che\cl\ot A)=0$ where $b=\dim C$. 
Consider the spectral sequence $E^{r,s}_2=H^r_c(C,\ch^s(A)\ot\che\cl)\imp H^{r+s}_c(C,A\ot\che\cl)$. Then 
$E^{r,s}_2=0$ if $s>i_0$ (by our choice of $i_0$) or if $r>2b$. It follows that 
$E^{2b,i_0}_2=E^{2b,i_0}_3=\do=E^{2b,i_0}_\iy$. But $E^{2b,i_0}_\iy$ is a subquotient of 
$H^{2b+i_0}_c(C,A\ot\che\cl)$ hence it is zero. It follows that 
$0=E^{2b,i_0}_2=H^{2b}_c(C,\ch^{i_0}(A)\ot\che\cl)$. Since $\cl$ is a direct summand of $\ch^{i_0}(A)|_C$ it 
follows that $H^{2b}_c(C,\cl\ot\che\cl)=0$. This is a contradiction. This proves (a).

\subhead 46.3\endsubhead
In this subsection we assume that $G^0$ is almost simple, that $m:=|G/G^0|>1$, and that $\cz_G\sub G^0$.
Let $A\in\hD^{unc}$. Let $S$ be the stratum of $D$ such that $\supp(A)$ is
the closure of $S$. Now $A|_S$ is (up to shift) an irreducible cuspidal local system $\ce$. Note that $m$ is $2$ 
or $3$. Let $s\in G$ be a semisimple element and let $u\in G$ be a unipotent element such that $su=us\in S$. Let 
$G'=Z_G(s)$. Let $\d$ be the connected component of $G'$ that contains $u$. Let $S'$ be the (isolated) stratum of
$\d$ that contains $u$.  Let $\ce'$ be the inverse image of $\ce$ under $S'@>>>S$, $g\m sg$. Let 
$A'=IC(\bS',\ce')[\dim S']$ extended by $0$ on $\d-\bS'$. By 23.4(c), $A'$ is a direct sum of cuspidal admissible
complexes $A'_j$ on $G'{}^0$.

We show:

(a) {\it If $p\ne m$ then $A$ is clean. }
\nl
By our assumption, the image of 
$u$ in $G/G^0$ is $1$. Thus $u\in Z_{G^0}(s)$. Since $Z_{G^0}(s)/Z_{G^0}(s)^0$ has order prime to $p$ we see that
$u\in Z_{G^0}(s)^0$. Hence $\d=Z_{G^0}(s)^0=G'{}^0$. By 23.4(a) it is enough to show that each $A'_j$ is clean 
with respect to $G'$. This follows from 46.1(b),(c) applied to $G',G'{}^0$. (Note that $G'{}^0$ does not have a 
factor $E_8$; it can have a factor $F_4$ only if $G^0$ is of type $E_6$ and $p\ne2$, in which case 46.1(b),(c) 
are applicable.) This proves (a).

We show:

(b) {\it Assume that $G^0$ is of type $A_{n-1}$ ($n\ge3$) or $D_n$ ($n\ge2$). Assume that $p=m=2$ and that for any
proper parabolic subgroup $P$ of $G^0$ such that $N_DP\ne\em$ the following condition is satisfied: any 
irreducible cuspidal admissible complex on $N_DP/U_P$ is a character sheaf on $N_DP/U_P$. Then $A$ is clean.}
\nl
In this case the image of $s$ in $G/G^0$ is $1$. Hence $s\in G^0$ and $u\in D$. There is at most one cuspidal 
admissible complex on $D$. (See 12.9.) This complex must be isomorphic to $A$. Now the conclusion follows from 
46.2(a).

\subhead 46.4\endsubhead
In this subsection we assume that $G^0$ is simple of type $A_{n-1}$ ($n\ge3$), that $|G/G^0|=2$, that 
$\cz_G=\{1\}$ and that $D\ne G^0$. In this case $\e:\WW@>>>\WW$ is given by $w\m w_0ww_0\i$. In particular we have
$\Irr^\e(\WW)=\Irr(\WW)$ (see 43.1). We show:

(a) {\it $D$ has property $\fA$;}

(b) {\it $D$ has property $\tfA$;}

(c) {\it if $p=2$ then any irreducible cuspidal admissible complex on $D$ is in $\hD^{unc}$;}

(d) {\it for any $E_0\in\Irr(\WW)$ there is a unique object $A_{E_0}\in\hD^{un}$ (up to isomorphism) which 
satisfies $R_E=s_EA_{E_0}$ in $\ck^{un}_\QQ(D)$ for any $E\in\Irr(\tWW)$ such that $E|_{\QQ[W]}=E_0$ (here 
$s_E=\pm1$); moreover, $E_0\m A_{E_0}$ is a bijection from the set of isomorphism classes in $\Irr(\WW)$ to 
$\uhD^{un}$.}
\nl
We can assume that (a)-(d) hold when $n$ is replaced by $n'$ where $3\le n'<n$. (This assumption is empty if 
$n=3$.)

Note that if $P$ is a proper parabolic subgroup of $G^0$ such that $N_DP\ne\em$ and such that (setting 
$D'=N_DP/U_P$) either $\hD'{}^{unc}\ne\em$ or (if $p=2$) there is at least one cuspidal admissible complex on 
$D'$, then $P/U_P$ is of type $A_r$ (or a torus) and the induction hypothesis shows that $D'$ satisfies property 
$\fA_0$ and (if $p=2$) any irreducible cuspidal admissible complex on $D'$ is in $\hD'{}^{unc}$.

Using 46.3(a) (if $p\ne2$) and 46.3(b) (if $p=2$) we see that (a) holds.

Now let $E_0\in\Irr(\WW)$. We can extend $E_0$ to a $\tWW$-module $E$ in which $\vp$ acts as $w_0\in\WW$. 
We set $e_E=(-1)^{a_{E_0\ot\sgn}+l(w_0)}$, $e'_E=(-1)^{a_{E_0}}$. From 
\cite{\ORA, (7.6.6)} we see that there exists $x\in\boc_{E_0}$ such that $\ale_{x\vp}=e_E\ph_E$, 
$(-1)^{l(x)-\aa(x)}=e_Ee'_E$. Using 44.15(c) (which is applicable in view of (a)) we deduce that $e_ER_E$ is a 
$\ZZ$-linear combination of objects $A\in\hD^{un}$ such that $\ee^A=e_Ee'_E$. Since $(R_E:R_E)=1$ we 
deduce that $R_E=s_EA_{E_0}$  for a well defined $A_{E_0}\in\hD^{un}$ and $s_E=\pm1$; moreover,
$\ee^{A_{E_0}}=e_Ee'_E$. Since any $A\in\hD^{un}$ satisfies $(A:R_E)\ne0$ for some $E$ as above we see that 
$A=A_{E_0}$ for some $E_0$. Also if $E_0,E'_0$ are non-isomorphic objects of $\Irr(\WW)$ and $E,E'$ are the 
corresponding extension to $\tWW$ then $(R_E:R_{E'})=0$ hence $(A_{E_0}:A_{E'_0})=0$ so that 
$A_{E_0}\not\cong A_{E'_0}$. We see that (d) holds.

Let $E_0,E$ be as above. For $w\in\WW$, we have 

$(A_{E_0}:gr_1(K^w_D))=\pm(R_E:gr_1(K^w_D))=\pm\tr(w\vp,E)=\pm\tr(ww_0,E_0)$
\nl
(see 44.7(p)). Hence, by 44.14(a), 
the condition that $A_{E_0}$ is cuspidal is that $\tr(ww_0,E_0)=0$ whenever $w\in\WW$ is not $D$-anisotropic. Now 
$w\in\WW$ is not $D$-anisotropic if and only if $ww_0$ has even order. Thus the condition that $A_{E_0}$ is 
cuspidal is that $\tr(w',E_0)=0$ whenever $w'\in\WW$ has even order. The last condition holds if and only if $n$ 
is of the form $1+2+\do+s$ and $E_0$ corresponds to the partition of $n$ with parts $1,2,\do,s$. (See 
\cite{\CLA, 9.2, 9.3, 9.4}.) In this case we have $a_{E_0}=a_{E_0\ot\sgn}$ hence 
$\ee^{A_{E_0}}=(-1)^{l(w_0)}=(-1)^{\II_\e}=(-1)^{\codim(\supp(A_0))}$. (For the last equality see 44.8(a).) Thus 
the equality $\ee^A=(-1)^{\codim(\supp(A))}$ holds for any cuspidal $A\in\hD^{un}$. The analogous equality holds 
for non-cuspidal $A$ in view of the induction hypothesis and 44.15(a). We see that (b) holds.

Now assume that $p=2$. Let $\cx_1$ be the set of isomorphism classes of irreducible cuspidal admissible complexes
on $D$. Let $\cx_2$ be the set of isomorphism classes of objects in $\hD^{unc}$. Using 12.9 we see that 
$|\cx_1|=1$ if $n\in\{3,6,10,\do\}$ and $|\cx_1|=0$ otherwise. By the arguments above we see that $|\cx_2|=1$ if 
$n\in\{3,6,10,\do\}$. Clearly, $\cx_2\sub\cx_1$. It follows that $\cx_2=\cx_1$. This proves (c).

This completes the inductive proof of (a)-(d).

Let $E_0,E,x$ be as above. By 44.17(d) (which is applicable in view of (a),(b)) we have 
$(A_{E_0}:R_{\ale_{x\vp}})\in\NN$ hence $(A_{E_0}:e_ER_E)\in\NN$ hence $(s_ER_E:e_ER_E)\in\NN$ hence 
$s_Ee_E\in\NN$ hence $s_E=e_E$. Thus we have 

(e) $A_{E_0}=e_ER_E$.

\subhead 46.5\endsubhead
Assume that $G^0$ is semisimple and that $A$ is a cuspidal admissible sheaf on $D$ such that $\supp(A)$ is 
contained in the unipotent variety of $D$. Assume also that $G^0$ is of type $A_n\T A_n\T\do\T A_n$ ($r$ factors,
$n=1$ or $n=2$). We show:

(a) {\it $A$ is clean.}
\nl
By arguments in 12.3-12.6 we are reduced to the case where $G^0$ is almost simple and $\cz_G\sub G^0$. If $G=G^0$,
the conclusion follows from 46.1. Thus we can assume that $G\ne G^0$. As in 12.7 we see that we must have $n=2$, 
$p=2$, $|G/G^0|=2$. By 46.4(c), we have $A\in\hD^{unc}$; using this and 46.4(a), we see that $A$ is clean. This 
proves (a).

\subhead 46.6\endsubhead
In the setup of 46.3 we assume that $G^0$ is of type $D_4$ and $p=m=3$ or of type $E_6$ and $p=m=2$. Let
$A\in\hD^{unc}$. We show:

(a) {\it $A$ is clean.}
\nl
By 12.9 there is exactly one cuspidal admissible complex on $D$ (say $A'$) whose support is contained in the 
variety of unipotent elements in $D$. If $A\cong A'$ then $A$ is clean by 46.2(a). Hence we may assume that 
$\supp(A)$ is not contained in the variety of unipotent elements in $D$. In this case $G'{}^0$ is of type 
$A_1\T A_1\T A_1\T A_1$ (if $G$ is of type $D_4$) and of type $A_2\T A_2\T A_2$ (if $G$ is of type $E_6$).
By 23.4(a) it is enough to show that each $A'_j$ (as in 46.3) is clean with respect to $G'$. This follows from 
46.5(a) with $r=4,n=1$ or $r=3,n=2$. This proves (a).

\subhead 46.7\endsubhead
In this subsection we assume that $G^0$ is simple of type $D_4$, that $|G/G^0|=3$, that $\cz_G=\{1\}$, hence
$D\ne G^0$. We show:

(a) {\it $D$ has property $\fA$.}
\nl
Note that if $P$ is a proper parabolic subgroup of $G^0$ such that $N_DP\ne\em$ and such that (setting 
$D'=N_DP/U_P$) we have $\hD'{}^{unc}\ne\em$ then $P$ is a Borel subgroup so that $D'$ satisfies property $\fA_0$.
Using 46.3(a) (if $p\ne3$) and 46.6(a) (if $p=3$) we see that (a) holds.

The objects of $\Irr^\e(\WW)$ can be listed as: $1,4,1',4',2,6,8$ (each number represents an object of the 
corresponding degree; moreover, $1$ is the unit representation, $1'$ is the sign representation, $4$ is the 
reflection representation, $4'=4\ot 1'$). Each of these objects is naturally defined over $\QQ$ and it can be 
viewed as an object of $\Irr(\tWW)$ which is also defined over $\QQ$ with $\vp^3=1$ on it; we denote this object 
of $\Irr(\tWW)$ in the same way as the corresponding object in $\Irr^\e(\WW)$. From \cite{\ORA, (7.6.5)} we see 
that each of the elements 
$$\ph_{1},\ph_{4},\ph_{1'},\ph_{4'},\ph_{8}+\ph_{2},\ph_{8}-\ph_{2},\ph_{8}+\ph_{6},\ph_{8}-\ph_{6}$$
is of the form $\ale_{x\vp}$ for some $x\in\WW$ such that $l(x)-\aa(x)=0\mod2$. From this we deduce using 44.15(c)
that each of the elements 
$$R_{1},R_{4},R_{1'},R_{4'},R_{8}+R_{2},R_{8}-R_{2},R_{8}+R_{6},R_{8}-R_{6}\tag b$$ 
is a $\ZZ$-linear combination of objects $A\in\hD^{un}$ such that $\ee^A=1$. Since the elements (b) span
over $\QQ$ the same vector space as that spanned by the $R_E$ with $E\in\Irr(\tWW)$ and since each each 
$A\in\hD^{un}$ satisfies $(A:R_E)\ne0$ for some $E\in\Irr(\tWW)$ we see that each $A\in\hD^{un}$ has non-zero 
inner product with some element in (b) hence it satisfies $\ee^A=1$. If $A\in\hD^{unc}$ then 
$\codim(\supp(A))=|\II_\e|\mod2$; we have $|\II_\e|=2$ hence $\codim(\supp(A))=0\mod2$. Thus 
$\ee^A=(-1)^{\codim(\supp(A))}$ if $A\in\hD^{unc}$. The analogous equality holds for non-cuspidal $A$ in view of 
44.15(a) since it trivially holds on $D'$ as above. We see that

(c) {\it $D$ has property $\tfA$.}
\nl
By 44.17(d) (which is applicable in view of (a),(b)), the inner product of any $A\in\hD^{un}$ with any element in 
(b) is in $\NN$. Since the inner product of any two elements in (b) is known (it is $0,1$ or $2$) we see that 
there exist mutually nonisomorphic objects 
$$A_1,A_4,A_{1'},A_{4'},a,b,c,d\tag d$$
of $\hD^{un}$ such that
$$\align&R_{1}=A_1,R_{4}=A_4,R_{1'}=A_{1'},R_{4'}=A_{4'},R_{8}+R_{2}=a+b,\\&R_{8}-R_{2}=c+d,R_{8}+R_{6}=a+c,
R_{8}-R_{6}=b+d.\endalign$$
The list (d) exhausts the isomorphism classes in $\hD^{un}$ since any $A\in\hD^{un}$ has nonzero inner product 
with some element in (b). Note that $R_8=(a+b+c+d)/2$, $R_2=(a+b-c-d)/2$, $R_6=(a-b+c-d)/2$.

\subhead 46.8\endsubhead
In this subsection we assume that $G^0$ is simple of type $E_6$, that $|G/G^0|=2$, that $\cz_G=\{1\}$, hence
$D\ne G^0$. We show:

(a) {\it $D$ has property $\fA$.}
\nl
Note that if $P$ is a proper parabolic subgroup of $G^0$ such that $N_DP\ne\em$ and such that (setting 
$D'=N_DP/U_P$) there is at least one cuspidal admissible complex on $D'$ then $P/U_P$ is either of type $A_5$ or a
torus. (The case where $P/U_P$ is of type $D_4$ is excluded using 23.4(a) when $p\ne2$ and 12.9(b) when $p=2$.) In
either case $D'$ satisfies property $\fA_0$. Using 46.3(a) (if $p\ne2$) and 46.6(a) (if $p=2$) we see that (a) 
holds.

In our case $\e:\WW@>>>\WW$ is given by $w\m w_0ww_0\i$. The objects of $\Irr(\WW)$ (up to isomorphism) can be 
listed as 
$$\align&1_0,6_1,20_2,30_3,15_3,\ti{15}_3,64_4,60_5,81_6,24_6,80_7,60_7,90_7,10_7,\\&
20_7,81_{10},60_{11},24_{12},64_{13},
30_{15},15_{15},\ti{15}_{15},20_{20},6_{25},1_{36}\endalign$$
where $N_n$ or $\ti N_n$ denotes an object $E_0\in\Irr(\WW)$ such that $\dim E_0=N,a_{E_0}=n$. Each object of 
$\Irr(\WW)$ can be regarded as an object of $\Irr(\tWW)$ on which $\vp$ acts as $w_0$; this object of $\Irr(\tWW)$
is denoted in the same way as the corresponding object in $\Irr(\WW)$. From \cite{\ORA, 7.10} we see that each of
the elements 

$\ph_{1_0},-\ph_{6_1},\ph_{20_2},-\ph_{60_5},\ph_{24_6},\ph_{81_6},\ph_{81_{10}},\ph_{24_{12}},
-\ph_{60_{11}},\ph_{20_{20}},-\ph_{6_{25}},\ph_{1_{36}}$,

$-\ph_{30_3}-\ph_{15_3},-\ph_{30_3}+\ph_{15_3},-\ph_{30_3}-\ph_{\ti{15}_3},-\ph_{30_3}+\ph_{\ti{15}_3}$,

$-\ph_{30_{15}}-\ph_{15_{15}},-\ph_{30_{15}}+\ph_{15_{15}},-\ph_{30_{15}}-\ph_{\ti{15}_{15}},
-\ph_{30_{15}}+\ph_{\ti{15}_{15}}$,

$-\ph_{80_7}+\ph_{60_7}+\ph_{10_7},-\ph_{80_7}-\ph_{60_7}+\ph_{10_7}, -2\ph_{80_7}-\ph_{10_7}$,

$-\ph_{80_7}+\ph_{60_7}+\ph_{90_7},-\ph_{80_7}-\ph_{60_7}+\ph_{90_7}, -2\ph_{80_7}-\ph_{90_7}$, 
$-\ph_{80_7}-\ph_{20_7}$
\nl
is of the form $\ale_{x\vp}$ ($x\in\WW,l(x)=\aa(x)\mod2$) and that each of the elements

$-\ph_{64_4},\ph_{64_{13}}$
\nl
is of the form $\ale_{x\vp}$ ($x\in\WW,l(x)\ne\aa(x)\mod2$). From this we deduce using 44.15(c) that each of the 
elements 

(b) $R_{1_0},-R_{6_1},R_{20_2},-R_{60_5},R_{24_6},R_{81_6},R_{81_{10}},R_{24_{12}},
-R_{60_{11}},R_{20_{20}},-R_{6_{25}},R_{1_{36}}$,

$-R_{30_3}-R_{15_3},-R_{30_3}+R_{15_3},-R_{30_3}-R_{\ti{15}_3},-R_{30_3}+R_{\ti{15}_3}$,

$-R_{30_{15}}-R_{15_{15}},-R_{30_{15}}+R_{15_{15}},-R_{30_{15}}-R_{\ti{15}_{15}},-R_{30_{15}}+R_{\ti{15}_{15}}$,

$-R_{80_7}+R_{60_7}+R_{10_7},-R_{80_7}-R_{60_7}+R_{10_7},-2R_{80_7}-R_{10_7}$,

$-R_{80_7}+R_{60_7}+R_{90_7},-R_{80_7}-R_{60_7}+R_{90_7},-2R_{80_7}-R_{90_7},-R_{80_7}-R_{20_7}$
\nl
is a $\ZZ$-linear combination of objects $A\in\hD^{un}$ such that $\ee^A=1$ and that each of the elements 

(c) $-R_{64_4},R_{64_{13}}$
\nl
is a $\ZZ$-linear combination of objects $A\in\hD^{un}$ such that $\ee^A=-1$. Since the elements in (c)
have self-inner product $1$, we have $R_{64_4}=\pm A$, $R_{64_{13}}=\pm A'$ where $A,A'\in\hD^{un}$. Since 
$(R_{64_4}:R_{64_{13}})=0$ we see that $A\not\cong A'$. By 44.8(c) we have $\dd(R_{64_4})=R_{64_{13}}$ hence 
$\dd(A)=\pm A'$. If $A$ were cuspidal we would have $\dd(A)=A$. Thus $A$ is not cuspidal. Similarly $A'$ is not 
cuspidal. If $A_1\in\hD^{unc}$ then $A_1$ must have non-zero inner product with some $R_E$ hence with at least 
one of the elements in (b),(c). But we have just seen that its inner product with any element in (c) is zero. 
Thus, $A_1$ must have non-zero inner product with at least one of the elements in (b). It follows that
$\ee^{A_1}=1$. We have $\codim(\supp(A_1))=|\II_\e|\mod2$; moreover $|\II_\e|=4$ hence 
$\codim(\supp(A_1))=0\mod2$. Thus, $\ee^A=(-1)^{\codim(\supp(A))}$ if $A\in\hD^{unc}$. The analogous equality 
holds for non-cuspidal $A$ in view of 44.15(a) since it holds on $D'$ as above, by 46.4(b). We see that:

(d) {\it $D$ has property $\tfA$.}
\nl
By 44.17(d) (which is applicable in view of (a),(d)), the inner product of any $A\in\hD^{un}$ with any element in
(b) or (c) is in $\NN$. Since the inner products of any two elements in (b) or (c) are known we see that there 
exist mutually nonisomorphic objects 

$A_{1_0},A_{6_1},A_{20_2},A_{60_5},A_{24_6},A_{81_6},A_{81_{10}},A_{24_{12}},
A_{60_{11}},A_{20_{20}},A_{6_{25}},A_{1_{36}}$,

$a_3,b_3,c_3,d_3,a_{15},b_{15},c_{15},d_{15}$, $a,b,c,d,e,f,g,h$
\nl
of $\hD^{un}$ such that

$R_{1_0}=A_{1_0}$, $-R_{6_1}=A_{6_1}$, $R_{20_2}=A_{20_2}$, $-R_{60_5}=A_{60_5}$, $R_{24_6}=A_{24_6}$, 

$R_{81_6}=A_{81_6}$, $R_{81_{10}}=A_{81_{10}}$, $R_{24_{12}}=A_{24_{12}},
-R_{30_3}-R_{15_3}=a_3+b_3$, 

$-R_{30_3}+R_{15_3}=c_3+d_3$, $-R_{30_3}-R_{\ti{15}_3}=a_3+c_3$, 
$-R_{30_3}+R_{\ti{15}_3}=b_3+d_3$,

$-R_{30_{15}}-R_{15_{15}}=a_{15}+b_{15}$, $-R_{30_{15}}+R_{15_{15}}=c_{15}+d_{15}$, 
$-R_{30_{15}}-R_{\ti{15}_{15}}=a_{15}+c_{15}$, 

$-R_{30_{15}}+R_{\ti{15}_{15}}=b_{15}+d_{15}$,

$-R_{80_7}+R_{60_7}+R_{10_7}=a+b+d$, $-R_{80_7}-R_{60_7}+R_{10_7}=d+e+f$, 

$-2R_{80_7}-R_{10_7}=b+c+f+g+h$, $-R_{80_7}+R_{60_7}+R_{90_7}=a+b+c$, 

$-R_{80_7}-R_{60_7}+R_{90_7}=c+e+f$, $,-2R_{80_7}-R_{90_7}=b+d+f+g+h$,

$-R_{80_7}-R_{20_7}=b+f$.
\nl
(We use \cite{\ORA, 7.7(iii)}.) Hence we have 

$-R_{80_7}=(a+3b+2c+2d+e+3f+2g+2h)/6$, $R_{60_7}=(a+b-e-f)/2$, 

$R_{90_7}=(a+2c-d+e-g-h)/3$, $R_{10_7}=(a-c+2d+e-g-h)/3$, 

$-R_{20_7}=(a-3b+2c+2d+e-3f+2g+2h)/6$.

\subhead 46.9\endsubhead
We fix an integer $n\ge1$. Let $W_n$ be the group of all permutations of $\{1,2,\do,n,n',\do,2',1'\}$ which 
commute with the involution $i\lra i'$. For each $j\in[1,n-1]$ let $s_j\in W_n$ be the involution which 
interchanges $j,j+1$ and also $j',(j+1)'$ and leaves the other elements unchanged. Let $s_n\in W_n$ be the 
permutation which interchanges $n,n'$ and leaves the other elements unchanged. Define a homomorphism 
$\c:W_n@>>>\{\pm1\}$ by the condition $\c(s_j)=1$ if $j\in[1,n-1]$, $\c(s_n)=-1$. 

We now assume that $n\ge2$. Then $W'_n:=\ker\c$ is a Coxeter group on the generators $s_j(j\in[1,n-1])$ and 
$s_ns_{n-1}s_n$.

For $h\in[2,n-1]$ let $W_{n,h}$ be the subgroup of $W_n$ consisting of the permutations in $W_n$ which carry each
of 

$\{1,2,\do,n-h\}$, $\{n-h+1,n-h+2,\do,n,n',\do,(n-h+2)',(n-h+1)'\}$, 

$\{1',2',\do,(n-h)'\}$
\nl
into itself. We may identify in an obvious way $W_{n,h}$ with $\fS_{n-h}\T W_h$ where $\fS_{n-h}$ is the 
symmetric group in $n-h$ letters. 

\subhead 46.10\endsubhead
Let $m\in\NN$. Let $X_n^m$ be the set of all ordered pairs $(S,T)$ ("symbols") of distinct subsets of $\NN$ (with
$|S|=|T|=m$) such that 

$\sum_{x\in S}x+\sum_{x\in T}x=n+m^2-m$.
\nl
We define a "shift" map $X_n^m@>>>X_n^{m+1}$ by 
$(S,T)\m(\{0\}\cup(S+1),\{0\}\cup(T+1))$. Using the shift maps we can form the direct limit 
$X_n=\lim_{m\to\iy}X_n^m$. We have an obvious map $X_n^m@>>>X_n$. If $m\ge n$ then any $(S,T)\in X_n^{m+1}$ 
satisfies $0\in S,0\in T$. Hence if $m\ge n$, the shift map $X_n^m@>>>X_n^{m+1}$ is a bijection. We shall 
sometimes identify $X_n$ with $X_n^m$ with some fixed $m\ge n$.
But some elements of $X_n$ can be represented by elements of $X_n^m$ where $m<n$. 

Note that if $(S,T)\in X_n^m$ then $S\cup T\sub[0,n+m-1]$. Thus $X_n^m$ is finite for any $m$ so that $X_n$ is 
finite. 

Let $\bX_n^m$ be the set of all pairs $(M,N)$ of disjoint subsets of $\NN$ such that $M\ne\em$, $|M|+2|N|=2m$ and

$\sum_{x\in M}x+2\sum_{x\in N}x=n+m^2-m$.
\nl
We define a "shift" map $\bX_n^m@>>>\bX_n^{m+1}$ by 
$(M,N)\m(M+1,\{0\}\cup(N+1))$. Using the shift maps we can form the direct limit $\bX_n=\lim_{m\to\iy}\bX_n^m$. We
have an obvious map $\bX_n^m@>>>\bX_n$. If $m\ge n$, then any $(M,N)\in\bX_n^{m+1}$ satisfies $0\in N$ (hence 
$0\n M$). Hence if $m\ge n$, the shift map $\bX_n^m@>>>\bX_n^{m+1}$ is a bijection. We shall sometimes identify 
$\bX_n$ with $\bX_n^m$ with some fixed $m\ge n$. 

For $(M,N)\in\bX^m_n$ let $\cv_M$ (resp. $V_M$) be the set of all subsets of $M$ with cardinal $|M|/2$ (resp. with
even cardinal); we regard $V_M$ as an $\FF_2$-vector space with addition $E,E'\m E*E'=(E\cup E')-(E\cap E')$. Let

$V'_M=\{\et:V_M@>>>\FF_2\text{-linear, }\et(M)=1\}$; 
\nl
here $M$ is viewed as an element of $V_M$. 

Define $t_M:M@>>>\FF_2$ by $t_M(x)=|\{x'\in M;x'<x\}|\mod2$. Define an injective map $\cv_M@>>>V_M$ by 

$H\m H^\sh:=t_M\i(1)*H$;
\nl
the image of this map is denoted by $\tcv_M$.

We define a (surjective) map $\z:X_n^m@>>>\bX_n^m$ by $(S,T)\m(S*T,S\cap T)$; if $(M,N)\in\bX_n^m$, then 
$H\m(N\cup H,N\cup(M-H))$ is a bijection $\cv_M\lra\z\i(M,N)$.

\subhead 46.11\endsubhead
An irreducible $\QQ[W_n]$-module is said to be {\it nondegenerate} if its restriction to $W'_n$ is irreducible. To
a nondegenerate irreducible $\QQ[W_n]$-module we associate an element $(S,T)$ of $X_n$ as in \cite{\CLA, 2.7(ii)}.
We obtain a bijection $[[S,T]]\lra(S,T)$ between the set of nondegenerate irreducible $\QQ[W_n]$-modules (up to 
isomorphism) and $X_n$. Note that $[[S,T]]$ and $[[T,S]]$ have the same restriction to $W'_n$.

\subhead 46.12\endsubhead
In 46.12-46.24 we assume that $G^0$ is adjoint of type $D_n$ ($n\ge2$), that $|G/G^0|=2$, that $\cz_G=\{1\}$ hence
$D\ne G^0$. We choose an isomorphism of $\WW$ with $W'_n$ as Coxeter groups and we use it to identify the two
groups. We define a surjective homomorphism $\tWW@>>>W_n$: it takes $\vp$ to $s_n$ and its restriction to $\WW$ is
the obvious imbedding $\WW=W'_n@>>>W_n$. Via this homomorphism any nondegenerate irreducible $\QQ[W_n]$-module can
be viewed as an object of $\Irr(\tWW)$ so that the set of isomorphism classes of objects of $\Irr(\tWW)$ can be 
identified with the set of isomorphism classes of nondegenerate irreducible $\QQ[W_n]$-modules, hence with the set
$\{[[S,T]];(S,T)\in X_n\}$. Note that for $(S,T),(S',T')$ in $X_n$ we have $\z(S,T)=\z(S',T')$ if and only if the
two sided cells attached to $[[S,T]]$ and to $[[S',T']]$ coincide. Thus $\bX_n$ may be viewed as as indexing set 
for the two-sided cells of $\WW$ which are $\e$-stable. We write $\boc_{M,N}$ for the two-sided cell of $\WW$ 
corresponding to $(M,N)\in\bX_n$.

\subhead 46.13\endsubhead
For any two-element subset $C$ of $\NN$ let $[C]$ be the closed interval in $\RR$ with extremities in $C$. Let $M$
be a finite non-empty subset of $\NN$ of even cardinal. An {\it admissible arrangement} of $M$ is a set $\Ph$ of 
two-element subsets of $M$ forming a partition of $M$ with the following property: for any four element subset of
$M$ of the form $C\sqc C'$ where $C\in\Ph$, $C'\in\Ph$, we have $[C]\sub[C']$ or $[C']\sub[C]$ or $[C]\cap[C']=0$.
(This agrees with the definition in \cite{\ORA, p.164}.) For example the admissible arrangements of 
$\{0,1,2,3,4,5\}$ are 

$\Ph_1=\{(0,1),(2,3),(4,5)\}$, $\Ph_2=\{(0,5),(1,2),(3,4)\}$, $\Ph_3=\{(0,3),(1,2),(4,5)\}$,
$\Ph_4=\{(0,1),(2,5),(3,4)\}$, $\Ph_5=\{(0,5),(1,4),(2,3)\}$. 
\nl
If $\Ps$ is a subset of $\Ph$ and $i\in\FF_2$ we denote by $\Ps^i$ the set of all $x\in t_M\i(i)$ such that $x$
belongs to some pair in $\Ps$.

Now let $(M,N)\in\bX_n^m$. Let $\Ph$ be an admissible arrangement of $M$ and let $\hPh\sub\Ph$ be a subset such 
that $|\hPh|$ is odd. We set
$$c(M,N,\Ph,\hPh)=\fra{1}{2}\sum_{\Ps\sub\Ph}(-1)^{|\hPh\cap\Ps|}
\ph_{[[\Ps^0\cup(\Ph-\Ps)^1\cup N,\Ps^1\cup(\Ph-\Ps)^0\cup N]]}\in\car(\tWW).$$
The last inclusion holds since for any $\Ps\sub\Ph$ we have $(-1)^{|\hPh\cap\Ps|}=-(-1)^{|\hPh\cap(\Ph-\Ps)|}$.
From \cite{\ORA, (5.18.1)} we see that

(a) {\it there exists $x\in\WW$ such that $c(M,N,\Ph,\hPh)=\ale_{x\vp}$ and $l(x)=\aa(x)\mod2$.}
\nl
From \cite{\ORT, 1.19} we see that 

(b) {\it if $H\in\cv_M$ then there exists an admissible arrangement $\Ph$ of $M$ and $\Ps\sub\Ph$ such that 
$H=\Ps^0\cup(\Ph-\Ps)^1$ that is,

$[[N\cup H,N\cup(M-H)]]=[[\Ps^0\cup(\Ph-\Ps)^1\cup N,\Ps^1\cup(\Ph-\Ps)^0\cup N]]$;
\nl
 moreover,}
$$\ph_{[[N\cup H,N\cup(M-H)]]}=
2^{-|M/2|+1}\sum_{\hPh\sub\Ph;|\hPh|=\text{odd}}(-1)^{|\hPh\cap\Ps'|}c(M,N,\Ph,\hPh).$$

\subhead 46.14\endsubhead
We now state some properties (a)-(d) of $D$.

(a) {\it $D$ has property $\fA$;}

(b) {\it $D$ has property $\tfA$.}
\nl
In view of (a),(b), the results in 44.17-44.21 are applicable to $D$. In particular for any $\e$-stable two-sided 
cell $\boc$ of $\WW$, the subcategory $\hD^{un}_{\boc}$ of $\hD^{un}$ is defined as in 44.19. We shall write 
$\hD^{un}_{M,N}$, $\uhD^{un}_{M,N}$ instead of $\hD^{un}_{\boc_{M,N}}$, $\uhD^{un}_{\boc_{M,N}}$ where 
$(M,N)\in\bX_n$. 

(c) {\it For any $m\ge n$ and any $(M,N)\in\bX_n^m$ there exists a bijection $\et\m A_\et$, 
$V'_M\lra\uhD^{un}_{M,N}$ such that 
$$(A_\et:R_{[[N\sqc H,N\sqc(M-H)]]})=2^{-|M|/2+1}(-1)^{\et(t_M\i(1)*H)}$$
for any $\et\in V'_M$, $H\in\cv_M$;}

(d) {\it if $p=2$ then any irreducible cuspidal admissible complex on $D$ is in $\hD^{unc}$; moreover, $\hD^{unc}$
is empty unless $n=s^2$ with $s$ odd, $s\ge3$, in which case $\hD^{unc}$ has exactly one object up to isomorphism;
its support is contained in the set of unipotent elements of $D$.}
\nl
The proofs for (a)-(d) are given in 46.15-46.23 under the induction hypothesis that (a)-(d) hold when $n$ is 
replaced by $n'$ with $2\le n'<n$. (This assumption is empty if $n=2$.)

\subhead 46.15\endsubhead
If $P$ is a proper parabolic subgroup of $G^0$ such that $N_DP\ne\em$ and such that (setting $D'=N_DP/U_P$) either
$\hD'{}^{unc}\ne\em$ or (if $p=2$) there is at least one cuspidal admissible complex on $D'$, then $P/U_P$ is of
type $D_r$ (or a torus) and the induction hypothesis shows that $D'$ satisfies property $\fA_0$ and (if $p=2$) any
irreducible cuspidal admissible complex on $D'$ is in $\hD'{}^{unc}$. 

Using 46.3(a) (if $p\ne2$) and 46.3(b) (if $p=2$) we see that 46.14(a) holds.

Using 46.13(a) and 44.15(c) (which is applicable in view of 46.14(a)) we see that for any $M,N,\Ph,\hPh$ as in 
46.13(a), $R_{c(M,N,\Ph,\hPh)}$ is a $\ZZ$-linear combination of objects $A\in\hD^{un}$ such that $\ee^A=1$. Using
46.13(b) we deduce that for any $E\in\Irr(\tWW)$, $R_E$ is a $\ZZ$-linear combination of objects $A\in\hD^{un}$ 
such that $\ee^A=1$. Since any $A\in\hD^{un}$ appears with non-zero coefficient in $R_E$ for some 
$E\in\Irr(\tWW)$, we see that any $A\in\hD^{un}$ satisfies $\ee^A=1$.

We show: 

(a) {\it if $\hD^{unc}\ne\em$ then $n$ is odd.}
\nl
If $p=2$ this follows from 12.9(b). If $p\ne2$ then we can find an isolated semisimple element $s\in D$ such that
$Z_G(s)^0$ carries a cuspidal admissible complex supported on the unipotent variety of $Z_G(s)^0$ (see 23.4(b)).
Now $Z_G(s)^0$ is either semisimple of type $B_{n-1}$ (and then $n-1$ must be even by the known theory for
connected classical groups) or is semisimple of type $B_a\T B_b$ with $a\ge1,b\ge1,a+b=n-1$ (and then $a,b$ must 
be even and $n-1$ must be even). Thus (a) holds.

Now if $A\in\hD^{unc}$, we have $(-1)^{\codim(\supp(A))}=(-1)^{|\II_\e|}=(-1)^{n-1}$ and this equals $1$ by (a).
Thus we have $\ee^A=(-1)^{\codim(\supp(A))}$ for any cuspidal $A\in\hD^{un}$. The analogous equality holds for
non-cuspidal $A$ in view of the induction hypothesis and 44.15(a). We see that 46.14(b) holds.

\subhead 46.16\endsubhead
For $h\in[2,n-1]$ let $P^h$ be the parabolic subgroup of $G^0$ which contains $B^*$ and is such that the Weyl
group of $P^h/U_{P^h}$ is the subgroup of $\WW_{I^h}:=W'_n\cap W_{n,h}$ of $\WW=W'_n$. Then $\tWW_{I^h}$ (the 
subgroup of $\tWW$ generated by $\WW_{I^h}$ and $\vp$, see 43.8) is the inverse image under $\tWW@>>>W_n$ of 
$W_{n,h}$ and $\Irr(\tWW_{I^h})$ can be identified under $\tWW_{I^h}@>>>W_{n,h}$ with the set of isomorphism 
classes of irreducible $\QQ[W_{n,h}]$-modules of the form $E\bxt E'$ where $E$ is an irreducible 
$\QQ[\fS_{n-h}]$-module and $E'$ is an irreducible nondegenerate $\QQ[W_h]$-module. Let $G^h=N_GP^h/U_{P^h}$. Then
$D^h=N_DP^h/U_{P^h)}$ is a connected component of $G^h$. We have $G^h/Z^h=PGL_{n-h}\T\bG^h$ where $Z^h$ is a one 
dimensional torus in the centre of $(G^h)^0$ and $\bG^h$ is a group like $G$ (with $n$ replaced by $h$). Hence 
46.14(a)-46.14(d) hold for $D^h$ instead of $D$ (by the induction hypothesis) and the objects in $(\hD^h)^{un}$ 
can be written in the form $A\bxt A'$ with $A\in\widehat{PGL}_{n-h}^{un}$ and $A'\in(\hat{\bD}^h)^{un}$ (where 
$\bD^h=D^h/Z^h$).

\subhead 46.17\endsubhead
Using 46.13(a) and 44.17(d) (which is applicable in view of 46.14(a), 46.14(b)) we see that for any 
$(M,N)\in\bX_n^m$, any admissible arrangement $\Ph$ of $M$ and any $\hPh\sub\Ph$ with $|\hPh|=\text{odd}$ we have
that $R_{c(M,N,\Ph,\hPh)}$ is a $\NN$-linear combination of objects in $\hD^{un}$ or equivalently that
$$\fra{1}{2}\sum_{\Ps\sub\Ph}(-1)^{|\hPh\cap\Ps|}R_{[[\Ps^0\cup(\Ph-\Ps)^1\cup N,\Ps^1\cup(\Ph-\Ps)^0\cup N]]}$$
is an $\NN$-linear combination of objects in $\hD^{un}$. 

\subhead 46.18\endsubhead
We prove 46.14(c) assuming that $|M|=2$. We have $M=\{x,y\}$ with $x<y$. From 46.17 we see that
$R_{[[N\sqc\{y\},N\sqc(\{x\})]]}$ is an $\NN$-linear combination of objects in $\hD^{un}_{M,N}$. Since 
$R_{[[N\sqc\{y\},N\sqc(\{x\})]]}$ has self inner product $1$ it must be equal to a single object of 
$\hD^{un}_{M,N}$ and the desired result follows.

\subhead 46.19\endsubhead
We prove 46.14(c) assuming that $|M|=4$. We have $M=\{x,y,z,u\}$ with $x<y<z<u$. From 46.17 we see that

(a) $R_{[[N\sqc\{y,u\},N\sqc(\{x,z\})]]}\pm R_{[[N\sqc\{x,u\},N\sqc(\{y,z\})]]}$,

 $R_{[[N\sqc\{y,u\},N\sqc(\{x,z\})]]}\pm R_{[[N\sqc\{z,u\},N\sqc(\{x,y\})]]}$
\nl
are $\NN$-linear combinations of objects in $\hD^{un}_{M,N}$. Since the inner products of any two elements in (a) 
are known (they are $0,1$ or $2$) we see that there exist four mutually non-isomorphic objects $a,b,c,d$ of 
$\hD^{un}_{M,N}$ such that

$R_{[[N\sqc\{y,u\},N\sqc(\{x,z\})]]}+R_{[[N\sqc\{x,u\},N\sqc(\{y,z\})]]}=a+b$,

$R_{[[N\sqc\{y,u\},N\sqc(\{x,z\})]]}-R_{[[N\sqc\{x,u\},N\sqc(\{y,z\})]]}=c+d$,

$R_{[[N\sqc\{y,u\},N\sqc(\{x,z\})]]}+R_{[[N\sqc\{z,u\},N\sqc(\{x,y\})]]}=a+c$,

$R_{[[N\sqc\{y,u\},N\sqc(\{x,z\})]]}-R_{[[N\sqc\{z,u\},N\sqc(\{x,y\})]]}=b+d$.
\nl
Hence we have 

$R_{[[N\sqc\{y,u\},N\sqc(\{x,z\})]]}=(a+b+c+d)/2$,

$R_{[[N\sqc\{x,u\},N\sqc(\{y,z\})]]}=(a+b-c-d)/2$,

$R_{[[N\sqc\{z,u\},N\sqc(\{x,y\})]]}=(a-b+c-d)/2$.
\nl
There are well defined elements $\et_a,\et_b,\et_c,\et_d$ of $V'_M$ such that
$$\align&\et_a(\{x,y\})=0,\et_a(\{y,z\})=0, \et_b(\{x,y\})=0,\et_b(\{y,z\})=1,\\&
\et_c(\{x,y\})=1,\et_c(\{y,z\})=0,\et_d(\{x,y\})=1,\et_d(\{y,z\})=1.\endalign$$
The assignment $\et_a\m a$, $\et_b\m b$, $\et_c\m c$, $\et_d\m d$ is a bijection $V'_M\lra\uhD^{un}_{M,N}$ which 
establishes 46.14(c) in our case. 

\subhead 46.20\endsubhead
We now assume that $|M|\ge4$ and that $(M,N)$ has the following property: there exists $k\in[0,\max(M\cup N)]$ 
such that $k\n M\cup N$. We set 

$h=n-|\{x>k;x\in M\}|-2|\{x>k;x\in N\}|$.
\nl
Clearly, $h<n$. Let 

$M'=\{x<k;x\in M\}\sqc\{x\ge k;x+1\in M\}$, 

$N'=\{x<k;x\in N\}\sqc\{x\ge k;x+1\in N\}$.
\nl
Note that $M',N'$ are disjoint subsets of $\NN$ such that $|M'|=|M|,|N'|=|N|$ and 
$$\sum_{x\in M'}x+2\sum_{x\in N'}x=\sum_{x\in M}x+2\sum_{x\in N}x-(n-h)=h+m^2-m.$$
In particular, $h\ge0$. If $h\le1$ we see that $|M'|=2h$ hence $|M|=2h<4$, a contradiction. Thus we have 
$h\in[2,n-1]$. We see also that $(M',N')\in X_h$. We define a bijection $M'@>\si>>M$ by $x\m x$ if $x<k$ and 
$x\m x+1$ if $x\ge k$. This induces a bijection $V_{M'}@>\si>>V_M$ hence a bijection $V'_M@>\si>>V'_{M'}$.
Consider the two-sided cell $\boc'=\boc_{M',N'}\T\boc_0$ of $\WW_{I^h}$ (see 46.16) where $\boc_0$ is the 
two-sided cell associated to the sign representation $\sgn_h$ of $\fS_{n-h}$. We have $\boc'\sub\boc$ where 
$\boc=\boc_{M,N}$. Moreover, $\boc',\boc$ satisfy the assumptions (i),(ii) of 44.21. Consider the composite 
bijection
$$V'_M@>\si>>V'_{M'}@>\si>>(\un{\hat{\bD}}^h)^{un}_{M',N'}@>\si>>(\uhD^h)^{un}_{\boc'}@>\si>>\uhD^{un}_{M,N};$$
here the first bijection is as above; the second bijection comes from the induction hypothesis; the third 
bijection is $A'\m A\bxt A'$ where $A=R_{\sgn_{n-h}}\in\widehat{PGL}_{n-h}^{un}$; the fourth bijection comes from
44.21(h). Using 44.21(h) we see that this composite bijection has the required properties. This proves 46.14(c) in
our case.

\subhead 46.21\endsubhead
We now assume that $|M|\ge4$ and that there exists $y>0$ such that $y\in N,y-1\n N$. Recall that 
$M\cup N\sub[0,m+n-1]$. We can assume that $m=n$ so that $M\cup N\sub[0,t]$ where $t=2n-1$. Let 

$M'=\{x;t-x\in M\}\sub\NN$, $N'=\{x\in[0,t];t-x\n M\cup N\}\sub\NN$.
\nl
We have $M'\cap N'=\em$, $|M'|+2|N'|=|M|+2(t+1)-2|M\cup N|=2n$,
$$\align&\sum_{x;t-x\in M}x+2\sum_{x\in[0,t];t-x\n M\cup N}x
=\sum_{x\in M}(t-x)+2\sum_{x\in[0,t]}x-2\sum_{x\in M\cup N}(t-x)\\&
=|M|t-\sum_{x\in M}x+t^2+t-2|M|t-2|N|t+2\sum_{x\in M}x+2\sum_{x\in N}x\\&
=t^2+t-|M|t-2|N|t+\sum_{x\in M}x+2\sum_{x\in N}x=n^2.\endalign$$
We see that $(M',N')\in X_n^n$. We have a bijection $M'@>\si>>M$, $x\m t-x$. This induces a bijection 
$V_{M'}@>\si>>V_M$ and a bijection $V'_M@>\si>>V'_{M'}$. Since $y\in N$, we have $y\n M$ hence $t-y\n M'$. Since 
$y\in N$, we have $t-y\n N'$. Thus, $t-y\n M'\cup N'$. If $y-1\in M$, then $t-y+1\in M'$ and $t-y<t-y+1$. If 
$y-1\n M$, then $y-1\n M\cup N$ (since $y-1\n N$) hence $t-y+1\in N'$ and $t-y<t-y+1$. In any case we have 
$t-y+1\in M'\cup N'$ and $t-y\in[0,\max(M'\cup N')]$. By 46.20, 46.14(c) holds when $(M,N)$ is replaced by 
$(M',N')$. Consider the composite bijection 
$$V'_M@>\si>>V'_{M'}@>\si>>\uhD^{un}_{M',N'}@>\si>>\uhD^{un}_{M,N};$$
here the first bijection is as above; the second bijection is as in 46.14(c) for $(M',N')$; the third bijection is
$A\m A^\circ$, see 44.19(a). (Note that for $A\in\hD^{un}$ we have $A^\circ=\dd(A)$ since $\ee^A=1$ by 46.15.) The
composite bijection above is denoted by $\et\m A_\et$. We have $A_\et=\dd(A_{\et'})$ where $\et\in V'_M$ 
corresponds to $\et'\in V'_{M'}$ and $A_{\et'}$ is attached to $\et'$ by 46.14(c) for $(M',N')$. For any 
$J\sub M$, let $J'\sub M'$ be the image of $J$ under $x\m t-x$. Let $H\in\cv_M$. Using 44.8(c) and 
\cite{\ORT, (1.4.1)} we have 
$$(A_\et:R_{[[N\sqc H,N\sqc(M-H)]]})=(\dd(A_{\et'}):\dd(R_{[[N'\sqc(M'-H'),N'\sqc H']])}).$$
(We have $[[N\sqc H,N\sqc(M-H)]]\ot\sgn=[[N'\sqc(M'-H'),N'\sqc H']]$.) This equals 
$$(A_{\et'}:R_{[[N'\sqc(M'-H'),N'\sqc H']]})=2^{-|M'|/2+1}(-1)^{\et'((M'-H')*t_{M'}\i(1))}.$$
(We have used 46.14(c) for $(M',N')$.) By definition we have 
$$\align&\et'((M'-H')*t_{M'}\i(1))=\et'((M-H)'*t_M\i(0)')=\et'(((M-H)*t_M\i(0))')\\&=\et((M-H)*t_M\i(0))
=\et(H*t_M\i(1))\endalign$$
so that 46.14(c) holds in our case. For the last equality we note that 
$$\align&\et((M-H)*t_M\i(0))+\et(H*t_M\i(1))=\et((M-H)*t_M\i(0)*H*t_M\i(1))\\&
=\et(M*M)=\et(\em)=0.\endalign$$

\subhead 46.22\endsubhead
We now assume that $(M,N)\in X_n^m$ does not satisfy the assumptions of 46.18, 46.19, 46.20 or 46.21. Then 
$|M|\ge6$ and there exist $r\ge0$, $s\ge3$ such that 

$N=\{0,1,\do,r-1\}$, $M=\{r,r+1,r+2,\do,r+2s-1\}$.
\nl
Note that $(M,N)$ has the same image in $\bX_n$ as $(M',N')=(\{0,1,2,\do,2s-1\},\em)$. Since the statements of 
46.14(c) for $(M,N)$ and $(M',N')$ are equivalent, it is enough to prove 46.14(c) for $(M',N')$ instead of 
$(M,N)$. Thus we may assume that $(M,N)=(\{0,1,2,\do,2s-1\},\em)$ with $s\ge3$. We have $(M,N)\in X^s_{s^2}$.

If $\Ph$ is an admissible arrangement of $M$ let $\cc_\Ph$ be the set of all subsets $E$ of $M$ with the following
property: if $(x,y)$ is a pair in $\Ph$ then $x\in E$ if and only if $y\in E$. Note that $\cc_\Ph$ is a subspace
of the vector space $V_M$ of dimension $s$ and containing $M$. Clearly, $\Ps\m(\Ps^0\cup(\Ph-\Ps)^1))^\sh$ is a 
bijection between the sets of subsets of $\Ph$ and $\cc_\Ph$. Via this bijection the function 
$\Ps\m|\hPh\cap\Ps|\mod2$ (for $\hPh\sub\Ph$ that $|\hPh|$ is odd) can be viewed as a linear function 
$\cc_\Ph@>>>\FF_2$. This gives a bijection between $\{\hPh;\hPh\sub\Ph,|\hPh|=\text{odd}\}$ and the set of linear
functions $\cc_\Ph@>>>\FF_2$ which take the value $1$ on $M$. Using the notation $\la E\ra$ instead of $[[S,T]]$ 
where $(S,T)\in\z\i(M,N)$ and $E=S^\sh\in\tcv_M$ we see that the elements $c(M,N,\Ph,\hPh)$ (see 46.13) are the 
same as the elements
$$c(M,N,\Ph;\x)=\fra{1}{2}\sum_{E\in\cc_\Ph}(-1)^{\x(E)}\ph_{\la E\ra}\in\car(\tWW)$$
for various linear functions $\x:\cc_\Ph@>>>\FF_2$ such that $\x(M)=1$.

Now let $\Ph'$ be another admissible arrangement of $M$ and let $\x':\cc_{\Ph'}@>>>\FF_2$ be a linear form such 
that $\x'(M)=1$. We have
$$\align&(R_{c(M,N,\Ph;\x)}:R_{c(M,N,\Ph';\x')})=
\fra{1}{4}\sum_{E\in\cc_\Ph,E'\in\cc_{\Ph'}}(-1)^{\x(E)+\x'(E')}(R_{\la E\ra}:R_{\la E'\ra})\\&
=\fra{1}{4}\sum_{E\in\cc_\Ph\cap\cc_{\Ph'}}(-1)^{\x(E)+\x'(E)}
-\fra{1}{4}\sum_{E\in\cc_\Ph\cap\cc_{\Ph'}}(-1)^{\x(E)+\x'(M-E)}\\&
=\fra{1}{2}\sum_{E\in\cc_\Ph\cap\cc_{\Ph'}}(-1)^{\x(E)+\x'(E)}
=|\{\et\in\Hom(V_M,\FF_2);\et|_{\cc_\Ph}=\x,\et|_{\cc_{\Ph'}}=\x'\}|.\endalign$$
Now let $k\in[0,2s-2]$ and let $M'=\{0,1,2,\do,k-1,k+1,\do,2s-2\}$, $N'=\{k\}$. We have 
$\sum_{x\in M'}x+\sum_{x\in N'}x=h+s^2-s$ where $h=s^2-(2s-k-1)$. Since $s\ge3$ and $k\in[0,2s-2]$, we have 
$h\in[4,s^2-1]$ and $(M',N')\in\bX^s_h$.

Consider the two-sided cell $\boc'=\boc_{M',N'}\T\boc_0$ of $\WW_{I^h}$ (see 46.16) where $\boc_0$ is the 
two-sided cell associated to the sign representation $\sgn_h$ of $\fS_{n-h}$. We have $\boc'\sub\boc$ where 
$\boc=\boc_{M,N}$. 

Define an imbedding $j:M'@>>>M$ by $j(x)=x$ if $x\in[0,k-1]$, $j(x)=x+1$ if $x\in[k+1,2s-2]$. Let 
$V_M^0=\{E\in V_M;|E\cap\{k,k+1\}|=\text{even}\}$, a hyperplane in $V_M$. If $E\in V_M^0$ then $j\i(E)\in V_{M'}$.

Let $\et_1,\et_2$ be two elements of $V'_M$ such that

(a) $\et_1(E)+\et_2(E)=|E\cap\{k,k+1\}|\mod2$ for all $E\in V_M$ and $\et_1(\{k,k+1\})=\et_2(\{k,k+1\})=0$.
\nl
We define a linear function $\et':V_{M'}@>>>\FF_2$ by $\et'(E')=\et_1(j(E'))=\et_2(j(E'))$ for $E'\in V_{M'}$. 
(The last equality follows from (a) and the fact that $j(E')\cap\{k,k+1\}=\em$.) We have $\et'(M')=1$. (We use 
that 

$1=\et_1(M)=\et_1(j(M')*\{k,k+1\})=\et_1(j(M'))$
\nl
which follows from (a).) Thus we have $\et'\in V'_{M'}$. Let
$A_{\et'}$ be the object of $(\hat{\bD^h})^{un}_{M',N'}$ associated to $\et'$ by the induction hypothesis applied
to $(M',N')$. Then $R_{\sgn_h}\bxt A_{\et'}\in(\hD^{un})_{\boc'}$ is defined. We set 
$\a_{\et_1,\et_2}=\tind_{D^h}^D(R_{\sgn_h}\bxt A_{\et'})$ (see 44.20). By definition, this is an element of 
$\ck^{un}(D)$ which is an $\NN$-linear combination of objects in $D^{un}_{M,N}$. Now let $(S,T)\in\z\i(M,N)$. 
Using 44.20(h) we see that $(\a_{\et_1,\et_2}:R_{[[S,T]]})$ is $0$ if $|S\cap\{k,k+1\}|\ne1$, while if
$|S\cap\{k,k+1\}|=1$, it is 

(b) $(A_{\et'}:R_{[[S',T']]})$
\nl
where $(S',T')\in\z\i(M',N')$ is given by 

$S'=\{x<k;x\in S\}\sqc\{k\}\sqc\{x>k;x+1\in S\}$, 

$T'=\{x<k;x\in T\}\sqc\{k\}\sqc\{x>k;x+1\in T\}$. 
\nl
By the induction hypothesis, the expression (b) is equal to 
$$2^{-|M'|/2+1}(-1)^{\et'(t_{M'}\i(1)*(S'-\{k\}))}=2^{-s+2}(-1)^{\et_1(S^\sh)}=2^{-s+2}(-1)^{\et_2(S^\sh)}.$$
Hence if $\Ph$ is an admissible arrangement of $M$ and $\x:\cc_\Ph@>>>\FF_2$ is a linear function such that 
$\x(M)=1$ then
$$\align&(\a_{\et_1,\et_2}:R_{c(M,N,\Ph;\x)})=
\fra{1}{2}\sum_{E\in\cc_\Ph}(-1)^{\x(E)}(\a_{\et_1,\et_2}:R_{\la E\ra})\\&
=\fra{1}{2}\sum\Sb E\in\cc_\Ph;\\|E\cap\{k,k+1\}|=\text{even}\eSb(-1)^{\x(E)}2^{-s+2}(-1)^{\et_1(E)}\\&
=\sum\Sb E\in\cc_\Ph;\\|E\cap\{k,k+1\}|=\text{even}\eSb 2^{-s+1}(-1)^{\et_1(E)+\x(E)}.\endalign$$
This is equal to the number of elements in $\{\et_1,\et_2\}$ whose restriction to $\cc_\Ph$ is equal to $\x$.
(It is $2,1$ or $0$.)
We now apply \cite{\ORA, 9.2} to $Y=V_M$ with its basis 

$\{\{0,1\},\{1,2\},\do,\{2s-2,2s-1\}\}$
\nl
and to the family
of elements $R_{c(M,N,\Ph,\x)}$ for various $\Ph,\x$ as above and the family of elements $\a_{\et_1,\et_2}$ for 
various $\et_1,\et_2,k$ as above. (These elements are $\NN$-linear combinations of objects in $\hD^{un}_{M,N}$.) 
We see that there exists a bijection $V'_M\lra\hD^{un}_{M,N}$, $\et\lra A_\et$ such that for any $\et\in V'_M$ we
have $R_{c(M,N,\Ph,\x)}=\sum_{\et\in V'_M;\et|_{\cc_\Ph}=\x}A_\et$ for any $\Ph,\x$ as above and
$\a_{\et_1,\et_2}=A_{\et_1}+A_{\et_2}$ for any $\et_1,\et_2,k$ as above.

Now let $E\in\tcv_M$. We can rephrase 46.13(b) as follows: there exists an admissible arrangement $\Ph$ of $M$ 
such that $E\in\cc_\Ph$; moreover, 
$$\ph_{\la E\ra}=2^{-s+1}\sum_{\x\in\Hom(\cc_\Ph,\FF_2);\x(M)=1}(-1)^{\x(E)}c(M,N,\Ph;\x).$$
For $\et\in V'_M$ we then have
$$(A_\et:R_{\la E\ra})=2^{-s+1}\sum\Sb \x\in\Hom(\cc_\Ph,\FF_2);\\\x(M)=1\eSb
(-1)^{\x(E)}(A_\et:R_{c(M,N,\Ph;\x)})=2^{-s+1}(-1)^{\et(E)}.$$
We see that 46.14(c) holds in our case. This completes the proof of 46.14(c).

\subhead 46.23\endsubhead
In this subsection we assume that $p=2$. Let $P$ be a proper parabolic subgroup of $G^0$ such that $N_DP\ne\em$
and such that (setting $G'=N_GP/U_P,D'=N_DP/U_P$) we have $\hD'{}^{unc}\ne\em$. Let $\bD',\bG'$ be the quotient of
$D',G'$ by the translation action of $\cz_{G'{}^0}^0$. Let $\p:D'@>>>\bD'$ be the obvious map. From the induction 
hypothesis we see that $P/U_P$ is of type $D_r$ (with $r$ an odd square $\ge9$) or a torus, that $\hD'{}^{unc}$ 
has exactly one object $A$ up to isomorphism and that $\supp(A)$ is contained in the inverse image under $\p$ of 
the variety of unipotent elements of $\bG'$ contained in $\bD'$. Let $\hD^{un,P}$ be the subcategory of $\hD^{un}$
consisting of objects which are isomorphic to direct summands of $\ind_{D'}^D(A)$. From 27.2 and 11.9 we see that
the set of isomorphism classes in $\hD^{un,P}$ is in bijection with the set of isomorphism classes of simple 
modules of $\QQ[W_{n-r}]$. Since any noncuspidal object of $\hD^{un}$ belongs to $\hD^{un,P}$ for a $P$ as above 
(unique up to $G^0$-conjugacy) we see that the number of non-cuspidal objects of $\uhD^{un}$ is equal to
$$\sum_{k>0,s\ge0,s\text{ odd, }s^2+k=n}p_2(k)\tag a$$
where $p_2(k)$ is the number of irreducible representations of $W_k$ up to isomorphism. Now let $x_n=|\uhD^{un}|$.
From 46.14(c) we see that $x_n=|X_n|$.  Since $|X_n|$ is known from \cite{\CLA} we see that 
$$x_n=|X_n|=\sum_{k\ge0,s\ge0,s\text{ odd, }s^2+k=n}p_2(k)$$
where $p_2(0)=1$. Comparing with (a) we see that the number of cuspidal objects of $\uhD^{un}$ is $1$ if $n=s^2$ 
for some odd $s\ge3$ and is $0$ otherwise. From 12.9 we see that the set of irreducible cuspidal
admissible complexes on $D$ (up to isomorphism) is empty unless $n=s^2$ for some odd $s\ge3$ in which case it has
exactly one object (whose support is necessarily contained in the unipotent variety). Since any object of 
$\uhD^{un}$ is an admissible complex on $D$ we see that 46.14(d) holds for $D$.

This completes the inductive proof of the statements 46.14(a)-(d).

\subhead 46.24\endsubhead
Let $(M,N)=(\{0,1,2,\do,2s-1\},\em)\in X^s_n$, $n=s^2$ with $s$ odd, $s\ge3$. Define a linear function 
$\et:V_M@>>>\FF_2$ by 

$\et(E)=|E\cap t_M\i(0)|\mod2=|E\cap t_M\i(1)|\mod2$. 
\nl
Since $s$ is odd we have $\et(M)=1$ hence $\et\in V'_M$. In the setup of 46.22 we show:

(a) $A_\et\in\hD^{unc}$.
\nl
For $w\in\WW$, we have (in view of 46.22 and 44.7(i)):
$$\align&(A_\et:gr_1(K^w_D))=(-1)^{\dim G}\fra{1}{2}\sum_{E\in\tcv_M}\tr(w\vp,\la E\ra)(A_\et:R_{\la E\ra})\\&
=(-1)^{\dim G}\fra{1}{2}\sum_{E\in\tcv_M}\tr(w\vp,\la E\ra)2^{-s+1}(-1)^{\et(E)}.\endalign$$
By 44.14(a), the condition that $A_\et$ is cuspidal is that $(A_\et:gr_1(K^w_D))=0$ whenever $w\in\WW$ is not 
$D$-anisotropic. Thus it is enough to show that
$$\sum_{E\in\tcv_M}\tr(w\vp,\la E\ra)(-1)^{|E\cap t_M\i(0)|}=0\tag b$$
whenever $w\in\WW=W'_n$ satisfies the condition: $w$ is not $D$-anisotropic or equivalently, the condition:
$ws_n\in W_n$ has no eigenvalue $1$ in the reflection representation of $W_n$. Note that (b) holds by 
\cite{\CS, V, (22.5.2)}. (In that reference the words: "elements of $W'$" should be replaced by: "elements of 
$W'-W$".)

\proclaim{Theorem 46.25} Assume that $p$ satisfies the following condition: if $G^0$ has a factor of type $E_8$ or
$F_4$ then $p\ne2$. Then:

(a) if $A$ is a unipotent cuspidal character sheaf on $D$ then $A$ is clean (see 44.7);

(b) if $A$ is a unipotent character sheaf on $D$ then for any $w\in\WW$, $i\in\ZZ$ such that 
$(A:H^i(\bK^w_D))\ne0$ we have $i=\dim\supp(A)\mod2$ (or equivalently $\ee^A=(-1)^{\codim(\supp(A))}$). 
\endproclaim
By the results in \S45 we are reduced to the case where $G^0$ is simple and $\cz_G=\{1\}$. If $D=G^0$, (a) is a 
special case of 46.1(b); the fact that (a) implies (b) is proved in this case as in \cite{\CS, IV,V}. If 
$D\ne G^0$ then (a) and (b) follow from 46.4(a),(b); 46.7(a),(c); 46.8(a),(d); 46.14(a),(b). This completes the 
proof.

\subhead 46.26\endsubhead
Let $e$ be a pinning (see 1.6) of $G^0$ which projects to $(B^*,T)$ (see 28.5) under the map $p$ in 1.6. We can 
find $d\in D$ such that $\b:=\Ad(d):G^0@>>>G^0$ preserves $e$. Moreover $\b$ depends only on $D$ (not on $d$). 
Note that $\b$ has finite order, say $r$.

Let $\Bbb G$ be a connected reductive algebraic group over $\CC$ with a fixed Borel subgroup $\Bbb B$, a fixed
maximal torus $\Bbb T\sub\Bbb B$ and a fixed pinning $\ue$ which projects to $(\Bbb B,\Bbb T)$ such that $\Bbb G$
is a Langlands dual of $G^0$. In particular, $T,\Bbb T$ are Langlands dual tori. There is a unique automorphism 
$\g:\Bbb G@>>>\Bbb G$ preserving $\ue$ such that the restriction of $\g$ to $\Bbb T$ corresponds to (is 
"contragredient of") the restriction of $\b$ to $T$ under the Langlands duality between $T$ and $\Bbb T$. Note 
that $\g$ has order $r$.

A $\Bbb G$-conjugacy class $C$ in $\Bbb G$ is said to be special if some/any $g\in C$ is 
such that $g_s$ has finite order not divisible by $p$, $g_u$ is a special unipotent element of the connected 
reductive group $Z_{\Bbb G}(g_s)^0$ (see \cite{\ORA, (13.1.1)}).

Let $C$ be a special $\Bbb G$-conjugacy class in $\Bbb G$ which is $\g$-stable. For $g\in C$ let $A(g_u)$ be the 
group of components of the centralizer of $g_u$ in $Z_{\Bbb G}(g_s)^0$, let $\bA(g_u)$ be the canonical quotient 
of $A(g_u)$ defined in \cite{\ORA, p.343} (in terms of $g_u,Z_{\Bbb G}(g_s)^0$ instead of $u, G_1$) and let 
$I(g_u)$ be the kernel of the canonical homomorphism $A(g_u)@>>>\bA(g_u)$. Let 

$\ti{\Bbb A}(g)=\{(a,j)\in\Bbb G\T\ZZ/r\ZZ; a\g^j(g)a\i=g\}/Z_{\Bbb G}(g)^0$, 
\nl
a group with multiplication $(a,j)(a',j')=(a\g^j(a'),j+j')$. We identify $Z_{\Bbb G}(g)^0$ with a (normal) 
subgroup of $\ti{\Bbb A}(g)$ by $a\m(a,0)$ and we set $\Bbb A(g)=\ti{\Bbb A}(g)/Z_{\Bbb G}(g)^0$ (a finite group).
Let $\Bbb A(g)@>>>\ZZ/r\ZZ$ be the (surjective) homomorphism induced by $(a,j)\m j$. Since 
$Z_{Z_{\Bbb G}(g_s)^0}(g_u)^0=Z_{\Bbb G}(g)^0$ we see that $I(g_u)$ is naturally a subgroup of $\Bbb A(g)$. From 
the definitions we see that that in fact $I(g_u)$ is normal in $\Bbb A(g)$. Let $\cg_g=\Bbb A(g)/I(g_u)$. The
homomorphism $\Bbb A(g)@>>>\ZZ/r\ZZ$ induces a surjective a homomorphism $\cg_g@>>>\ZZ/r\ZZ$. For $j\in\ZZ/r\ZZ$ 
let $\cg_g^j$ be the inverse image of $j$ under this homomorphism. Let $\cg_C=\sqc_{g\in C}\cg_g$. Now $\Bbb G$ 
acts on $\cg_C$: if $x\in\Bbb G$, $g\in C$, then $\Ad(x)$ induces an isomorphism $\cg_g@>\si>>\cg_{xgx\i}$. Let 
$\cg_C^1=\sqc_{g\in C}\cg_g^1$, a $\Bbb G$-stable subset of $\cg_C$. For any $g\in C$, the set of $\Bbb G$-orbits
on $\cg_C^1$ is in natural bijection with the (finite) set of $\cg_g$-conjugacy classes in $\cg_g^1$. Thus 
$\Bbb G$ acts on $\cg_C^1$ with finitely many orbits. This makes $\cg_C^1$ into an algebraic variety (a finite 
union of homogeneous spaces for $\Bbb G$).

Let $\fP_\g$ be the set of all triples $(C,X,\ce)$ where $C$ is a $\g$-stable special $\Bbb G$-conjugacy class in 
$\Bbb G$, $X$ is a $\Bbb G$-orbit in $\cg_C^1$ and $\ce$ is an irreducible $\Bbb G$-equivariant local system on 
$X$ (up to isomorphism). Let $\fP_\g^{un}$ be the set of all $(C,X,\ce)\in\fP_\g$ such that $C$ is a unipotent 
$\Bbb G$-conjugacy class in $\Bbb G$.

\subhead 46.27\endsubhead
We have $\fP_\g^{un}=\sqc_{C}\fP_{\g,C}^{un}$ where $C$ runs over the set of $\g$-stable special unipotent classes
in $\Bbb G$ and $\fP_{\g,C}^{un}$ is the set of triples in $\fP_\g^{un}$ whose first component is $C$. Under the 
Springer correspondence, the set of $\g$-stable special unipotent classes in $\Bbb G$ is in bijection with the set
of special irreducible representations $E_0$ (up to isomorphism) of the Weyl group of $\Bbb G$ or of $G^0$ whose 
character is fixed by $\e:\WW@>>>\WW$ and hence in bijection (via $E_0\m\boc_{E_0}$, see 43.6) with the set of
$\e$-stable two-sided cells of $\WW$; let $C_\boc$ be the special unipotent class corresponding to the two-sided
cell $\boc$. Assume that $p$ is as in 46.25. We have the following result:

(a) {\it For any $\e$-stable two-sided cell $\boc$ in $\WW$ there is a natural bijection 
$\hD^{un}_\boc\lra\fP_{\g,C_\boc}^{un}$.}
\nl
By the results in \S45 we are reduced to the case where $G^0$ is simple and $\cz_G=\{1\}$. If $D=G^0$, (a) is 
established in \cite{\CS, IV,V}. If $D\ne G^0$ then (a) follows from 46.4(d), 46.7, 46.8, 46.14(c). 

By taking disjoint union over the various $\boc$ we obtain a bijection $\hD^{un}\lra\fP_\g^{un}$. We will show
elsewhere that this extends to a natural bijection $\hD\lra\fP_\g$. (See \cite{\CS, IV,V} for the case where 
$G=G^0$.)

\enddocument